\documentclass[12pt]{amsart} \usepackage{latexsym, amssymb}

\newtheoremstyle{slthm}
  {9pt}
  {5pt}
  {\slshape}
  {}
  {\bfseries}
  {.}
  {.5em}
  {\thmname{#1}\thmnumber{ #2}\thmnote{ (#3)}}

\newtheoremstyle{prcl}
  {9pt}
  {5pt}
  {\slshape}
  {}
  {\bfseries}
  {.}
  {.5em}
  {\thmname{#3}\thmnumber{ #2}}

\theoremstyle{slthm}
\newtheorem{thm}{Theorem}[section]
\newtheorem{lemma}[thm]{Lemma}
\newtheorem{prop}[thm]{Proposition}
\newtheorem{cor}[thm]{Corollary}

\theoremstyle{definition}

\newtheorem{df}[thm]{Definition}
\newtheorem{dfs}[thm]{Definitions}
\newtheorem{rmkdf}[thm]{Remark and Definition}
\newtheorem{nrmk}[thm]{Remark}
\newtheorem{nrmks}[thm]{Remarks}
\newtheorem{expl}[thm]{Example}

\theoremstyle{remark}
\newtheorem*{rmk}{Remark}
\newtheorem*{rmks}{Remarks}

\theoremstyle{prcl}
\newtheorem*{prclaim}{Proclaim}

\newenvironment{renumerate}
        {\renewcommand{\labelenumi}{\rm (\roman{enumi})}
         \begin{enumerate}}
        {\end{enumerate}}

\newenvironment{lenumerate}[2]
        {\renewcommand{\labelenumi}{\rm (#1\theenumi)}
         \begin{enumerate}{\setcounter{enumi}{#2}}}
        {\end{enumerate}}

\newcounter{flexnummark}

\DeclareMathOperator{\cl}{cl}
\DeclareMathOperator{\fr}{fr}
\DeclareMathOperator{\bd}{bd}
\DeclareMathOperator{\ir}{int}

\DeclareMathOperator\gr{gr}
\DeclareMathOperator{\Ma}{M}
\DeclareMathOperator\dist{dist}
\DeclareMathOperator\lnest{\la_{\text{nest}}}

\newcommand{\rest}[1]{\!\!\upharpoonright_{#1}}

\newcommand{\into}{\longrightarrow}

\newcommand{\set}[1]{\left\{#1\right\}}

\newcommand{\NN}{\mathbb{N}}

\newcommand{\QQ}{\mathbb{Q}}
\newcommand{\RR}{\mathbb{R}}

\newcommand{\curly}[1]{\mathcal{#1}}
\newcommand{\A}{\curly{A}}
\newcommand{\B}{\curly{B}}
\newcommand{\C}{\curly{C}}
\newcommand{\D}{\curly{D}}
\newcommand{\E}{\curly{E}}

\newcommand{\G}{\curly{G}}

\newcommand{\K}{\curly{K}}

\newcommand{\N}{\curly{N}}
\renewcommand{\P}{\curly{P}}
\newcommand{\Q}{\curly{Q}}
\newcommand{\R}{\curly{R}}
\renewcommand{\S}{\curly{S}}

\newcommand{\V}{\curly{V}}

\newcommand{\la}{\curly{L}}

\numberwithin{equation}{section}

\allowdisplaybreaks[2]

\title {The theorem of the complement for nested subpfaffian sets}

\author {J.-M. Lion and P. Speissegger}

\address {UFR math\'ematiques, Universit\'e de Rennes I, Campus de
  Beaulieu, 35042 Rennes cedex, France}

\email{jean-marie.lion@univ-rennes1.fr}

\address {Department of Mathematics \& Statistics, McMaster
  University, 1280 Main Street West, Hamilton, Ontario L8S 4K1,
  Canada}

\email {speisseg@math.mcmaster.ca}

\date{\today}

\subjclass {Primary 14P10, 58A17; Secondary 03C99}

\keywords {O-minimal structures, pfaffian systems, analytic
  stratification}

\thanks {Supported by the CNRS of France and NSERC of Canada grant
  RGPIN 261961.  We thank the Newton Institute at Cambridge and the
  Fields Institute in Toronto, where parts of this paper were
  written.}

\begin{document}

\begin{abstract}
  Let $\R$ be an o-minimal expansion of the real field, and let
  $\lnest(\R)$ be the language consisting of all nested Rolle leaves
  over $\R$.  We call a set nested subpfaffian over $\R$ if it is the
  projection of a positive boolean combination of definable sets and
  nested Rolle leaves over $\R$.  Assuming that $\R$ admits analytic
  cell decomposition, we prove that the complement of a nested
  subpfaffian set over $\R$ is again a nested subpfaffian set over
  $\R$.  As a corollary, we obtain that if $\R$ admits analytic cell
  decomposition, then the pfaffian closure $\P(\R)$ of $\R$ is
  obtained by adding to $\R$ all nested Rolle leaves over $\R$, a
  one-stage process, and that $\P(\R)$ is model complete in the
  language $\lnest(\R)$.
\end{abstract}

\maketitle
\markboth{}{}

\section*{Introduction}

The basic objects we study in this paper are nested pfaffian sets over
a given o-minimal expansion of the real field.  Before defining them,
let us briefly recall some of the history around the notion of
pfaffian functions: roughly speaking, pfaffian functions are maximal
solutions of triangular systems of partial differential equations with
polynomial coefficients, see Khovanskii \cite{kho:fewnomials},
Gabrielov \cite{gab:pfaffmult} and Wilkie \cite{wil:exp}.  In his
thesis \cite{kho:fewnomials}, Khovanskii proves (among other things)
that any set defined by finitely many equations and inequalities
between pfaffian functions has a finite number of connected
components.  In the early 1980s, Van den Dries conjectured that the
expansion of the real field by all pfaffian functions was model
complete, which, together with Khovanskii's theorem, would imply that
this expansion is o-minimal.  (For generalities on o-minimal
structures, we refer the reader to Van den Dries and Miller
\cite{vdd-mil:cat}.)

Wilkie \cite{wil:complement} used a different approach to obtain 
the first general o-minimality result for pfaffian functions, showing
that the real field expanded by all \textit{totally defined} pfaffian
functions is o-minimal.  Based on Lion and Rolin
\cite{lio-rol:feuilletages}, this theorem was strengthened in the
following way: given an o-minimal expansion $\R$ of the real field, we
call a function pfaffian over $\R$ if it is a maximal solution of a
triangular system of partial differential equations with coefficients
definable in $\R$.  Then \cite{spe:pfaffian} there is an o-minimal
expansion $\P(\R)$ of $\R$, called the pfaffian closure of $\R$, such
that every pfaffian function over $\P(\R)$ is definable in $\P(\R)$.

However, none of the above o-minimality proofs establish the model
completeness of the respective structures in any meaningfull language,
and Gabrielov restated the question in \cite{gab:invent}.  Based on
techniques used in \cite{lio-spe:an_strat}, we give here a natural
language $\lnest(\R)$ in which $\P(\R)$ is model complete, in the case
where $\R$ admits analytic cell decomposition.  To do this, we draw
inspiration from the setting in \cite{kho:fewnomials}, where pfaffian
functions are replaced by nested separating integral manifolds.  We
need a few definitions to state the precise theorem.

We denote by $G^l_n$ the \textbf{Grassmannian} of all $l$-dimensional
linear subspaces of $\RR^n$.  Let $M \subseteq \RR^n$ be a
$C^1$-submanifold of dimension $m$.  Throughout this paper, we
identify the tangent space $T_xM$ of $M$ at $x \in M$ with an element
of $G^m_n$ in the usual way, and we define an
\textbf{$l$-distribution} on $M$ to be a map $d:M \into G^{l}_n$ such
that $d(x) \subseteq T_xM$ for all $x \in M$.  For example, the
$m$-distribution $g_M$ on $M$ defined by $g_M(x):= T_xM$ is the
\textbf{Gauss map of $M$}.

Let $d$ be an $l$-distribution on $M$.  An immersed manifold $V
\subseteq M$ of dimension $l$ is called an \textbf{integral
  manifold} of $d$ if $T_xV = d(x)$ for all $x \in V$.  A
\textbf{leaf} of $d$ is a maximal connected integral manifold of $d$.

Assume now that $l = m-1$.  Then a leaf $V$ of $d$ is \textbf{Rolle}
(see Moussu and Roche \cite{mou-roc:khovtheory}) if $V$ is a closed
submanifold of $M$ and for every $C^1$-curve $\gamma:[0,1] \into M$
such that $\gamma(0), \gamma(1) \in V$, there exists a $t \in [0,1]$
such that $\gamma'(t) \in d(\gamma(t))$.  The following criterion for
the Rolle property is crucial to our paper:

\begin{prclaim}[Haefliger's Theorem \cite{hae:theorem,mou-roc:finitude}] 
  Assume that $M$ and $d$ as above are analytic and that $M$ is simply
  connected and $d$ is integrable (Definition \ref{integrable_dist}).
  Then every leaf of $d$ is a Rolle leaf.
\end{prclaim}

Let $d = (d_0, \dots, d_k)$ be a tuple of distributions on $M$. We
call $d$ \textbf{nested} if each $d_j$ is an $(m-j)$-distribution on
$M$ and $d_k(x) \subseteq d_{k-1}(x) \subseteq \cdots \subseteq d_0(x)
= g_M(x)$ for all $x \in M$.

Assume that $d$ is nested, and let $V = (V_0, \dots, V_k)$ be a tuple
of immersed manifolds contained in $M$.  We call $V$ a \textbf{nested
  integral manifold} of $d$ if each $V_j$ is an integral manifold of
$d_j$ and $V_0 \supseteq \cdots \supseteq V_k$.  Moreover, $V$ is a
\textbf{nested leaf} (respectively, \textbf{nested Rolle leaf}) of
$d$, if $V_0 = M$ and for $j=1, \dots, k$, the set $V_j$ is a leaf
(respectively, Rolle leaf) of the restriction $d_j \rest{V_{j-1}}$
of $d_j$ to $V_{j-1}$.  Note that in this situation, $d_j
\rest{V_{j-1}}$ is of class $C^1$ and $d_j(x) \subseteq T_x V_{j-1}$
is of codimension $1$ for all $x \in V_{j-1}$; in particular, $\dim
V_j = m-j$.

\begin{prclaim}[Example 1 \cite{spe:pfaffian}] \rm Let $\Omega =
  (\omega_1, \dots, \omega_k)$ be a family of differential $1$-forms
  on $M$, and assume that $\Omega$ is \textbf{nonsingular}, that
  is, $\omega_1(x) \wedge \cdots \wedge \omega_k(x) \ne 0$ for all $x
  \in M$.  For $j = 1, \dots, k$, put $d_j(x) := \ker \omega_1(x)
  \cap \cdots \cap \ker\omega_j(x)$; then $d:= (g_M, d_1, \dots, d_k)$
  is a nested distribution on $M$.

  Conversely, let $d = (d_0, \dots, d_k)$ be a nested distribution on
  $M$ and assume that $M$ is simply connected.  Define unit vector
  fields $a_j = (a_{j1}, \dots, a_{jn})$ on $M$, for $j=1, \dots, k$,
  by induction on $j$ as follows: let $a_1$ be one of the two unit
  vector fields orthogonal to $d_1$, and for $j>1$ let $a_j$ be one of
  the two unit vector fields orthogonal to the vector space field
  spanned by $d_j \cup \{a_1, \dots, a_{j-1}\}$.  Finally, put
  $\omega_j:= a_{j1} dx_1 + \cdots + a_{jn} dx_n$ for $j=1, \dots, k$.
  Then $\Omega:= (\omega_1, \dots, \omega_k)$ is a nonsingular family
  of differential $1$-forms on $M$.

  In the notation used before Example 1, if $M$ is simply connected,
  then $V$ is a nested integral manifold (leaf, Rolle leaf) of this
  $d$ if and only if $V_1$ is an integral manifold (leaf, Rolle leaf)
  of $\omega_1$---as defined in \cite{spe:pfaffian}---and for $j = 2,
  \dots, k$, $V_j$ is an integral manifold (leaf, Rolle leaf) of the
  pull-back of $\omega_j$ on $V_{j-1}$.
\end{prclaim}

Let $\R$ be an o-minimal expansion of the real field, and assume that
$M$ is definable in $\R$ and of class $C^2$ (the reason for the latter
assumption is to be consistent with pull-backs, see the conventions
below).  Let $d = (d_0, \dots, d_k)$ be a nested distribution on $M$
and $V = (V_0, \dots, V_k)$ be a nested integral manifold of $d$.  For
$l \le n$ we identify $G^l_n$ with an algebraic (and hence definable
in $\R$) subvariety of $\,\RR^{n^2}$ (see the conventions below for
details).  We call $d$ \textbf{definable} if $d$ is a definable map
under this identification.  If $d$ is definable, we call $V$ a
\textbf{nested integral manifold over $\R$}, and if $V$ is a nested
Rolle leaf of $d$, we call $V$ a \textbf{nested Rolle leaf over $\R$}.
Note that in the latter situation, the leaves $V_0, \dots, V_k$ are
uniquely determined by $d$ and $V_k$, but that $V_k$ is not definable
in $\R$ in general.  For convenience, we call a set $W \subseteq
\RR^n$ an \textbf{integral manifold over $\R$} if there is a nested
integral manifold $(W_0, \dots, W_k)$ over $\R$ with $W_k = W$, and if
in addition $(W_0, \dots, W_k)$ is a nested Rolle leaf over $\R$, we
call $W$ a \textbf{Rolle leaf over $\R$}.

\begin{prclaim}[Example 2] 
  \rm Let $C \subseteq \RR^n$ be a definable $C^2$-cell.  Taking $M =
  C$ and $d$ to be the Gauss map on $M$ makes $C$ trivially into a
  Rolle leaf over $\R$.
\end{prclaim}

A set $X \subseteq \RR^n$ is a \textbf{basic nested pfaffian set over
  $\R$}, if there are a definable set $A \subseteq \RR^n$ and a Rolle
leaf $W \subseteq \RR^n$ over $\R$ such that $X = A \cap W$.  A
\textbf{nested pfaffian set over $\R$} is a finite union of basic
nested pfaffian sets over $\R$, and a \textbf{nested subpfaffian set
  over $\R$} is the image under a coordinate projection of a nested
pfaffian set over $\R$.

We let $\lnest(\R)$ be the collection of all Rolle leaves over $\R$,
and we denote by $\N(\R)$ the expansion of $\R$ by all $W \in
\lnest(\R)$.  It follows from $C^2$-cell decomposition in $\R$ and
Example 2 that every set definable in $\R$ is quantifier-free
definable in $\N(\R)$.

Khovanskii theory as in \cite{mou-roc:khovtheory,spe:pfaffian}
generalizes in a straightforward way to the setting of Rolle
leaves over $\R$ (Sections \ref{nested} and \ref{khovsection}).  It
follows in particular that every nested subpfaffian set over $\R$ has
finitely many connected components, and every  Rolle leaf over
$\,\R$ is definable in $\P(\R)$ (Proposition \ref{definability}).
Hence, every nested subpfaffian set $X$ over $\R$ is definable in the
\hbox{o-minimal} structure $\P(\R)$; we denote its dimension by $\dim
X$.  Building on these observations, we prove:

\begin{prclaim}[Main Theorem]
  Assume that $\,\R$ admits analytic cell decomposition.  Then the
  complement of every nested subpfaffian set over $\,\R$ is again
  nested subpfaffian over $\,\R$; in particular, $\N(\R)$ is model
  complete.
\end{prclaim}

The Main Theorem implies that, in the construction of $\P(\R)$ in
\cite{spe:pfaffian}, every Rolle leaf added to $\,\R$ is a nested
subpfaffian set over $\R$ (Proposition \ref{interdefinable}). 

\begin{prclaim}[Corollary 1]
  If $\,\R$ admits analytic cell decomposition, then $\N(\R)$ is
  interdefinable with $\P(\R)$; in particular, $\P(\R)$ is model
  complete in $\lnest(\R)$. \qed
\end{prclaim}

The model completeness of $\P(\R)$ in the language $\lnest(\R)$
remains an open problem if $\R$ does not admit analytic cell
decomposition.  Also, even in the analytic case, we do not know
whether the reduct of $\P(\R)$ generated by all pfaffian functions
over $\R$ is model complete.  

To prove the Main Theorem, we use Corollary 2.9 of
\cite{vdd-spe:genpower}, with $\Lambda$ there equal to the collection
of all nested pfaffian sets over $\R$ contained in $[-1,1]^n$, for $n
\in \NN$.  This means that it suffices to establish Axioms (I)--(IV)
there; they follow easily from the following list of statements, which
correspond to theorems proved in this paper:
\begin{itemize}
\item[(P1)] every set definable in $\R$ is nested pfaffian over $\R$;
\item[(P2)] the union, the intersection and the cartesian product of
  two nested pfaffian sets over $\R$ are nested pfaffian over $\R$,
  and each connected component of a nested pfaffian set over $\R$ is
  nested pfaffian over $\R$;
\item[(P3)] if $X \subseteq \RR^n$ is nested pfaffian over $\R$ and $1
  \le m \le n$, there is a finite collection $\P$ of analytic
  manifolds contained in $X$ such that $\Pi_m(X) = \bigcup_{Y \in \P}
  \Pi_m(Y)$ and for each $Y \in \P$, the set $Y$ is nested pfaffian
  over $\R$, $\dim Y \le m$ and there is a strictly increasing
  $\lambda:\{1, \dots, \dim Y\} \into \{1, \dots, m\}$ such that
  $\Pi_\lambda \rest{Y}:Y \into \RR^{\dim Y}$ is an immersion;
\item[(P4)] if $X \subseteq \RR^n$ is nonempty and nested pfaffian
  over $\R$, there is a closed, nested subpfaffian set $Y \subseteq
  \RR^n$ over $\R$ such that $\fr X \subseteq Y$ and $\dim Y < \dim
  X$.
\end{itemize}
Here, for $l \ge m$, the map $\Pi_m:\RR^{l} \into \RR^m$ denotes the
projection on the first $m$ coordinates, and for every strictly
increasing $\lambda:\{1, \dots, m\} \into \{1, \dots, l\}$, the map
$\Pi_\lambda:\RR^{l} \into \RR^{m}$ denotes the projection
$\Pi_\lambda(x_1, \dots, x_{l}):= (x_{\lambda(1)}, \dots,
x_{\lambda(m)})$.  Also, for any set $S\subseteq \RR^n$, we denote by
$\cl S$ the topological closure of $S$, and we define the
frontier of $S$ as the set $\fr S:= \cl S \setminus S$.

Statement (P1) follows from Example 2, and (P2) follows from
Khovanskii theory for nested pfaffian sets over $\R$ (Corollary
\ref{pfaff_intersection}).  Statement (P3) follows from the Fiber
Cutting Lemma for nested pfaffian sets over $\R$ (Corollary
\ref{fiber_cutting}), which is obtained using an approach similar to
Gabrielov's in \cite{gab:frontier}.  The main contribution of this
paper is to establish (P4).  To explain how this is done, we let $W
\subseteq \RR^n$ be a Rolle leaf over $\R$.  By Khovanskii theory
again (Corollary \ref{khovprop}), statement (P4) follows from

\begin{prclaim}[Theorem 1]
  Assume that $\R$ admits analytic cell decomposition.  Then there is
  a closed, nested subpfaffian set $Y \subseteq \RR^n$ over $\R$ such
  that $\fr W \subseteq Y$ and $\dim Y < \dim W$.
\end{prclaim}

Theorem 1 was proved in the special case $\dim W = n-1$ by Cano et
al. in \cite{can-lio-mou:hyper}.  For the proof of the general case
(Section \ref{main}), we let $d = (d_0, \dots, d_k)$ be a definable,
nested distribution on a bounded, definable $C^2$-manifold $M
\subseteq \RR^n$ and $V = (V_0, \dots, V_k)$ be a nested Rolle leaf of
$d$ such that $W = V_k$.  We consider $\fr W$ as a Hausdorff limit
of a certain sequence of integral manifolds of a definable nested
distribution $d'$ on $M$ derived from $d$ (Section
\ref{limits_section}).  We then use the method of blowing up along
$d'$ (Section \ref{jet}), similar to \cite{lio-spe:an_strat}, to
recover---roughly speaking---distributions on the frontier of $M$,
such that $\fr W$ is almost everywhere an integral manifold of one of
these distributions.  (Strictly speaking, these distributions are
recovered on the frontier of the manifold obtained from $M$ by blowing
up and have the described property for the corresponding lifting of
$W$; to keep notations simple, we continue using $M$ and $W$ in the
introduction.)  The main problems solved in this paper are the
following: we did not know in \cite{lio-spe:an_strat} if
\begin{itemize}
\item[(a)] the distributions recovered in this way were components of
  definable \textit{nested} distributions;
\item[(b)] the integral manifolds in question were contained in Rolle
  leaves over $\R$ of the same dimension.
\end{itemize}
Here we deal with (a) and (b) separately; we establish (a) for the
case that $\R$ is any o-minimal expansion of the real field, but we
need to assume that $\R$ admits analytic cell decomposition to
establish (b).

For (a), we define the degree of $d$ to be the number of component
distributions of $d$ whose associated foliation of $M$ is \textit{not}
definable in $\R$ (Section \ref{nested}).  We show in Section
\ref{limits_section} that this degree behaves well in the following
sense: the nested distribution $d'$ derived from $d$ used to describe
$\fr W$ as a Hausdorff limit, as mentioned above, has degree less than
or equal to that of $d$.  Moreover, we also prove that the negligible
set, off which $\fr W$ is a finite union of integral manifolds of the
recovered distributions, is a union of Hausdorff limits of the same
type obtained from distributions of degree less than or equal to that
of $d$.  These observations and a refinement of the blowing-up method
in \cite{lio-spe:an_strat} yield the following version of (a),
combining Propositions \ref{frontier_reduction} and \ref{recover_1}
below:

\begin{prclaim}[Proposition 1]
  There is a $q \in \NN$ and, for $p=1, \dots, q$, there are $n_p \ge
  n$ and an integral manifold $U_p \subseteq \RR^{n_p}$ over $\R$ such
  that $\fr W \subseteq \Pi_n(U_1) \cup \cdots \cup \Pi_n(U_q)$ and
  for each $p$, the set $U_p$ is definable in $\P(\R)$ and $\dim
  \Pi_n(U_p) < \dim W$.
\end{prclaim}

We call an integral manifold $U \subseteq \RR^n$ over $\R$
\textbf{definable in $\P(\R)$} if there is a nested integral manifold
$Z = (Z_0, \dots, Z_l)$ over $\R$ such that $U = Z_l$ and each $Z_j$
is definable in $\P(\R)$.  For (b) it now remains to show, for each of
the integral manifolds $U_p$ of Proposition 1, that $\Pi_n(U_p)$ is in
turn contained in a subpfaffian set over $\R$ of dimension at most
$\dim \Pi_n(U_p)$ (implied by Proposition \ref{Rolle}).  

To do so, we let $d = (d_0, \dots, d_k)$ be a definable nested
distribution on some manifold $M$ and $V = (V_0, \dots, V_k)$ be a
nested integral manifold of $d$ definable in $\P(\R)$, and we try to
reduce to a situation where, up to finite union and projection, the
leaf $L_k$ of $d_k$ containing $V_k$ is a Rolle leaf $d$.  To
establish the Rolle property of $L_k$, we want to use Haefliger's
Theorem; this is one of the reasons for our assumption that $\R$
admits analytic cell decomposition.  Thus, if $k=1$, we can easily
recover the Rolle property from Haefliger's Theorem using analytic
cell decomposition of $M$.  If $k>1$, however, we can only apply
Haefliger's Theorem if $V_{k-1}$ is \textit{simply connected}.
Proceeding by induction on $k$, we may assume that $(V_0, \dots,
V_{k-1})$ is a nested Rolle leaf over $\R$; therefore, what we need to
establish is (see Corollary \ref{pfaffian_complement}):

\begin{prclaim}[Proposition 2]
  Assume that $\R$ admits analytic cell decomposition.  Then $V_{k-1}$
  is a finite union of simply connected nested subpfaffian sets over
  $\R$ that are analytic manifolds.
\end{prclaim}

To prove this, we introduce in Section \ref{proper} the notion of
\textit{proper} nested subpfaffian set over $\R$.  These are certain
projections of nested pfaffian sets $X \subseteq [-1,1]^n$ over $\R$
that are \textit{restricted off} some closed set $Z$; if $Z = \{0\}$,
this means that for every $r>0$, the set $X \setminus (-r,r)^n$ is a
restricted nested pfaffian set similar to Gabrielov's in
\cite{gab:frontier} or \cite{wil:exp}.  Remarkably, based on the ideas
in \cite{gab:frontier}---adapted to our situation in Sections
\ref{analytic} and \ref{proper}---we obtain a cell decomposition
theorem for certain proper nested subpfaffian sets over $\R$
(Proposition \ref{proper_cells}).  Proposition 2 then follows from the
observation that, up to an analytic inversion of the ambient space,
$V_{k-1}$ is a restricted pfaffian set off $\{0\}$ (Example
\ref{proper_example} and Proposition \ref{inside_out}).

Unfortunately, Haefliger's Theorem is false without the analyticity
assumption, even in the o-minimal context, see for instance
\cite[Section 3]{lio-spe:haefliger}.  We prove there a weaker version
of Haefliger's Theorem in the general o-minimal context, but we do not
know if the proof of the Main Theorem goes through with this weaker
version of Haefliger's Theorem.

\begin{prclaim}[Conventions] 
  \rm Throughout this paper, all cells, manifolds, functions, maps,
  etc. are of class $C^1$, and manifolds are embedded, unless
  otherwise specified.  We write $\NN = \{0,1,2,3,\dots\}$ for the
  set of all natural numbers.  We sometimes abbreviate ``analytic'' as
  ``$C^\omega$'', and we extend the usual linear ordering on $\NN \cup
  \{\infty\}$ to $\NN \cup \{\infty,\omega\}$ by putting $\omega >
  \infty$.  We shall use ``component'' in place of ``connected
  component'' whenever the meaning is clear from context.

  $\R$ denotes a fixed, but arbitrary, o-minimal expansion of the real
  field, and ``definable'' means ``definable in $\R$ with parameters
  from $\RR$'' unless indicated otherwise.

  A \textbf{box} in $\RR^n$ is a subset of the form $I_1 \times \cdots
  \times I_n$, where each $I_j$ is a nonempty open interval in $\RR$.
  For $x \in \RR^n$, we put $|x|:= \sup\{|x_1|, \dots, |x_n|\}$, and
  for $r>0$, we let $B(x,r) := \set{y \in \RR^n:\ |y-x| < r}$.

  For any set $S\subseteq \RR^n$, we denote by $|S|$ the cardinality
  of $S$, by $\cl S$ the topological closure of $S$ and by $\ir S$
  the interior of $S$, and we define the \textbf{boundary} of $S$ as
  $\bd S := \cl S \setminus \ir S$ and the \textbf{frontier} of $S$
  as the set $\fr S:= \cl S \setminus S$.  A family $\S$ of subsets
  of $\RR^n$ is a \textbf{stratification} if the members of $\S$ are
  pairwise disjoint and for all $S_1, S_2 \in \S$, we have either $S_1
  \cap \cl S_2 = \emptyset$ or $S_1 \subseteq \cl S_2$.  In this
  paper, we also use \textbf{Whitney stratifications}; their
  definition is more technical, and we refer the reader to Sections 1
  and 4 of Miller and Van den Dries \cite{vdd-mil:cat} for a thorough
  discussion in the o-minimal context.

  We let $\Sigma_n$ be the collection of all permutations of $\{1,
  \dots, n\}$.  For $\sigma \in \Sigma_n$, we write $\sigma:\RR^n
  \into \RR^n$ for the map defined by $\sigma(x_1, \dots, x_n):=
  \left(x_{\sigma(1)}, \dots, x_{\sigma(n)}\right)$.

  For $l \ge m$, the map $\Pi^l_m:\RR^{l} \into \RR^m$ denotes the
  projection on the first $m$ coordinates; and for every strictly
  increasing $\lambda:\{1, \dots, m\} \into \{1, \dots, l\}$, the map
  $\Pi^l_\lambda:\RR^{l} \into \RR^{m}$ denotes the projection
  $\Pi^l_\lambda(x_1, \dots, x_{l}):= (x_{\lambda(1)}, \dots,
  x_{\lambda(m)})$.  When $l$ is clear from context, we usually write
  $\Pi_m$ and $\Pi_\lambda$ in place of $\Pi^l_m$ and $\Pi^l_\lambda$,
  respectively.

  We let $\K_n$ be the space of all compact subsets of $\RR^n$
  equipped with the Hausdorff metric.  (We consider $\emptyset \in
  \K_n$ with $d(A,\emptyset) = \infty$ for all nonempty $A \in \K_n$.)
  Given a sequence $(A_\iota)_{\iota \in \NN}$ of bounded subsets of
  $\RR^n$, we say that $(A_\iota)$ \textbf{converges to} $C \in \K_n$
  if the sequence $(\cl A_\iota)$ converges in $\K_n$ to $C$, and in
  this situation we write $C = \lim_\iota A_\iota$.  We refer the
  reader to Kuratowski \cite{kur:topology} for the classical results
  about $\K_n$; in particular, we shall often use without reference
  the fact that every bounded sequence in $\K_n$ contains a convergent
  subsequence.

  Let $l \leq n$.  We denote by $G^l_n$ the \textbf{Grassmannian} of
  all $l$-dimensional vector subspaces of $\RR^n$.  This $G^l_n$ is an
  analytic, real algebraic variety with a natural analytic embedding
  into the vector space $M_n$ of all real valued $(n \times
  n)$-matrices: each $l$-dimensional vector space $E$ is identified
  with the unique matrix $A_E$ (with respect to the standard basis of
  $\RR^n$) corresponding to the orthogonal projection on the
  orthogonal complement of $E$ (see Section 3.4.2 of
  \cite{boc-cos-roy:geometrie}); in particular, $E = \ker(A_E)$.  We
  shall identify $M_n$ with $\RR^{n^2}$ via the map $A = (a_{ij})
  \mapsto z_A = (z_1, \dots, z_{n^2})$ defined by $a_{ij} =
  z_{n(i-1)+j}$, and we identify $G^l_n$ with its image in $M_n$ under
  this map.  Note that the sets $G^{0}_{n}, \dots ,G^{n}_{n}$ are the
  components of $G_{n}:= \bigcup_{p=0}^{n} G^{p}_{n}$.

  Let $M \subseteq \RR^n$ be a manifold of dimension $m \le n$.  For
  $\eta > 0$, we say that $M$ is \textbf{$\eta$-bounded} if for every
  $x \in M$ there is a matrix $L = (l_{ij}) \in \Ma_{n-m,m}(\RR)$ such
  that $|L| := \sup_{i,j} |l_{i,j}| \le \eta$ and $T_x M = \set{(u,L
    u):\ u \in \RR^m}$.

  We call a map $d:M \into G_n$ a \textbf{distribution on $M$} if
  $d(x) \subseteq T_xM$ for all $x \in M$.  Given two maps $d,e:M
  \into G_n$, we write $d \cap e:M \into G_n$ for the map defined by
  $(d \cap e)(x):= d(x) \cap e(x)$, and we write $d \subseteq e$ if
  $d(x) \subseteq e(x)$ for all $x \in M$.  Note that, by linear
  algebra, if $M$ and $d,e:M \into G_n$ are of class $C^p$, with $p
  \in \NN \cup \{\infty,\omega\}$, then so is $d \cap e$.  If $d:M
  \into G_n$ is a map, we say that \textbf{$d$ has dimension} if $d(M)
  \subseteq G^m_n$ for some $m \le n$; in this situation, we put $\dim
  d:= m$.

  Assume that $M$ is of class $C^2$, let $N \subseteq \RR^l$ be a
  $C^2$-manifold and $f:N \into M$ a $C^2$-map, and let $g$ be a
  distribution on $M$.  The \textbf{pull-back} of $g$ on $N$ by $f$ is
  the distribution $f^*g$ on $N$ defined by
  $
    f^*g(y) := (df_y)^{-1} \big(g(f(y))\big),
  $
  where $df_y:T_yN \into T_{f(y)}M$ is the linear map defined by the
  jacobian matrix of $f$ at $y$ and $(df_y)^{-1}(S)$ denotes the
  inverse image of $S$ under this map for any $S \subseteq
  T_{f(y)}M$.  If $N$ is a $C^2$-submanifold of $M$ and $f$ is the
  inclusion map, we write $g^N$ in place of $f^*g$.
\end{prclaim}

\section{Preliminaries}
\label{prelim}

This section introduces some terminology and contains several basic
lemmas needed later on.  Below, for any map $f$ we denote by $\gr f$
the graph of $f$.

\begin{lemma} 
  \label{o-min_lemma} 
  Let $p \ge 1$ be finite, and let $M \subseteq \RR^n$ be a definable
  $C^p$-manifold of dimension $d$.  Let also $m \le n$, and assume
  that $\Pi_m\rest{M}$ has constant rank $\nu$.  Then $M$ is the union
  of finitely many definable open subsets $N$ such that $\Pi_m(N)$ is
  a $C^p$-submanifold of $\,\RR^m$ of dimension $\nu$.
\end{lemma}

\begin{proof}
  Given a permutation $\sigma \in \Sigma_m$ and denoting by
  $\sigma:\RR^n \into \RR^n$ the map defined by $\sigma(x) :=
  (x_{\sigma(1)}, \dots, x_{\sigma(m)}, x_{m+1}, \dots, x_n)$, the set
  \begin{equation*}
    M_\sigma := \set{y \in M:\ \Pi_\nu\rest{\sigma(M)}
      \text{ is a submersion at } \sigma(y)} 
  \end{equation*}
  is an open subset of $M$ and $M$ is the union of all $M_\sigma$ with
  $\sigma \in \Sigma_m$.  Thus by replacing $M$ with each
  $\sigma(M_\sigma)$, we may assume that $\Pi_\nu\rest{M}$ is a
  submersion; in particular, $U:= \Pi_\nu(M)$ is open.  Since $M$ is
  definable, the hypotheses now imply that there is a $K \in \NN$ such
  that for every $y \in \RR^\nu$, the fiber $\Pi_m(M)_y$ has at most
  $K$ elements.  For $k \in \{0,\dots, K\}$ we let $D_k \subseteq
  \RR^\nu$ be the set of all $y$ such that $\Pi_m(M)_y$ has exactly
  $k$ elements.  Let $\C$ be a $C^p$-cell decomposition of $\RR^m$
  compatible with $\Pi_m(M)$ such that $\D:= \set{\Pi_\nu(C):\ C \in
    \C}$ is a stratification compatible with $D_0, \dots, D_K$, and
  let $C \in \C$ be such that $C \subseteq \Pi_m(M)$. \medskip

  \noindent\textbf{Claim:} There are open $U_1, \dots, U_l \subseteq
  U$ and definable $C^p$-maps $f_i:U_i \into \RR^{m-\nu}$ such that
  $\gr f_i \subseteq \Pi_m(M)$ for each $i$ and $C \subseteq \gr f_1
  \cup \cdots \cup \gr f_l$. \medskip

  Assuming the claim, we obtain a finite covering of $\Pi_m(M)$ by
  finitely many definable $C^p$-submanifolds $V$ of $\RR^m$ of
  dimension $\nu$.  For each such $V$, we let $N$ be the set of all $x
  \in M$ for which there exists an open neighbourhood $M_x$ in $M$
  such that $\Pi_m(M_x) \subseteq V$. By the Rank Theorem, $N$ is an
  open subset of $M$ and $\Pi_m(N) = V$, so the lemma is proved.

  To see the claim, we write $C = \gr g$ with $g:D \into \RR^{m-\nu}$
  definable and $C^p$ and $D:= \Pi_\nu(C) \in \D$.  We also put
  $\D_D:= \set{D' \in \D:\ D \subseteq \fr D'}$ and $\C_D:= \set{C'
    \in \C:\ \Pi_\nu(C') \in \D_D \text{ and } C' \subseteq
    \Pi_m(M)}$, and we call a map $s:\D_D \into \C_D$ a
  \textbf{section} if $\Pi_\nu(s(D')) = D'$ for every $D' \in
  \D_D$. To every section $s:\D_D \into \C_D$, we associate a set
  $$A_s:= \gr g \cup \bigcup_{D' \in \D_D} s(D');$$  then $A_s$ is the
  graph of a definable function $g_s:D \cup \bigcup \D_D \into
  \RR^{m-\nu}$.  Since $\D$ is a stratification, the set $D \cup
  \bigcup \D_D$ is open, and we let $D_s$ be the set of all $y \in D$
  such that $g_s$ is of class $C^p$ in a neighbourhood of $y$.  The
  sets $D_s$ are definable open subsets of $D$, and by the hypotheses
  and the Rank Theorem they form a covering of $D$.  Now use definable
  choice to obtain open definable neighbourhoods $U_s$ of $D_s$ such
  that $g_s \rest{U_s}$ is $C^p$, and the claim is proved.
\end{proof}

Next, we let $M \subseteq \RR^n$ be a manifold of dimension $m \le n$.
We also let $d$ be a $p$-distribution on $M$, with $p \le m$, and we
fix $\eta > 0$.

\begin{df} 
  \label{eta-bdd}
  The distribution $d$ is called \textbf{$\eta$-bounded at $x \in M$}
  if there is a matrix $L \in\Ma_{n-p,p}(\RR)$ such that $|L| \le
  \eta$ and $d(x) = \set{(u,L u):\ u \in \RR^p}$.  The distribution
  $d$ is \textbf{$\eta$-bounded} if $d$ is $\eta$-bounded at every $x
  \in M$.
\end{df} 

\begin{rmk}
  If $d$ is $\eta$-bounded, then every integral manifold of $d$ is
  $\eta$-bounded (as defined in our conventions).
\end{rmk}

Given a permutation $\sigma \in \Sigma_n$, the set $\sigma^{-1}(M)$ is
a manifold and the pull-back $\sigma^*d$ is a distribution on
$\sigma^{-1}(M)$; we define
\begin{equation*}
  M_{\sigma,\eta}:= \set{x \in M:\ \sigma^*d
    \text{ is } \eta\text{-bounded at } \sigma^{-1}(x)}.
\end{equation*}
Note that $M_{\sigma,\eta}$ is open in $M$.

\begin{lemma}  
  \label{eta-bdd-cover}
  \begin{enumerate}
  \item If $\eta > 1$, then $M = \bigcup_{\sigma \in \Sigma_n}
    M_{\sigma,\eta}$.  
  \item If $M$ is definable, then so is each $M_{\sigma,\eta}$.
  \end{enumerate}
\end{lemma}

\begin{proof}
  Part (1) follows from the following elementary observation (see
  Lemma 3 of \cite{lio-spe:hausdorff} for details): let $E \subseteq
  \RR^n$ be a linear subspace of dimension $p$.  Then there exist
  $\sigma \in \Sigma_n$ and $L \in M_{n-p,p}(\RR)$ such that $|L|
  \le 1$ and
  $$\sigma^{-1}(E) = \set{(u,Lu) \in \RR^n:\ u \in \RR^p}.$$
  
  For part (2), note that the set $\E_\eta$, consisting of all $E \in
  G^p_n$ for which there exists an $L \in M_{n-p,p}(\RR)$ such that
  $|L| < \eta$ and $E = \{(u,Lu):\ u \in \RR^p\}$, is open and
  semialgebraic.
\end{proof}

The next two lemmas are crucial tools in our use of Hausdorff limits.
For $x \in \RR^n$ and $p \le n$, we set $x_{\le p}:= (x_1, \dots,
x_p)$ and $x_{>p}:= (x_{p+1}, \dots, x_n)$.  Recall that we are
working in the topology of the norm $|\cdot|$ below.

\begin{lemma}
  \label{local_single}
  Let $\eta > 0$, and let $V$ be an $\eta$-bounded submanifold of $M$
  of dimension $p \le m$.  Let $x \in V$, and let $\epsilon > 0$ be
  such that $$\big(B(x_{\le p},\epsilon) \times
  B(x_{>p},p\eta\epsilon)\big) \cap \fr V = \emptyset.$$ Then the
  component of $\,V \cap \big(B(x_{\le p},\epsilon) \times
  B(x_{>p},p\eta\epsilon)\big)$ that contains $x$ is the graph of a
  $p\eta$-Lipschitz function $g:B(x_{\le p},\epsilon) \into B(x_{>p},
  p\eta\epsilon)$.
\end{lemma}

\begin{proof}
  We set $W:= B(x_{\le p},\epsilon)$ and $W':=
  B(x_{>p},p\eta\epsilon)$, and we denote by $C$ the component of $\,V
  \cap (W \times W')$ that contains $x$.  Since $C$ is
  \hbox{$\eta$-bounded}, the map $\Pi_p\rest{C}:C \into W$ is a local
  homeomorphism onto its image.  By general topology, it is therefore
  enough to show that $\Pi_p(C) = W$; we do this by showing that there
  is a function $g:W \into W'$ such that $\gr g \subseteq C$.

  Since $V$ is $\eta$-bounded, there are $\delta > 0$ and a
  $p\eta$-Lipschitz function $g:B(x_{\le p},\delta) \into W'$ such
  that $\gr g \subseteq C$.  We extend $g$ to all of $W$ as follows:
  for $v \in \bd W$, we let $v'$ be the point in the closed line
  segment $[x_{\le p},v]$ closest to $v$ such that $g$ extends to a
  $p\eta$-Lipschitz function $g_v$ along the half-open line segment
  $[x_{\le p},v')$ satisfying $\gr g_v \subseteq V \cap (W \times
  W')$.  Then the proportion of the sidelengths of $W$ and $W'$, the
  $\eta$-boundedness of $V$ and the fact that $(W \times W') \cap \fr
  V = \emptyset$ imply that $v' = v$ for all $v \in \bd W$.  Since the
  graph of the resulting function $g:W \into W'$ is connected and
  contains $x$, it follows that $\gr g \subseteq C$, as required.
\end{proof}

\begin{lemma}  
  \label{local}
  Let $(V_\iota)$ be a sequence of submanifolds of $\,M$ of dimension
  $p\le m$.  Let $\eta>0$, and assume that each $V_\iota$ is
  $\eta$-bounded.  Moreover, assume that both $\lim_\iota V_\iota$ and
  $\lim_\iota \fr V_\iota$ exist and there is an $\nu \in \NN$ such
  that for every $\iota$ and every open box $U \subseteq \RR^n$, the
  set $V_\iota \cap U$ has at most $\nu$ components.  Then for every
  $x \in \lim_\iota V_\iota \setminus \lim_\iota \fr V_\iota$, there
  are a box $U \subseteq \RR^n$ containing $x$ and $p\eta$-Lipschitz
  functions $f_1, \dots, f_\nu:\Pi_p(U) \into \RR^{n-p}$ such that
  \begin{equation*}
    \lim_\iota V_\iota \cap U = (\gr f_1 \cap U) \cup \cdots \cup (\gr
    f_\nu \cap U). 
  \end{equation*}
\end{lemma}

\begin{proof}
  We write ``$\lim$'' in place of ``$\lim_\iota$'' throughout this
  proof.  Let $x \in \lim V_\iota \setminus \lim \fr V_\iota$, and
  choose $\epsilon > 0$ such that $\big(B(x_{\le p},3\epsilon) \times
  B(x_{>p},3p\eta\epsilon)\big) \cap \fr V_\iota = \emptyset$ for all
  $\iota$ (after passing to a subsequence if necessary).  We let $U :=
  B(x_{\le p},\epsilon) \times B(x_{>p},p\eta\epsilon)$, $W:= B(x_{\le
    p},\epsilon)$ and $W':= B(x_{>p},3p\eta\epsilon)$.  Then for each
  $\iota$, the assumptions and Lemma \ref{local_single} imply, with
  $2\epsilon$ in place of $\epsilon$ and each $z \in U \cap V_\iota$
  in place of $x$, that there are definable $p\eta$-Lipschitz
  functions $f_{1,\iota}, \dots, f_{\nu,\iota}:W \into \RR^{n-p}$ such
  that every component of $V_\iota \cap (W \times W')$ intersecting
  $U$ is the graph of some $f_{\lambda,\iota}$.  Moreover, we have
  either $f_{\lambda,\iota} = f_{\lambda',\iota}$ or $\gr
  f_{\lambda,\iota} \cap \gr f_{\lambda',\iota} = \emptyset$, for all
  $\lambda,\lambda' \in \{1, \dots, \nu\}$, and
  \begin{equation*}
    V_\iota \cap U = (\gr f_{1,\iota} \cap U) \cup \cdots \cup (\gr
    f_{\nu,\iota} \cap U).
  \end{equation*}
  Passing to a subsequence if necessary, we may assume that each
  sequence $(f_{\lambda,\iota})_\iota$ converges to a
  $p\eta$-Lipschitz function $f_\lambda:W \into \RR^{n-p}$; then
  $\gr f_\lambda \subseteq \lim V_\iota$.  On the other hand, if $x'
  \in \lim V_\iota \cap U$, then $x' \in \lim (V_\iota \cap U)$, so by
  the above $x' \in \lim (\gr f_{\lambda,\iota} \cap U)$ for some
  $\lambda$, that is, $x' \in \gr f_\lambda$.
\end{proof}

The following lemma is also central to our use of Hausdorff
limits.  We recall that a function $\phi:M \into (0,\infty)$ is a
\textbf{carpeting function on $M$} if $\phi$ is proper and satisfies
$\lim_{x \to y} \phi(x) = 0$ whenever $y \in \fr M$, where the
frontier is taken inside the one-point compactification $\RR^n \cup
\{\infty\}$ of $\RR^n$.  For instance, given positive real numbers
$u_1, \dots, u_n$, the function
\begin{equation*}
  x \mapsto \phi_u(x) := \frac1{1+u_1 x_1^2 + \dots + u_n x_n^2}
\end{equation*}
is a real analytic carpeting function on $\RR^n$.  

\begin{lemma}
  \label{frontier_as_limit}
  Assume that $M$ is bounded and has a carpeting function $\phi$.  Let
  $V$ be a closed subset of $M$, and assume that $V \cap U$ has
  finitely many components for every open box $U \subseteq \RR^n$.
  Then for every sequence $(r_\kappa)_{\kappa \in \NN}$ of positive
  real numbers satisfying $r_\kappa \to 0$ as $\kappa \to \infty$, we
  have
  $$\fr V = \lim_\kappa \left(V \cap
    \phi^{-1}(r_\kappa)\right).$$
\end{lemma}

\begin{proof}
  Let $r_\kappa \to 0$ as $\kappa \to \infty$.  It suffices to show
  that $$\fr V = \lim_{j} \left(\phi^{-1}(r_{\kappa(j)}) \cap
    V\right)$$ for every subsequence $(\kappa(j))_{j \in \NN}$ of
  $(\kappa)$ such that the limit on the right-hand side exists in
  $\K_n$, that is, we may assume that the sequence
  $\left(\phi^{-1}(r_\kappa) \cap V\right)$ converges in $\K_n$.  The
  properties of $\phi$ then imply that $\fr V \supseteq \lim_\kappa
  \left(\phi^{-1}(r_\kappa) \cap V\right)$.  Conversely, let $x \in
  \fr V$.  Since $V \cap B(x,1)$ has finitely many components, there
  is a component $C$ of $V \cap B(x,1)$ such that $x \in \fr C$.  Then
  $C \cup \{x\}$ is connected, so there is a continuous curve
  $\gamma:[0,1] \into C \cup \{x\}$ such that $\gamma([0,1)) \subseteq
  C$ and $\gamma(1) = x$.  Hence $\phi \circ \gamma:[0,1) \into
  (0,\infty)$ is continuous and satisfies $\lim_{t \to 1}
  \phi(\gamma(t)) = 0$, so the intermediate value theorem implies that
  the image $\gamma([0,1))$ intersects $\phi^{-1}(r_\kappa)$ for all
  sufficiently large $\kappa$, so that $x \in \lim_\kappa
  \left(\phi^{-1}(r_\kappa) \cap V\right)$.  Hence $\fr V = \lim_\kappa
  \left(\phi^{-1}(r_\kappa) \cap V\right)$.
\end{proof}

\begin{df}
  \label{convergence_def}
  We abbreviate the conclusion of the lemma by the statement
  $$\fr V = \lim_{r \to 0} \big(V \cap \phi^{-1}(r)\big).$$
\end{df}

\section{Nested distributions}  \label{nested}

We fix an o-minimal expansion $\R$ of the real field.  The goal of
this section and the next is to develop Khovanskii theory for nested
pfaffian sets over $\R$.  We closely follow the exposition of Sections
1 and 2 in \cite{spe:pfaffian}.

Let $M \subseteq \RR^n$ be a $C^2$-manifold of dimension $m$.

\begin{df}[\cite{mou-roc:khovtheory}]
  Let $\D$ be a set of distributions on $M$.  A $C^2$-submanifold $N$
  of $M$ is \textbf{compatible with $\D$} if the pull-back
  $\left(\bigcap_{e \in \E} e\right)^N$ has dimension for every $\E
  \subseteq \D$.  A collection $\C$ of $C^2$-submanifolds of $M$ is
  \textbf{compatible with $\D$} if every $C \in \C$ is compatible with
  $\D$.
\end{df}

\begin{prop}
  \label{partition}
  Let $A_1, \dots, A_k \subseteq \RR^n$ be definable and $p \ge 2$,
  and assume that $M$ is definable and $\D$ is a finite set of
  definable distributions on $M$.  Then there is a finite partition
  (stratification, Whitney stratification) $\P$ of $M$ into definable
  $C^p$-cells such that $\P$ is compatible with each $A_j$ as well as
  $\D$.
\end{prop}

\begin{proof}
  We proceed by induction on $m = \dim M$; the case $m=0$ is trivial.
  So we assume $m > 0$ and the lemma holds for lower values of $m$.
  By $C^2$-cell decomposition, we may assume that $A_1, \dots, A_k$ is
  a partition of $M$ into definable $C^2$-cells.  Thus, for $x \in M$,
  there is a unique $l(x) \in \{1, \dots, k\}$ such that $x \in
  A_{l(x)}$.

  For $x \in M$ and $\E \subseteq \D$ we set $T_x\E := T_x A_{l(x)}
  \cap \bigcap_{d \in \E} d(x)$.  For $\E \subseteq \D$, $j \in \{1,
  \dots, k\}$ and $i \in \{0, \dots, m\}$, we define the set
  \begin{equation*}
    M_{\E,j,i} := \set{x \in A_j:\ \dim T_x\E = i}.
  \end{equation*}
  For each $\E$, the sets $M_{\E,j,i}$ form a covering of $M$, and
  since each $d \in \D$ is definable, each set $M_{\E,j,i}$ is
  definable.  

  Let $\C$ be a partition (stratification, Whitney stratification) of
  $M$ into definable $C^p$-cells compatible with each $M_{\E,j,i}$.
  Then for $C \in \C$, there is a unique $j(C) \in \{1, \dots, k\}$
  such that $C \subseteq A_{j(C)}$.

  Fix a $C \in \C$.  If $\dim C = m$, then for $\E \subseteq \D$ there
  is a unique $i(C,\E) \in \{0, \dots, m\}$ such that $C \subseteq
  M_{\E,j(C),i(C,\E)}$.  Since $C$ is open in $M$, it follows that
  $\left(\bigcap_{d \in \E} d\right)^C$ has dimension $i(C,\E)$ for
  every $\E \subseteq \D$.  On the other hand, if $\dim C < m$, then
  the inductive hypothesis applied to $C$ and $\D^C:= \big\{d^C:\ d
  \in \D\big\}$ in place of $M$ and $\D$ produces a partition
  (stratification, Whitney stratification) $\P_\C$ of $C$ compatible
  with each $A_j$ as well as $\D^C$.  Now it is straightforward to see
  that the collection $$\P := \set{C \in \C:\ \dim C = m} \cup
  \bigcup_{C \in \C,\ \dim C < m} \P_\C$$ is a partition
  (stratification, Whitney stratification) of $M$ compatible with each
  $A_j$ as well as $\D$ (we leave the details to the reader).
\end{proof}

\begin{df}
  \label{integrable_dist}
  Let $d$ be a distribution on $M$ that has dimension.  Then $d$
  gives rise to a vector subbundle $T_{M}d$ of $TM$.  A section
  \hbox{$v:M \into TM$} is a \textbf{vector field on $M$}; $v$ is
  \textbf{tangent to} $d$ if \hbox{$v(M) \subseteq T_Md$} (or
  equivalently, if $v(x) \in d(x)$ for all $x \in M$).  Let
  $\V^1(M,d)$ be the collection of all vector fields on $M$ tangent to
  $d$, and put
  \begin{equation*}
    I(d) := \set{x \in M:\ [v,w](x) \in d(x) \text{ for all } v,w
      \in \V^1(M,d)},
  \end{equation*}
  where $[v,w]$ denotes the Lie bracket of the vector fields $v$ and
  $w$.  We say that $d$ is {\bf integrable} if $I(d) = M$ and that $d$
  is \textbf{nowhere integrable} if $I(d) = \emptyset$.
\end{df}

\begin{nrmks}
  \label{integrable_rmks}
  Let $d$ be a distribution on $M$ that has dimension.
  \begin{enumerate}
  \item The Gauss map $g_M$ is integrable.
  \item Every $1$-distribution on $M$ is integrable.
  \item Let $e$ be a distribution on $M$ such that $e$ and $d \cap e$
    have dimension.  Then $I(d) \cap I(e) \subseteq I(d \cap e)$.  In
    particular, if $d$ and $e$ are integrable, then $d \cap e$ is
    integrable.
  \item Let $N$ be a $C^2$-submanifold of $M$ compatible with $d$.
    Then $I(d) \cap N \subseteq I\big(d^N\big)$; in particular, if $d$
    is integrable, so is $d^N$.
  \item If $M$ and $d$ are definable, then the set $I(d)$ is
    definable.
  \item Let $V$ be an integral manifold of $d$.  Then $V \subseteq
    I(d)$.
  \end{enumerate}
\end{nrmks}

Let $g$ be an integrable $p$-distribution on $M$ with $0 \le p \le m$.
Then by the theorem of Frobenius (see for instance Camacho and Lins
Neto \cite[p. 36]{cam-lin:foliations}), every $x \in M$ belongs to a
unique leaf $L_x^g$ of $g$, and there is an equivalence relation
$\sim_g$ on $M$ associated to $g$ given by
\begin{equation*}
  x \sim_g y \quad\text{iff}\quad L^g_x = L^g_y.
\end{equation*}
Clearly, if $\sim_g$ is definable, then so is $g$; however, the
converse is not true in general.

For the rest of this section, we assume that $M$ is definable, and we
fix a definable nested distribution $d = (d_0, \dots, d_k)$ on $M$.

\begin{dfs}
  \begin{enumerate}
  \item If $N$ is a $C^2$-submanifold of $M$, we say that $N$ is
    \textbf{compatible with $d$} if $N$ is compatible with the set
    $\{d_0, \dots, d_k\}$.
  \item Let $N$ be a $C^2$-submanifold of $M$ compatible with $d$.
    Then $\dim d_i^N \le \dim d_{i-1}^N \le 1+ \dim d_i^N$ for $i = 1,
    \dots, k$.  Thus, we define the \textbf{pull-back $d^N$ of $d$ on
      $N$} as the nested distribution on $N$ obtained by listing the
    set $\set{d_0^N, \dots, d_k^N}$ in decreasing order of dimension.
  \item We call $d$ \textbf{integrable} if each $d_i$ is integrable,
    and we call $d$ \textbf{nowhere integrable} if each $d_i$ is
    nowhere integrable.
  \item We put $\dim d := m-k$, and if $d$ is integrable, we define
    \begin{equation*}
      \deg d := \left|\set{i \in \{0, \dots, k\}:\ \sim_{d_i} \text{ is not
            definable}}\right| \le k,
    \end{equation*}
    called the \textbf{degree} of $d$, where $|\cdot|$ denotes
    cardinality. 
  \end{enumerate}
\end{dfs}

\begin{rmks}
  \begin{enumerate}
  \item If $N$ is a $C^2$-submanifold of $M$ compatible with $d$, then
    $\deg d^N \le \deg d$.
  \item Put $d':= (d_0, \dots, d_{k-1})$.  If $\sim_{d_k}$ is not
    definable, then $\deg d' < \deg d$.
  \end{enumerate}
\end{rmks}

We assume for the rest of this section that $d$ is integrable.

\begin{df}
  \label{core_distribution}
  Let $e = (e_0, \dots, e_l)$ be an integrable, definable nested
  distribution on $M$ with $l \le k$.  We call $e$ a \textbf{core
    distribution of $d$} if
  \begin{renumerate}
  \item $\sim_{d_i}$ is definable for $i=1, \dots, k-l$, and
  \item $d_i = d_{k-l} \cap e_{i-k+l}$ for $i=k-l+1, \dots, k$.
  \end{renumerate}
\end{df}

\begin{rmks}
  \begin{enumerate}
  \item Let $e$ be a core distribution of $d$.  Then $\deg d \le \deg
    e$.  Moreover, if $f$ is a core distribution of $e$, then $f$ is
    also a core distribution of $d$.
  \item Let $N$ be a definable $C^2$-submanifold of $M$ compatible
    with $d$, and let $e$ be a core distribution of $d$.  If $N$ is
    compatible with $e$, then $d^N$ has core distribution $e^N$.
  \end{enumerate}
\end{rmks}

In our setting, core distributions typically arise in the following
way:

\begin{expl} 
  \label{Morse_function} 
  Let $\phi:M \into \RR$ be definable and $C^2$, and define $g_\phi:M
  \into G_n$ by $g_\phi(x):= \ker d\phi(x) \subseteq T_xM$; note that
  $g_\phi(x)$ has codimension at most 1 in $T_xM$.  We set $\D_\phi :=
  \{g_M, g_\phi, d_1 \cap g_\phi, \dots, d_k \cap g_\phi\}$, and for a
  $C^2$-submanifold $N$ of $M$ compatible with $\D_\phi$, we let
  $d_\phi^N$ be the definable nested distribution on $N$ obtained by
  listing the pull-backs of all elements of $\D_\phi$ to $N$ in order
  of decreasing dimension.  If, in addition, $N$ is compatible with
  $d$, then $d_\phi^N$ has core distribution $d^N$ and $\dim d^N \ge
  \dim d_\phi^N \ge \dim d^N -1$; in particular, $\deg d_\phi^N \le
  \deg d$.

  We now let $\C$ be a Whitney stratification of $M$ by definable
  $C^2$-cells compatible with both $d$ and $\D_\phi$, as obtained from
  Proposition \ref{partition}.  Let $\C'$ be the set of all $C \in \C$
  such that $g_C \cap d_k \nsubseteq g_\phi$. Then $\dim d_\phi^C <
  \dim d$ for $C \in \C'$, and we claim that the union of all cells in
  $\C'$ is an open subset $M'$ of $M$.  To see this, note first that
  if $C,D \in \C$ are such that $D \subseteq \fr C$, then the Whitney
  property of the pair $(C,D)$, as defined on p. 502 of
  \cite{vdd-mil:cat}, implies that for every sequence $(x_i)_{i \in
    \NN}$ of points in $C$ that converges to a point $y \in D$ and for
  which $T:= \lim_i T_{x_i}C$ exists in $G_n^{\dim C}$, we have $T_yD
  \subseteq T$.  Since $d$ and $g_\phi$ are continuous, it follows
  that the union of all cells in $\C \setminus \C'$ is a closed subset
  of $M$; hence $M'$ is an open subset of $M$.  Finally, note that
  $d_k(x) \nsubseteq g_\phi(x)$ for all $x \in M'$, so that $M'$ is
  compatible with $\D_\phi$ and $\dim d_\phi^{M'} < \dim d^{M'}$.
\end{expl}

\section{Khovanskii theory}
\label{khovsection}

Let $M \subseteq \RR^n$ be a definable $C^2$-manifold of dimension
$m$.  We fix a finite family $\Delta = \{d^1, \dots, d^q\}$ of
definable nested distributions on $M$; we write $d^p = (d^p_0, \dots,
d^p_{k(p)})$ for $p = 1, \dots, q$.  We associate to $\Delta$ the
following set of distributions on $M$:
\begin{equation*}
  \D_\Delta := \set{d^0_0 \cap d^1_{k(1)} \cap \cdots \cap
    d^{p-1}_{k(p-1)} \cap d^p_j:\ p = 1, \dots, q \text{ and } j = 0,
    \dots, k(p)}, 
\end{equation*}
where $d^0_0 = g_M$.  If $N$ is a $C^2$-submanifold of $M$ compatible
with $\D_\Delta$, we let $d^{\Delta,N} = \left(d^{\Delta,N}_0, \dots,
  d^{\Delta,N}_{k(\Delta,N)}\right)$ be the nested distribution on $N$
obtained by listing the set $\set{g^N:\ g \in \D_\Delta}$ in order of
decreasing dimension.  In this situation, if $V_p$ is an integral
manifold of $d^p_{k(p)}$, for $p=1, \dots, q$, then the set $N \cap
V_1 \cap \cdots \cap V_q$ is an integral manifold of
$d^{\Delta,N}_{k(\Delta,N)}$.

\begin{lemma}
  \label{lowering dimension}
  Let $N \subseteq M$ be a definable $C^2$-cell compatible with
  $\D_\Delta$, and suppose that $\dim d^{\Delta,N}_{k(\Delta,N)} > 0$.
  Then there is a definable carpeting function $\phi$ on $N$ of class
  $C^2$ such that the definable set
  \begin{equation*}
    B := \set{x \in N:\ d^{\Delta,N}_{k(\Delta,N)}(x) \subseteq \ker
      d\phi(x)}
  \end{equation*}
  has dimension less than $\dim N$. 
\end{lemma}

\begin{proof}
  By \cite{vdd-mil:cat} there is a definable diffeomorphism
  $\sigma:\RR^{\dim N} \into N$ of class $C^2$.  Replacing $n$ by
  $\dim N$, $N$ by $\RR^{\dim N}$ and each $d^{\Delta,N}_i$ by its
  pull-back $\sigma^*d^{\Delta,N}_i$, we reduce to the case where $N =
  M = \RR^n$ and write $k$ and $d$ in place of $k(\Delta,N)$ and
  $d^{\Delta, N}$.  Then for $u = (u_1, \dots, u_n) \in (0,\infty)^n$
  we put
  \begin{equation*}
    B_u := \set{x \in \RR^n:\ d_k(x) \subseteq
      \ker d\phi_u(x)},
  \end{equation*}
  where $\phi_u$ is the carpeting function defined on $\RR^n$ before
  Lemma \ref{frontier_as_limit}.  If $\dim B_u < n$ for some $u$ as
  above, the proof is finished.  So assume for a contradiction that
  $\dim B_u = n$ for all $u$ as above.  Then $\dim B = 2n$, where
  \begin{equation*}
    B := \set{(u,x) \in (0,\infty)^n \times \RR^n:\ x \in B_u},
  \end{equation*}
  so there are nonempty open $V \subseteq (0,\infty)^n$ and $W
  \subseteq \RR^n$ such that $V \times W \subseteq B$.  Fix some $x
  \in W$ with all $x_i \neq 0$ and let $u$ range over $V$.  Note that
  \begin{equation*}
    d\phi_u(x) = -\phi_u(x)^2 \big(2 u_1 x_1 dx_1 + \dots + 2 u_n
    x_n dx_n\big).
  \end{equation*}
  Therefore the vector space generated by all $d\phi_u(x)$ as $u$
  ranges over $V$ has dimension $n$, that is, the intersection of all
  $\ker d\phi_u(x)$ as $u$ ranges over $V$ is trivial, which
  contradicts $\dim d > 0$.
\end{proof}

For convenience, if $e = (e_0, \dots, e_l)$ is a nested distribution
on $M$ and $W \subseteq M$, we call $W$ a \textbf{Rolle leaf of $e$}
if there is a nested Rolle leaf $V = (V_0, \dots, V_l)$ of $e$ such
that $W = V_l$.

\begin{prop}
  \label{components}
  Let $A \subseteq \RR^n$ be a definable set.  Then there exists a $K
  \in \NN$ such that, whenever $L_p$ is a Rolle leaf of $d^p$ for $p =
  1, \dots, q$, then $A \cap L_1 \cap \dots \cap L_q$ is a union of at
  most $K$ connected manifolds.
\end{prop}

\begin{proof}
  We proceed by induction on $\dim A$ and $k:= k(1) + \cdots + k(q)$.
  The cases $\dim A =0$ or $k=0$ being trivial, we assume that $\dim
  A>0$ and $k>0$ and that the result holds for lower values of $\dim
  A$ or $k$.  After shrinking $q$, we may also assume that $k(p)>0$
  for each $p$.  By Proposition \ref{partition}, it suffices to
  consider the case where $A = N$ is a $C^2$-cell contained in $M$ and
  compatible with $\D_\Delta$.  For $p = 1, \dots, q$, we let $L_p$ be
  a Rolle leaf of $d^p$, and we put $L:= L_1 \cap \dots \cap L_q$;
  then $N \cap L$ is an integral manifold of $d:=
  d^{\Delta,N}_{k(\Delta,N)}$.  \smallskip

  \noindent\textbf{Case $\dim d = 0$.}
  Let $\Delta':= \set{d^1, \dots, d^{q-1}, \left(d^q_0, \dots,
      d^q_{k(q)-1}\right)}$, let $L'_q$ be the Rolle leaf of
  $d^q_{k(q)-1}$ containing $L_q$ and put $L' := L_1 \cap \dots \cap
  L_{q-1} \cap L'_q$.  Then $N$ is compatible with $\D_{\Delta'}$ and
  $N \cap L'$ is an integral manifold of
  $d':=d^{\Delta',N}_{k(\Delta',N)}$ of dimension at most $1$.  By the
  inductive hypothesis, there is a $K \in \NN$ (depending only on $N$
  and $\Delta'$, but not on the particular Rolle leaves) such that the
  manifold $N \cap L'$ has at most $K$ components.  Thus, if $\dim(N
  \cap L') = 0$, we are done by the inductive hypothesis, so we assume
  that $\dim(N \cap L') = 1$.  Since $N$ is compatible with
  $\D_{\Delta'}$, it follows that $\dim d' = 1$ as well.

  Let $C$ be a component of $N \cap L'$.  If $C \cap L_q$ contains
  more than one point, then by the Rolle property of $L_q$ in $L'_q$
  and the fact that $C$ is a connected $C^1$-submanifold of $L'_q$ of
  dimension $1$, $C$ is tangent at some point to
  $\big(d^q_{k(q)}\big)^{L'_q}$, which contradicts the assumption that
  $\dim d = 0$.  So $C \cap L_q$ contains at most one point for each
  component $C$ of $N \cap L'$.  Hence $N \cap L$ consists of at most
  $K$ points.  \medskip

  \noindent\textbf{Case $\dim d > 0$.}
  Let $\phi$ and $B$ be obtained from Lemma \ref{lowering dimension}.
  Then $\dim B < \dim A$; so by the inductive hypothesis, there is a
  $K \in \NN$, independent of the particular Rolle leaves chosen, such
  that $B \cap L$ has at most $K$ components.  Since $N \cap L$ is a
  closed, embedded submanifold of $N$, $\phi$ attains a maximum on
  every component of $N \cap L$, and any point in $N \cap L$ where
  $\phi$ attains a local maximum belongs to $B$.  Hence $N \cap L$ has
  at most $K$ components.
\end{proof}

\begin{cor}
  \label{khovprop}
  \begin{enumerate}
  \item Let $\C$ be a partition of $M$ into definable $C^2$-cells
    compatible with $\D_\Delta$.  Then there is a $K \in \NN$ such
    that, for every $C \in \C$ and every Rolle leaf $L_p$ of $d^p$
    with $p=1, \dots, q$, the set $C \cap L_1 \cap \cdots \cap L_q$ is
    a union of at most $K$ Rolle leaves of $d^{\Delta,C}$.
  \item Let $\A$ be a definable family of sets.  Then there is a $K
    \in \NN$ such that whenever $A \in \A$ and $L_p$ is a Rolle leaf
    of $d^p$ for each $p$, the set $A \cap L_1 \cap \dots \cap L_q$ is
    a union of at most $K$ connected manifolds.
  \end{enumerate}
\end{cor}

\begin{proof}
  Part (1) follows from the proof of Proposition \ref{components}.
  For (2), we let $A \subseteq \RR^{m+n}$ be definable such that $\A =
  \set{A_z:\ z \in \RR^m}$ where, for $z \in \RR^m$, $A_z:= \set{x \in
    \RR^n:\ (z,x) \in A}$.  We replace $M$ by $M':= \RR^m \times M$
  and each $d^p_i$ by the distribution $e^p_i$ on $M'$ defined by
  $e^p_i(z,x):= \RR^m \times d^p_i(x)$.  Moreover, we let $e$ be the
  nested distribution on $\RR^{m+n}$ obtained as in Example 1 from the
  family $\Omega = (dz_1, \dots, dz_m)$, and we let $e^{q+1}$ be the
  restriction of $e$ to $M'$ and put $\Delta':= \{e^1, \dots,
  e^{q+1}\}$.  By Proposition \ref{components}, there is a $K'$ such
  that whenever $L'_p$ is a Rolle leaf of $e^p$, for $p = 1, \dots,
  q+1$, then $A \cap L'_1 \cap \cdots \cap L'_{q+1}$ is the union of
  at most $K'$ connected manifolds.  But for every Rolle leaf $L_p$ of
  $d^p$, with $p \in \{1, \dots, q\}$, the set $\RR^m \times L_p$ is a
  Rolle leaf of $e^p$; and for every $z \in \RR^m$ and each component
  $C$ of $M$, the set $\{z\} \times C$ is a Rolle leaf of $e^{q+1}$.
  Thus, we can take $K = K' \cdot l$, where $l$ is the number of
  components of $M$.
\end{proof}

\begin{cor}
  \label{pfaff_intersection}
  Let $X \subseteq \RR^n$, $X_1 \subseteq \RR^{n_1}$ and $X_2
  \subseteq \RR^{n_2}$ be nested pfaffian over $\R$.
  \begin{enumerate}
  \item Each component of $X$ is nested pfaffian over $\R$.
  \item If $n_1 = n_2$, then $X_1 \cap X_2$ is nested pfaffian over $\R$.
  \item The product $X_1 \times X_2$ is nested pfaffian over $\R$.
  \end{enumerate}
\end{cor}

\begin{proof}
  Part (1) follows from Corollary \ref{khovprop}(1).  For (2), we may
  assume that $X_1$ and $X_2$ are basic nested pfaffian over $\R$.
  Let $M_1, M_2 \subseteq \RR^n$ be definable $C^2$-manifolds with $n
  = n_1 = n_2$, and for $p=1,2$, let $d^p = (d^p_0, \dots,
  d^p_{k(p)})$ be definable nested distributions on $M_p$, $L_p$ be a
  Rolle leaf of $d^p$ and $A_p \subseteq \RR^n$ be definable such that
  $X_p = A_p \cap L_p$.  Let $\C$ be a $C^2$-cell decomposition of
  $\RR^n$ compatible with $M_1, M_2, M_1 \cap M_2, A_1$ and $A_2$.
  Refining $\C$ if necessary, we may also assume that if $C \in \C$ is
  such that $C \subseteq M_1 \cap M_2$, then $C$ is compatible with
  both $d^1$ and $d^2$.  Then it follows from Corollary
  \ref{khovprop}(1) that we may assume that $M_1 = M_2 = C$ for each
  such $C \in \C$.  In this case, we put $\Delta := \{d^1,d^2\}$ and
  obtain again from Corollary \ref{khovprop}(1) that $X_1 \cap X_2$ is
  a finite union of basic nested pfaffian sets over $\R$.

  For (3), we argue as in the proof of Corollary \ref{khovprop}(2),
  but without adding the extra nested distribution $e^{q+1}$ there, to
  see that $\RR^{n_1} \times X_2$ and $X_1 \times \RR^{n_2}$ are
  nested pfaffian over $\R$.  Hence $X_1 \times X_2 = \left(X_1 \times
    \RR^{n_2}\right) \cap \left(\RR^{n_1} \times X_2\right)$ is nested
  pfaffian over $\R$ by part (2).
\end{proof}

\begin{cor}  
  \label{pfaff_rolle}
  Let $A \subseteq M$ be a definable set and $d = (d_0, \dots, d_k)$
  be a definable nested distribution on $M$.  Let $p \in \NN \cup
  \{\infty,\omega\}$ with $p \ge 2$, and assume that $\R$ admits
  $C^p$-cell decomposition.  Then there are $N,s \in \NN$ and a
  collection $\set{(C_j, \psi_j, e^j):\ j = 1, \dots, s}$ such that
  $\{C_1, \dots, C_s\}$ is a $C^p$-cell decomposition of $A$ and for
  $j=1, \dots, s$,
  \begin{renumerate}
  \item $\psi_j:\RR^{n_j} \into C_j$ is a definable
    $C^p$-diffeomorphism, where $n_j := \dim C_j$;
  \item $e^j = \left(e^j_0, \dots, e^j_{k(j)}\right)$ is a definable
    nested distribution on $\RR^{n_j}$ of class $C^{p-1}$;
  \item whenever $\,V$ is a Rolle leaf of $d_k$, there are (not
    necessarily pairwise distinct) Rolle leaves $V^j_r$ of $e^j$, for
    $j=1, \dots, s$ and $r = 1, \dots, N$, such that $A \cap V =
    \bigcup_{j,r} \psi_j\left(V^j_r\right)$.
  \end{renumerate}
\end{cor}

\begin{proof}
  Let $\C$ be a definable $C^p$-cell decomposition of $\RR^n$
  compatible with $M$, $A$ and $d$.  Note that every $C \in \C$ is
  definably $C^p$-diffeomorphic to $\RR^{\dim C}$.  The corollary now
  follows from Corollary \ref{khovprop}(1).
\end{proof}

Finally, we address the issue of definability in the pfaffian closure
$\P(\R)$ of $\R$; we adopt here the terminology of Section 4 in
\cite{spe:pfaffian}. 

\begin{prop}  
  \label{definability}
  Let $W$ be a Rolle leaf over $\P(\R)$.  Then $W$ is definable in
  $\P(\R)$.
\end{prop}

\begin{proof}
  Let $N \subseteq \RR^n$ be a $C^2$-manifold definable in $\P(\R)$,
  and let $d = (d_0, \dots, d_k)$ be a nested distribution on $N$
  definable in $\P(\R)$ such that $W$ is a Rolle leaf of $d$.  We
  proceed by induction on $k$; the case $k=0$ is trivial, so we assume
  $k>0$ and the proposition holds for lower values of $k$.  By
  definition of $\P(\R)$ and $\R_i$, there is an $i \in \NN$ such that
  $N$ and $d$ are definable in $\R_i$.  By Corollary
  \ref{pfaff_rolle}, with $\R_i$ in place of $\R$, we may assume that
  $N = \RR^n$.  Let $\Omega$ be associated to $d$ as in Example 1
  (with $\RR^n$ in place of $M$), and let $V = (V_0, \dots, V_k)$ be a
  nested Rolle leaf of $d$ such that $W = V_k$.  Then $V_1 \in
  \la(\R_i)$, so $V_1$ is definable in $\R_{i+1}$.  Now let $\C$ be a
  decomposition of $V_1$ into $C^2$-cells definable in $\R_{i+1}$ and
  compatible with $d$.  Then $d^C = (d^C_0, \dots, d^C_{k(C)})$ is
  definable in $\P(\R)$ with $k(C) < k$, and Corollary
  \ref{khovprop}(1), with $\R_{i+1}$ in place of $\R$, implies that $C
  \cap W$ is a finite union of Rolle leaves of $d^C$.  The proposition
  therefore follows from the inductive hypothesis.
\end{proof}

\section{Pfaffian limits}  
\label{limits_section}

Let $M \subseteq \RR^n$ be a \textit{bounded}, definable
$C^2$-manifold of dimension $m$ and $d = (d_0, \dots, d_k)$ be a
definable and integrable nested distribution on $M$.

\begin{df}
  \label{admissible_manifold}
  A nonempty integral manifold $V$ of $d_k$ is an \textbf{admissible
    integral manifold of $d$} if $d$ has a core distribution $e =
  (e_0, \dots, e_l)$ and there are a definable, closed integral
  manifold $B = B(V)$ of $d_{k-l}$ and a Rolle leaf $W = W(V)$ of $e$
  such that $V = W \cap B$.  In this situation, $W$ (but not
  necessarily $B$) is uniquely determined by $V$ and $e$, and we call $W$ the
  \textbf{core} of $V$ \textbf{corresponding to $e$} and $B$ a
  \textbf{definable part} of $V$ \textbf{corresponding to $W$}.
\end{df}

Whenever they are clear from context, we shall omit the phrases
``corresponding to $e$'' or ``corresponding to $W$''.

\begin{rmks}
  (1) Since $d$ is integrable, an integral manifold $L$ of $d_k$
  definable in $\P(\R)$ is a leaf of $d_k$ if and only if $L$ is
  connected and closed in $M$.  Hence by Corollary \ref{khovprop}(1),
  every admissible integral manifold of $d$ is a closed submanifold of
  $M$ and a finite union of leaves of $d_k$.

  (2) Let $V$ be an admissible integral manifold of $d$.  Then there
  are admissible integral manifolds $V_j$ of $(d_0, \dots, d_j)$, for
  $j=0, \dots, k$, such that $V_0 \supseteq \cdots \supseteq V_k$ and
  $V = V_k$.  To see this, let $e = (e_0, \dots, e_l)$, $W$ and $B$ be
  such that $W$ is the core of $V$ corresponding to $e$ and $B$ is a
  definable part of $V$ corresponding to $W$.  Let $(W_0, \dots, W_l)$
  be the nested Rolle leaf of $e$ such that $W = W_l$.  Then for
  $j=k-l+1, \dots, k$, we put $V_j := W_{j-k+l} \cap B$, an admissible
  integral manifold of $d_j$ with core $W_{j-k+l}$ corresponding to
  $(e_0, \dots, e_{j-k+l})$ and definable part $B$.  For $j=0, \dots,
  k-l$, we let $V_j$ be the smallest union of leaves of $d_j$
  containing $B$, a definable, closed integral manifold of $d_j$.
\end{rmks}

\begin{nrmk}
  \label{admissible_finiteness}
  Let $e = (e_0, \dots, e_l)$ be a core distribution of $d$, and let
  $W$ be a Rolle leaf of $e$.  Let $N$ be a definable
  $C^2$-submanifold of $M$ compatible with $d$ and $e$, so that
  $e^{N}$ is a core distribution of $d^{N}$.  By Corollary
  \ref{khovprop}(1), $W \cap N$ is a finite union of Rolle leaves
  $W^N_1, \dots, W^N_q$ of $e^N$.

  Let now $V$ be an admissible integral manifold of $d$ with core $W$
  corresponding to $e$ and definable part $B$ corresponding to $W$.
  Then $V \cap N = V^N_1 \cup \cdots \cup V^N_q$ where, for $p=1,
  \dots, q$, the set $V^N_p:= W^N_p \cap B$ is an admissible integral
  manifold of $d^N$ with core $W^N_p$ corresponding to $e^N$ and
  definable part $B \cap N$.
\end{nrmk}

\begin{df}
  \label{admissible_sequence}
  A sequence $(V_\iota)_{\iota \in \NN}$ of integral manifolds of
  $d_k$ is an \textbf{admissible sequence of integral manifolds of
    $d$} if there are a core distribution $e = (e_0, \dots, e_l)$ of
  $d$, a Rolle leaf $W$ of $e$ and a definable family $\B$ of closed
  integral manifolds of $d_{k-l}$ such that each $V_\iota$ has core
  $W$ corresponding to $e$ and definable part in $\B$ corresponding to
  $W$.  In this situation, we call $W$ the \textbf{core} of the
  sequence $(V_\iota)$ \textbf{corresponding to $e$} and $\B$ a
  \textbf{definable part} of the sequence $(V_\iota)$
  \textbf{corresponding to $W$}.
\end{df}

\begin{rmks}
  \begin{lenumerate}{}{2}
  \item In the previous definition, we think of the core of
    $(V_\iota)$ as representing the ``non-definable part'' of
    $(V_\iota)$.  
  \item Let $(V_\iota)$ be an admissible sequence of integral
    manifolds of $d$.  Arguing as in Remark (2) above, we see that
    there is an admissible sequence $(U_\iota)$ of integral manifolds
    of $(d_0, \dots, d_{k-1})$ such that $V_\iota \subseteq U_\iota$
    for $\iota \in \NN$.
  \end{lenumerate}
\end{rmks}

\begin{df}
  \label{pfaffian_limit}
  Let $(V_\iota)$ be an admissible sequence of integral manifolds of
  $d$.  If $(V_\iota)$ converges to $K \in \K_n$, we call $K$ a
  \textbf{pfaffian limit over $\R$}.  In this situation, we say that
  $K$ \textbf{is obtained from} $d$, and we put
  \begin{equation*}
    \deg K := \min\set{\deg f:\ K \text{ is obtained from } f}.
  \end{equation*}
\end{df}

\begin{lemma} 
   \label{NPL_dimension} 
   Let $K \subseteq \RR^n$ be a pfaffian limit obtained from $d$.
   Then $K$ is definable in $\P(\R)$ and $\dim K \le \dim d$.
\end{lemma}

\begin{proof}
  Let $(V_\iota)$ be an admissible sequence of integral manifolds of
  $d$ such that $K = \lim_\iota V_\iota$.  Let $e$, $W$ and $\B$ be
  such that $W$ is the core of $(V_\iota)$ corresponding to $e$ and
  $\B$ is a definable part of $(V_\iota)$ corresponding to $W$.  Since
  $W$ is definable in $\P(\R)$, the family of all admissible integral
  manifolds of $d$ with core $W$ corresponding to $e$ and definable
  part in $\B$ corresponding to $W$ is definable in $\P(\R)$.  Hence
  the lemma follows from the versions of the Marker-Steinhorn theorem
  \cite{mar-ste:deftypes} found in \cite[Theorem 3.1]{vdd:hausdorff}
  and \cite[Theorem 1]{lio-spe:hausdorff}.
\end{proof}

For the remainder of this section, we let $e = (e_0, \dots, e_l)$ be a
core distribution of $d$ and $W$ be a Rolle leaf of $e$. 

One reason for working with pfaffian limits over $\R$ is that they can
be used to describe the frontiers of admissible integral manifolds.
To see how this is done, we assume for Propositions
\ref{frontier_reduction} and \ref{frontier_lemma} below that $M$ has a
definable $C^2$-carpeting function $\phi$, and we adopt the
corresponding notations introduced in Example \ref{Morse_function}.
We assume that the Whitney stratification $\C$ is also compatible with
$e$, and we let $N \in \C$ or $N= M'$ and adopt here the corresponding
notations introduced in Remark \ref{admissible_finiteness}.  By not
requiring that all $W^N_p$ be distinct, we may assume that $q$ is
independent of $N$.  Thus, for $r>0$, $p=1, \dots, q$ and every
admissible integral manifold $V$ of $d$ with core $W$ and definable
part $B$, each set $\phi^{-1}(r) \cap V^N_p$ is an admissible integral
manifold of $d_\phi^N$ with core $W^N_p$ and definable part
$\phi^{-1}(r) \cap B \cap N$.

Let now $V$ be an admissible integral manifold of $d$ with core $W$
and definable part $B$, and let $(r_\kappa)_{\kappa \in \NN}$ be a
sequence of positive real numbers such that $r_\kappa \to 0$ and
$\lim_\kappa \left(\phi^{-1}(r_\kappa) \cap V\right)$ exists.  Then
for $p = 1, \dots, q$, the sequence $\left(\phi^{-1}(r_\kappa) \cap
  V^N_p\right)$ is an admissible sequence of integral manifolds of
$d_\phi^N$ with core $W^N_p$ and definable part $\set{\phi^{-1}(r)
  \cap B \cap N:\ r>0}$.  Passing to a subsequence, we may assume that
$K^N_p(V):= \lim_\kappa \left(\phi^{-1}(r_\kappa) \cap V^N_p\right)$
exists for each $p$ and each $N$.  Thus, each $K^N_p(V)$ is a pfaffian
limit obtained from $d^N$.

\begin{prop}  
  \label{frontier_reduction}
  $ \fr V = K^{M'}_1(V) \cup \cdots \cup K^{M'}_q(V)$, and each
  $K^{M'}_p(V)$ is a pfaffian limit over $\R$ such that $\dim
  K^{M'}_p(V) < \dim d$ and $\deg K^{M'}_p(V) \le \deg d$.
\end{prop}

\begin{proof}
  Since $V$ is definable in $\P(\R)$, we get from Lemma
  \ref{frontier_as_limit} that $\fr V = \lim_\kappa
  \left(\phi^{-1}(r_\kappa) \cap V\right)$.  On the other hand, we have
  \begin{equation*}
    \bigcup_{p=1}^q K^{M'}_p(V) = \bigcup_{C \in \C'}
    \bigcup_{p=1}^q K^N_p(V),
  \end{equation*}
  and both unions are contained in $\fr V$.  Thus, to finish our
  argument, we let $z \in \lim_\kappa \left(\phi^{-1}(r_\kappa) \cap
    V\right)$ and show that $z \in K^N_p(V)$ for some $N \in \C'$ and
  $p \in \{1, \dots, q\}$.  Let $x_\kappa \in \phi^{-1}(r_\kappa) \cap
  V$ be such that $\lim x_\kappa = z$.  Let $N \in \C$ be such that
  infinitely many $x_\kappa$ belong to $N$; passing to a subsequence,
  we may assume that $x_\kappa \in N$ for all $\kappa$.  Then $N \in
  \C'$: otherwise, we have $g_N \cap d_k \subseteq d_\phi$, so the
  definability in $\P(\R)$ implies that $\phi^{-1}(r) \cap V \cap N =
  \emptyset$ for all but finitely many $r$, which contradicts
  $\lim_\kappa \phi(x_\kappa) = 0$.  Thus, passing again to a
  subsequence, we may assume that there is a $p$ such that $x_\kappa
  \in \phi^{-1}(r_\kappa) \cap V^N_p$ for all $\kappa$.  Hence $z \in
  K^N_p(V)$, as required.
\end{proof}

\begin{prop}  
  \label{frontier_lemma}
  Let $(V_\iota)$ be an admissible sequence of integral manifolds of
  $d$ such that $K':= \lim_\iota \fr V_\iota$ exists.  Then there are
  $q \in \NN$ and pfaffian limits $K_1, \dots, K_q \subseteq \RR^n$
  over $\R$ such that $K' = K_1 \cup \cdots \cup K_q$ and $\dim K_p <
  \dim d$ and $\deg K_p \le \deg d$ for each $p$.
\end{prop}

\begin{proof}
  We may assume that $(V_\iota)$ has core $W$, and we adopt the
  notations introduced before Proposition \ref{frontier_reduction},
  with $V_\iota$ in place of $V$ and $B(V_\iota) \in \B$, where $\B$
  is a definable part of $(V_\iota)$.  Passing to a subsequence if
  necessary, we may assume that $\lim_\iota K^{M'}_p(V_\iota)$ exists
  for each $p$.  Then by Proposition \ref{frontier_reduction}, we have
  \begin{equation*}
    K' = \bigcup_{p=1}^q
    \lim_\iota K^{M'}_p(V_\iota). 
  \end{equation*} 
  Hence $K' = \bigcup_{p=1}^q \lim_\iota
  \left(\phi^{-1}(r_{\kappa(\iota)}) \cap (V_\iota)^{M'}_p\right)$ for
  some subsequence $(\kappa(\iota))_\iota$.  But each sequence
  $\left(\phi^{-1}(r_{\kappa(\iota)}) \cap
    (V_\iota)^{M'}_p\right)_\iota$ is an admissible sequence of
  integral manifolds of $d_\phi^{M'}$ with core $W^{M'}_p$ and
  definable part $\{\phi^{-1}(r) \cap B \cap N:\ r>0,\ B \in \B\}$, as
  required.
\end{proof}

Another reason for working with pfaffian limits over $\R$ is that they
are themselves well behaved with respect to taking frontiers after
intersecting with definable cells.  To see this, we define the
following distributions on $\mathbf{M}:= M \times \RR^2$, where we
write $(x,r,\epsilon)$ for the typical element of $\mathbf{M}$ with $x
\in M$ and $r,\epsilon \in \RR$: we set $\mathbf{d}_0 := g_{\mathbf
  M}$, $\mathbf{d}_1 := \ker d\epsilon \cap \mathbf{d}_0$,
$\mathbf{d}_2 := \ker dr \cap \mathbf{d}_1$ and
\begin{align*}
  \mathbf{d}_{2+j}(x,r,\epsilon) &:= \left(d_j(x) \times
    \RR^{2}\right) \cap
  \mathbf{d}_{2}(x,r,\epsilon) &\text{ for } j=1, \dots, k, \\
  \mathbf{e}_j(x,r,\epsilon) &:= e_j(x) \times \RR^{2} &\text{for }
  j=0, \dots, l.
\end{align*}
We also put $\mathbf{W}:= W \times \RR^2$.  Then $\mathbf{d}:=
(\mathbf{d}_0, \dots, \mathbf{d}_{2+k})$ is a definable nested
distribution on $\mathbf{M}$ with core distribution $\mathbf{e} :=
(\mathbf{e}_0, \dots, \mathbf{e}_l)$.  Thus $\deg \mathbf{d} \le \deg
d$, and whenever $(V_\iota)$ is an admissible sequence of integral
manifolds of $d$ with core $W$ and $(r_\iota,\epsilon_\iota) \in
\RR^2$ for $\iota \in \NN$, the sequence $\big(V_\iota \times
\{(r_\iota,\epsilon_\iota)\}\big)$ is an admissible sequence of
integral manifolds of $\mathbf{d}$ with core $\mathbf{W}$.

\begin{lemma}  
  \label{frontier_of_NPL}
  Let $K$ be a pfaffian limit obtained from $d$, and let $C \subseteq
  \RR^n$ be a definable cell.  Then there is a definable open subset
  $\mathbf{N}$ of $\,\mathbf{M}$ and there are $q \in \NN$ and
  pfaffian limits $K_1, \dots, K_q$ obtained from
  $\mathbf{d}^{\mathbf{N}}$ such that $$\fr(K \cap C) = \Pi_n(K_1)
  \cup \cdots \cup \Pi_n(K_q);$$ in particular, $\deg K_p \le \deg d$
  for each $p$.
\end{lemma}

The following general observation is needed for the proof of this
lemma:

\begin{nrmk}
  \label{limit_rmks}
  In the situation of Remark \ref{admissible_finiteness}, let
  $(V_\iota)$ be an admissible sequence of integral manifolds of $d$
  with core $W$, and assume that $K^N := \lim_\iota (V_\iota \cap N)$
  exists.  After passing to a subsequence if necessary, we may assume
  that the sequence $\big((V_\iota)^N_p\big)_\iota$ converges to a set
  $K^N_p \in \K_n$, for $p=1, \dots, q$; then $$K^N = K^N_1 \cup
  \cdots \cup K^N_q.$$ Moreover, if $\B$ is the definable part of
  $(V_\iota)$, then each sequence $\big((V_\iota)^N_p\big)_\iota$ is
  an admissible sequence of integral manifolds of $d^N$ with core
  $W^N_p$ corresponding to $e^N$ and definable part $\B^N := \set{B
    \cap N:\ B \in \B}$.  Thus, each $K^N_p$ is a pfaffian limit
  obtained from $d^N$; in particular, $\deg K^N_p \le \deg d$.
\end{nrmk}

\begin{proof}[Proof of Lemma \ref{frontier_of_NPL}]
  We let $\phi$ be a definable carpeting function on $C$ and put
  \begin{equation*}
    \mathbf{N}:= \set{(x,r,\epsilon) \in \mathbf{M}:\
      d(x,\phi^{-1}(r)) < \epsilon},
  \end{equation*}
  where we set $d(x,\emptyset) := \infty$ for all $x \in M$.  Then
  $\mathbf{N}$ is an open, definable subset of $\mathbf{M}$, and since
  $K$ is compact and definable in $\P(\R)$, we have $\fr(K \cap C) =
  \lim_{r \to 0} (\phi^{-1}(r) \cap K)$ by Lemma
  \ref{frontier_as_limit}.  Moreover, we let $(V_\iota)$ be an
  admissible sequence of integral manifolds of $d$ such that $K =
  \lim_\iota V_\iota$; we may assume that $(V_\iota)$ has core $W$.
  Then for every $r>0$, the family of sets $\big\{\lim_\iota (V_\iota
  \cap N^{r,\epsilon}):\ \epsilon > 0\big\}$ is decreasing in
  $\epsilon$, where $\mathbf{N}^{r,\epsilon}:= \set{x \in M:\
    (x,r,\epsilon) \in \mathbf{N}}$, so we have $$\phi^{-1}(r) \cap K
  = \lim_{\epsilon \to 0} \lim_\iota (V_\iota \cap
  \mathbf{N}^{r,\epsilon}).$$ Hence, after passing to a subsequence of
  $(V_\iota)$ if necessary, there are $r_\iota \to 0$ and
  $\epsilon_\iota \to 0$ such that
  \begin{equation*}
    \fr(K \cap C) = \lim_\iota \left(V_\iota \cap
      \mathbf{N}^{r_\iota,\epsilon_\iota}\right) = \lim_\iota
    \Pi_n\big((V_\iota \times \{r_\iota,\epsilon_\iota)\}) \cap
    \mathbf{N}\big). 
  \end{equation*}
  Since $\lim_\iota(r_\iota,\epsilon_\iota) = (0,0)$, the right-hand
  side in the previous equality is equal to $\Pi_n \big( \lim_\iota
  \big((V_\iota \times \{(r_\iota,\epsilon_\iota)\}) \cap
  \mathbf{N}\big)\big)$.  Since the sequence $\big(V_\iota \times
  \{(r_\iota,\epsilon_\iota)\}\big)$ is an admissible sequence of
  integral manifolds of $\mathbf{d}$ with core $\mathbf{W}$, the lemma
  now follows from Remark \ref{limit_rmks} with $\mathbf{d}$ and
  $\mathbf{W}$ in place of $d$ and $W$.
\end{proof}

One problem with the previous lemma is that $\dim\mathbf{d} = \dim d$,
so it is possible that $\dim K_p > \dim \fr(K \cap C)$ for some $p$.
To remedy this, we need a fiber cutting lemma for pfaffian limits over
$\R$.

\begin{df} 
  \label{proper_NPL} 
  Let $K \subseteq \RR^n$ be a pfaffian limit obtained from $d$.
  We say that $K$ is \textbf{proper} if $\dim K = \dim d$.
\end{df}

\begin{prop} 
  \label{NPL_cutting} 
  Let $K \subseteq \RR^n$ be a pfaffian limit obtained from $d$ and
  $\nu \le n$.  Then there are $q \in \NN$ and proper pfaffian limits
  $K_1, \dots, K_q \subseteq \RR^n$ over $\R$ such that $$\Pi_\nu(K) =
  \Pi_\nu(K_1) \cup \cdots \cup \Pi_\nu(K_q)$$ and $\deg K_p \le \deg
  d$ and $\dim K_p = \dim \Pi_\nu(K_p) \le \dim K$ for each $p$.
\end{prop}

The following remark is needed in the proof of this proposition and
already appeared as Remark 3.5 in \cite{lio-spe:an_strat}; we restate
it here for the convenience of the reader.

\begin{nrmk}
  \label{isolated_points}
  Let $S \subseteq \RR^k$ be definable in an o-minimal expansion $\S$
  of the real field and put $l:= \dim S$.  Then there is a set $Y
  \subseteq S$, definable in $\S$, such that $S \subseteq \cl Y$,
  and for every $x \in Y$ there is a strictly increasing $\lambda:\{1,
  \dots, l\} \into \{1, \dots, k\}$ such that $x$ is isolated in $S
  \cap \Pi_\lambda^{-1}(\Pi_\lambda(x))$.
\end{nrmk}

\begin{proof}[Proof of Proposition \ref{NPL_cutting}]
  We proceed by induction on $m = \dim M$; the case $m=0$ is trivial,
  so we assume that $m>0$ and that the proposition holds for lower
  values of $m$.  Let $(V_\iota)$ be an admissible sequence of
  integral manifolds of $d$ such that $K = \lim_\iota V_\iota$.
  Choosing a suitable $C^2$-cell decomposition of $M$ compatible with
  $d$, and using Remark \ref{limit_rmks} and the inductive hypothesis,
  we reduce to the case where $M$ is a definable $C^2$-cell such that
  for every $s \leq \nu$ and every strictly increasing map
  $\lambda:\{1, \dots, s\} \into \{1, \dots, \nu\}$, the rank of
  $\Pi^n_{\lambda} \rest{d_l(x)}$ is constant for $x \in M$ and $l=1,
  \dots, k$; for $l=k$ we denote this rank by $r_\lambda$.  Putting
  $\Delta(\lambda) := \set{d, (\ker dx_{\lambda(1)})^{M}, \dots, (\ker
    dx_{\lambda(s)})^{M}}$ as at the beginning of Section
  \ref{khovsection} and using the associated notation, this means that
  each $g \in \D_{\Delta(\lambda)}$ has dimension, and we let
  $d^\lambda := d^{\Delta(\lambda),M} = \big(d^\lambda_0, \dots,
  d^\lambda_{k+r_\lambda}\big)$ be the corresponding definable nested
  distribution on $M$ of dimension $m-k-r_\lambda$.  It follows from
  the rank theorem and the fact that admissible integral manifolds of
  $d$ are closed in $M$ that $V_\iota \cap (\Pi^n_\lambda)^{-1}(y)$ is
  a closed integral manifold of $d^\lambda_{k+r_\lambda}$, for $\iota
  \in \NN$ and $y \in \Pi^n_\lambda(V_\iota)$.

  Let $s:= \dim \Pi_\nu(K)$; then $s \le \dim d$ by Lemma
  \ref{NPL_dimension}.  If $s = \dim d$, we are done, so we assume
  from now on that $s < \dim d$.  Let $\lambda:\{1, \dots, s\} \into
  \{1, \dots, \nu\}$ be strictly increasing; since $s<\dim d$, we
  have $$\dim d^\lambda \ge \dim d-s > 0;$$ in particular, $r_\lambda
  < \dim d$.  Hence by Lemma \ref{lowering dimension} and because each
  fiber $V_\iota \cap (\Pi^n_\lambda)^{-1}(y)$ is a closed submanifold
  of $M$, there is a closed, definable set $B_\lambda \subseteq M$
  such that $\dim B_\lambda < m$ and for $y \in \RR^s$ and $\iota \in
  \NN$, each component of the fiber $V_\iota \cap
  (\Pi^n_\lambda)^{-1}(y)$ intersects the fiber $B_\lambda \cap
  (\Pi^n_\lambda)^{-1}(y)$.  

  In particular, $\Pi^n_\lambda(V_\iota \cap B_\lambda) =
  \Pi^n_\lambda(V_\iota)$ for all $\iota$, and for all $y \in \RR^s$,
  every component of $\Pi_\nu(V_\iota) \cap (\Pi^\nu_\lambda)^{-1}(y)$
  intersects the fiber $\Pi_\nu(V_\iota \cap B_\lambda) \cap
  (\Pi^\nu_\lambda)^{-1}(y)$.

  We now denote by $\Lambda$ the set of all strictly increasing
  $\lambda:\{1, \dots, s\} \into \{1, \dots, \nu\}$.  Passing to a
  subsequence if necessary, we may assume for $\lambda \in \Lambda$
  that the sequence $(V_\iota \cap B_\lambda)_\iota$ converges to a
  compact set $K^\lambda$.  Choosing a suitable $C^2$-cell
  decomposition of $B_\lambda$ and using again Remark
  \ref{limit_rmks}, we see from the inductive hypothesis that the
  proposition holds with each $K^\lambda$ in place of $K$.  It
  therefore remains to show that $\Pi_\nu(K) = \bigcup_{\lambda \in
    \Lambda} \Pi_\nu\left(K^\lambda\right)$.  To see this, we fix a
  $\lambda \in \Lambda$; since each $\Pi_\nu\left(K^\lambda\right)$ is
  closed, it suffices by Remark \ref{isolated_points} to establish the
  following \medskip
  \\
  {\bf Claim.} Let $y \in \Pi^n_\lambda(K)$, and let $x \in \Pi_\nu(K)
  \cap (\Pi^\nu_\lambda)^{-1}(y)$ be isolated.  Then $x \in
  \Pi_\nu\left(K^\lambda\right)$. \medskip
  \\
  To see this, note that $\Pi_\nu(K) = \lim_\iota \Pi_\nu(V_\iota)$
  since $M$ is bounded.  Let $x_\iota \in \Pi_\nu(V_\iota)$ be such
  that $\lim_\iota x_\iota = x$, and put $y_\iota:=
  \Pi^\nu_\lambda(x_\iota)$.  Let $C_\iota \subseteq \RR^\nu$ be the
  component of $\Pi_\nu(V_\iota) \cap
  \left(\Pi^\nu_\lambda\right)^{-1}(y_\iota)$ containing $x_\iota$,
  and let $x_\iota'$ belong to $C_\iota \cap \Pi_\nu(V_\iota \cap
  B_\lambda)$.  Since also $\Pi_\nu\left(K^\lambda\right) = \lim_\iota
  \Pi_\nu(V_\iota \cap B_\lambda)$, we may assume, after passing to a
  subsequence if necessary, that $x' := \lim_\iota x_\iota' \in
  \Pi_\nu\left(K^\lambda\right)$.  We show that $x' = x$, which then
  proves the claim.  Assume for a contradiction that $x' \neq x$, and
  let $\delta > 0$ be such that $\delta \leq |x-x'|$ and
  \begin{equation} 
    \label{single} B(x,\delta) \cap \Pi_\nu(K) \cap
    (\Pi^\nu_\lambda)^{-1}(y) = \{x\}.
  \end{equation}
  Then for all sufficiently large $\iota$, there is an $x''_\iota \in
  C_\iota$ such that $\delta / 3 \leq |x''_\iota-x_\iota| \leq
  2\delta /3$, because $x_\iota, x_\iota' \in C_\iota$ and $C_\iota$
  is connected.  Passing to a subsequence if necessary, we may assume
  that $x'':= \lim_\iota x''_\iota \in \Pi_\nu(K)$.  Then $x'' \in
  B(x,\delta)$ with $x'' \neq x$, and since $x''_\iota \in C_\iota$
  implies that $\Pi^\nu_\lambda(x''_\iota) = y_\iota$, we get
  $\Pi^\nu_\lambda(x'') = y$, contradicting \eqref{single}.
\end{proof}

Combining Lemma \ref{frontier_of_NPL} with Proposition
\ref{NPL_cutting} gives:

\begin{cor}
  \label{frontier_of_pfaffian_limit}
  Let $K$ be a pfaffian limit obtained from $d$, and let $C \subseteq
  \RR^n$ be a definable cell.  Then there are $q \in \NN$ and proper
  pfaffian limits $K_1, \dots, K_q \subseteq \RR^{n+2}$ over $\R$ such
  that
  $$\fr(K \cap C) = \Pi_n(K_1) \cup \cdots \cup \Pi_n(K_q),$$ and $\deg K_p
  \le \deg d$ and $\dim K_p < \dim(K \cap C)$ for each $p$. \qed
\end{cor}

\section{Blowing-up along a nested distribution}  \label{jet}

In this section, we establish a criterion for generic portions (in the
sense of dimension and degree) of pfaffian limits over $\R$ to be
integral manifolds of definable nested distributions.  We fix a
bounded, definable manifold $M \subseteq \RR^n$ of dimension $m$ and a
definable nested distribution $d = (d_0, \dots, d_k)$ on $M$, and we
assume both are of class $C^2$.

\begin{df}  
  \label{blow-up}
  Put $n_1:= n+n^2$ and let $\Pi:\RR^{n_1} \into \RR^n$ denote the
  projection on the first $n$ coordinates.  We define
  \begin{align*}
    M^1 &:= \gr d_k \subseteq M \times G^{m-k}_n \subseteq \RR^{n_1},
    \text{ the graph of the distribution } d_k, \\ d^1_l &:= (\Pi
    \rest{M^1})^* d_l, \text{ the pull-back to } M^1 \text{ of } d_l
    \text{ via } \Pi, \text{ for } l=0, \dots, k.
  \end{align*}
  We call $d^1:= (d^1_0, \dots, d^1_k)$ the \textbf{blowing-up of $d$
    (along $d_k$)}; note that $M^1$ is of class $C^2$, while $d^1$ is
  of class $C^1$.  Finally, for $l = 0, \dots, k$ and an
  integral manifold $V$ of $\,d_l$, we define
  \begin{equation*}
    V^1 := (\Pi\rest{M^1})^{-1}(V),
  \end{equation*}
  the \textbf{lifting of $V$ (along $d_k$)}.  Note that, in this
  situation, $V^1$ is an integral manifold of $d_l^1$, and if $l=k$,
  then $V^1$ is also the graph of the Gauss map $g_V$.
\end{df}

Next, we write $M = \bigcup M_{\sigma}$, where $\sigma$ ranges over
$\Sigma_{n}$ and the $M_\sigma := M_{\sigma,2}$ are as before Lemma
\ref{eta-bdd-cover} with $d$ and $\eta$ there equal to $d_k$ and $2$
here.

\begin{df}  
  \label{frontier_blow-up}
  For an integral manifold $V$ of $\,d_k$ and $\sigma \in \Sigma_n$, we
  put $V_\sigma:= V \cap M_{\sigma}$.  Then $V_\sigma$ is an integral
  manifold of $d_k$, and we define
  \begin{equation*}
    F^1 V := \bigcup_{\sigma \in \Sigma_{n}} \fr\, V_\sigma^1.
  \end{equation*}
\end{df}

For our criterion, we let $D \subseteq \cl M^{1}$ be a definable
$C^2$-cell such that $C := \Pi(D)$ has the same dimension as $D$ and
$C$ is compatible with $M_\sigma$ and $\fr M_\sigma$ for every $\sigma
\in \Sigma_n$.  Then $D = \gr g$, where $g:C \into G^{m-k}_{n}$ is a
definable map, and we assume that the following hold:
\begin{renumerate}
\item the map $g \cap g_C$ has dimension and hence is a distribution
  on $C$;
\item if $g = g \cap g_C$, then either $g$ is integrable or $g$ is
  nowhere integrable.
\end{renumerate}  
We also assume that there is a definable set $W \subseteq \cl M^{1}$
such that $W \cap D = \emptyset$ and both $W$ and $W \cup D$ are open
in $\cl M^{1}$.  In this situation, for any sequence $(V_\iota)$ of
integral manifolds of $d_k$ such that $K := \lim_\iota V^{1}_\iota$
and $K' := \lim_\iota F^{1} V_\iota$ exist, we put
\begin{equation*}
  L_{(V_\iota)} := \left(D \cap K\right)
  \setminus \left(K'
    \cup \fr\left(W \cap K\right)\right).
\end{equation*}

\begin{rmk}
  Assume that $(V_\iota)$ is an admissible sequence of integral
  manifolds of $d$ such that $K := \lim_\iota V^1_\iota$ and $K' :=
  \lim_\iota F^1V_\iota$ exist, and assume that $K$ is proper.  Then
  $L_{(V_\iota)}$ is a generic subset of $K$ in the following sense:
  $\big(V^1_\iota\big)$ is an admissible sequence of integral
  manifolds of $d^1$, and each $(V_\iota)_\sigma^1$ is a finite union
  of admissible integral manifolds of $(d^1)^{M_\sigma}$.  By Lemma
  \ref{NPL_dimension}, Proposition \ref{frontier_lemma} and Remark
  \ref{limit_rmks}, $K'$ is a finite union of pfaffian limits over
  $\R$ of dimension less than $\dim K$ and degree at most $\deg d$.
  Moreover, by cell decomposition and Corollary
  \ref{frontier_of_pfaffian_limit}, there is a finite union $F
  \subseteq \RR^{n_1+2}$ of pfaffian limits over $\R$ of dimension
  less than $\dim K$ and degree at most $\deg d$ such that $\fr(W \cap
  K) = \Pi_{n_1}(F)$.
\end{rmk}

Finally, we let $g^1:D \into G^{m-k}_{n_{1}}$ be the pull-back of $g
\cap g_C$ to $D$ via $\Pi\rest{D}$.

\begin{prop}
  \label{local_2}
  In this situation, exactly one of the following holds:
  \begin{enumerate}
  \item $L_{(V_\iota)} = \emptyset$ for every admissible sequence
    $(V_\iota)$ of integral manifolds of $d$ such that $\lim_\iota
    V^{1}_\iota$ and $\lim_\iota F^{1} V_\iota$ exist;
  \item $g$ is an integrable distribution on $C$, and for every
    admissible sequence $(V_\iota)$ of integral manifolds of $d$
    such that $\,\lim_\iota V^{1}_\iota$ and $\lim_\iota
    F^{1} V_\iota$ exist, the set $L_{(V_\iota)}$ is an embedded
    integral manifold of $g^1$ and an open subset of $\,\lim_\iota
    V^{1}_\iota$.
  \end{enumerate}
  In particular, if $D$ is an open subset of $M^{1}$ and $(V_\iota)$
  is an admissible sequence of integral manifolds of $d$ such that
  $\,\lim_\iota V^{1}_\iota$ exists, then $D \cap \lim_\iota
  V_\iota^{1}$ is a finite union of leaves of $\big(d^{1}_k\big)^{D}$.
\end{prop}

We need the following observation for the proof of this proposition:

\begin{rmk}
  Let $\sigma \in \Sigma_{n}$.  Then $\sigma$ induces a diffeomorphism
  $\sigma:G_n \into G_n$ defined, in the notation of our conventions,
  by $\sigma(y) := A_{\sigma(\ker y)}$, and we define $\sigma^1:
  \RR^{n} \times G_{n} \into \RR^{n} \times G_{n}$ by $\sigma^1(x,y)
  := (\sigma(x), \sigma(y))$.  Note that $\sigma^1$ is also a
  permutation of coordinates.  The map $g^\sigma : \sigma(C) \into
  G^{m-k}_{n}$ defined by $g^\sigma(\sigma(x)) := \sigma(g(x))$
  satisfies $(g^\sigma)^1 = \sigma^1 \circ g^1 \circ (\sigma^1)^{-1}$.
  Moreover, if $(V_\iota)$ is a sequence of integral manifolds of
  $d_k$ such that $\lim_\iota V^{1}_\iota$ exists, then $\lim_\iota
  \sigma\left(V^{1}_\iota\right)$ also exists and $ \sigma^1(D) \cap
  \lim_\iota \sigma^1\left(V^{1}_\iota\right) = \sigma^1\left(D \cap
    \lim_\iota V^{1}_\iota\right)$.
\end{rmk}

\begin{proof}[Proof of Proposition \ref{local_2}]
  By the previous remark and Remark \ref{admissible_finiteness}, after
  replacing $M$ by $\sigma(M_{\sigma})$ and $W$ by $\sigma^1(W \cap
  \cl M_\sigma^1)$ for every $\sigma \in \Sigma_{n}$ satisfying $C
  \subseteq \cl M_\sigma$, we may assume for the rest of this proof
  that $d_k$ is $2$-bounded and prove the proposition with $\fr
  V_\iota^1$ in place of $F^1(V_\iota)$.  Thus, we let $(V_\iota)$ be
  an admissible sequence of integral manifolds of $d$ such that $K :=
  \lim_\iota V^{1}_\iota$ and $K' := \lim_\iota \fr V^{1}_\iota$
  exist, and we put
  \begin{equation*}
    L := (D \cap K)
    \setminus \big(K'
      \cup \fr(W \cap K)\big).
  \end{equation*}
  For the remainder of this proof, we simply write ``$\lim$'' in place
  of ``$\lim_\iota$''.  By definition of admissible sequence of
  integral manifolds and Corollary \ref{khovprop}(2), there is a $\nu
  \in \NN$ such that for every open box $U \subseteq \RR^n$ and every
  $\iota$, the set $U \cap V_\iota$ has at most $\nu$ components.  We
  assume that $L \ne \emptyset$; we need to show that $g$ is an
  integrable distribution on $C$ and that $L$ is an embedded integral
  manifold of $g^1$ and an open subset of $K$.  

  To do so, we choose an arbitrary $(x,y) \in L$ with $x \in \RR^{n}$
  and $y \in G_n$.  Since $W \cup D$ is open in $\cl M^{1}$, there is
  a bounded open box $B \subseteq\RR^{n_{1}}$ such that $(x,y) \in B$
  and
  \begin{equation*}
    \cl B \cap K \subseteq D \setminus
    \big(K' \cup\fr(W \cap K)\big); 
  \end{equation*}  
  in particular, $(x,y) \in D$.  We write $B = B_0 \times B_1$ with
  $B_0 \subseteq \RR^{n}$ and $B_1 \subseteq \RR^{n^2}$.  Since $D$ is
  the graph of the continuous map $g$ and $C$ is locally closed, we
  may also assume, after shrinking $B_0$ if necessary, that $D \cap
  \big(\cl B_0 \times \fr B_1\big) = \emptyset$.

  On the other hand, after passing to a subsequence if necessary, we
  may assume that $\lim \left(B \cap V^{1}_\iota\right)$, $\lim
  V_{\iota,B}^1$ and $\lim \fr V_{\iota,B}$ exist, where $V_{\iota,B}
  := \{x \in V_\iota:\ (x,T_x V_\iota) \in B\}$.  Then $$B \cap K =
  B \cap \lim \left(B \cap V^{1}_\iota\right) = B \cap \lim
  V_{\iota,B}^1.$$ We now claim that $x \notin \lim \fr V_{\iota,B}$:
  in fact, since $\fr V_\iota^1 \cap \cl B = \emptyset$ for all
  sufficiently large $\iota$, we have $\fr V_{\iota,B}^1 \subseteq \fr
  B$ for all sufficiently large $\iota$.  Also, $\lim V_{\iota,B}^1
  \subseteq \cl B \cap \lim V_\iota^1$ is disjoint from $\cl B_0
  \times \fr B_1$ by the previous paragraph, so $\cl V_{\iota,B}^1$ is
  disjoint from $\cl B_0 \times \fr B_1$ for all sufficiently large
  $\iota$.  Hence $\fr V_{\iota,B}^1 \subseteq \fr B_0 \times B_1$ for
  all sufficiently large $\iota$.  Since $B$ is bounded, we also have
  that $\fr V_{\iota,B} \subseteq \Pi_n\left(\fr
    V_{\iota,B}^1\right)$, and it follows that $\fr V_{\iota,B}
  \subseteq \fr B_0$ for all sufficiently large $\iota$, which proves
  the claim.

  Since each $V_\iota$ is a closed submanifold of $M$, we now apply
  Lemma \ref{local} with $V_{\iota,B}$ in place of $V_\iota$ and $\eta
  = 2$, to obtain a corresponding open neighbourhood $U \subseteq B_0$
  of $x$ and $f_{1}, \dots, f_{\nu}:\Pi_{m-k}(U) \into \RR^{n-m+k}$.
  We let $\lambda \in \{1, \dots, \nu\}$ be such that $x \in \gr
  f_{\lambda}$.  We claim that for every $x' \in \gr f_\lambda \cap
  U$, the map $f_{\lambda}$ is differentiable at $z' := \Pi_{m-k}(x')$
  with $T_{x'} \gr f_{\lambda} = g(x')$; since $x'$ is arbitrary, this
  claim implies that $\gr f_{\lambda}$ is an embedded, connected
  integral manifold of $g$.  Assumption (ii) and Remark
  \ref{integrable_rmks}(6) then imply that $g$ is an integrable
  distribution on $C$.  Since $(x,y) \in L$ was arbitrary, it follows
  that $L$ is an embedded integral manifold of $g^1$, as desired.
  
  To prove the claim, let $f_{\lambda,\iota}:\Pi_{m-k}(U) \into
  \RR^{n-m+k}$ be the functions corresponding to $f_{\lambda}$ as in
  the proof of Lemma \ref{local}.  After a linear change of
  coordinates if necessary, we may assume that $g(x') = \RR^{m-k}
  \times \{0\}$ (the subspace spanned by the first $m-k$ coordinates).
  It now suffices to show that $f_{\lambda}$ is $\eta$-Lipschitz at
  $x'$ for every $\eta > 0$, since then $T_{x'} \gr f_{\lambda} =
  \RR^{m-k} \times \{0\}$.  So let $\eta > 0$; since $\lim
  V^1_{\iota,B} \subseteq D = \gr g$ and $x' \in C$, and since $C$
  is locally closed and $g$ is continuous, there is a neighborhood $U'
  \subseteq U$ of $x'$ such that $\gr f_{\lambda,\iota} \cap U'$ is
  $\frac{\eta}{m-k}$-bounded for all sufficiently large $\iota$.
  Thus by Lemma \ref{local} again, $f_{\lambda}$ is $\eta$-Lipschitz
  at $x'$, as required.

  Finally, if $D$ is open in $M^{1}$, then $g = d_k\rest{C}$ and we
  can take $W:= \emptyset$.  Since $C$ is open in $M$ and compatible
  with $M_\sigma$ and $\fr M_\sigma$ for $\sigma \in \Sigma_n$, we have
  $C \cap \fr M_{\sigma} = \emptyset$ for each $\sigma$.  Hence
  $F^{1}(V_\iota) \cap D = \emptyset$, and it follows that
  $L_{(V_\iota)} = D \cap \lim V_\iota$ in this case.
\end{proof}

\section{Lifting a distribution to a nested
  distribution} \label{liftsection}

In the proofs of Propositions \ref{recover_1} and \ref{Rolle} below,
we will encounter individual distributions, such as the distribution
$g$ in the previous section, that need to be lifted to a given nested
distribution to produce a new nested distribution of lower dimension.
More precisely, we will encounter the following situation: we are
given $p \in \NN \cup \{\infty,\omega\}$ such that $p \ge 2$, a
definable $C^p$-cell $N \subseteq \RR^m$ with $m \ge n$ and a
definable, integrable nested $C^p$-distribution $f = (f_0, \dots,
f_l)$ on $N$.  We are also given a definable $C^p$-cell $D \subseteq
\RR^n$ such that $\Pi_n(N) \subseteq D$, a $k \le \dim D$ and a
definable, integrable distribution $h:D \into G^k_n$ on $D$ of class
$C^p$.  We let $\nu \le n$, and we assume that for all $\mu \le \nu$
and all strictly increasing $\lambda:\{1, \dots, \mu\} \into \{1,
\dots, \nu\}$ the dimension of the spaces $$F_\lambda(y):=
\Pi_\lambda(f_l(y)) \quad\text{and}\quad F^h_\lambda(y):=
\Pi_\lambda\big(\Pi_\nu(f_l(y)) \cap \Pi_\nu(h(\Pi_n(y)))\big)$$ is
constant as $y$ ranges over $N$; we denote these dimensions below by
$\dim F_\lambda$ and $\dim F^h_\lambda$, respectively.  For the
identity map $\lambda:\{1, \dots, \nu\} \into \{1, \dots, \nu\}$, we
set $F:= F_\lambda$ and $F^h := F^h_\lambda$ and put $\mu:= \dim F^h
\le k$.  Finally, we assume that $\dim F \ge \mu+1$.

\begin{rmkdf}
  \label{lift} 
  The assumptions in the previous paragraph imply that there is a
  strictly increasing $\lambda:\{1, \dots, \mu+1\} \into \{1, \dots,
  \nu\}$ such that $\dim F_\lambda = \mu+1$ and
  $\dim\left(F^h_\lambda\right) = \mu$.  We let $f_{l+1}:N \into G_m$ be
  the map on $N$ defined by
  \begin{equation*}
    f_{l+1}(y):= f_l(y) \cap \left(\Pi^m_\lambda\right)^{-1}
    \left(F^h_\lambda(y)\right).
  \end{equation*}
  The map $f_{l+1}$ is definable, $\dim f_{l+1}(y) = \dim f_l(y) - 1$
  and $\dim \Pi_\nu(f_{l+1}(y)) = \dim F - 1$ for $y \in N$; in
  particular, $f' := (f_0, \dots, f_{l+1})$ is a definable nested
  distribution on $N$.

  Next, let $Z \subseteq N$ be an embedded integral manifold of $f_l$,
  and let $L \subseteq D$ be an embedded integral manifold of $h$.  By
  our assumption on $f_l$, $\Pi_\nu \rest{Z}$ has constant rank
  $\dim F$; we assume here in addition that $\Pi_\nu(Z)$ is a
  submanifold of $\RR^\nu$.  Similarly, $\Pi_\nu \rest{L}$ has
  constant rank, and we assume here in addition that $\Pi_\nu(L)$ is a
  submanifold of $\RR^\nu$.  Then by our assumptions and the Rank
  Theorem, $\Pi_\nu(Z) \cap \Pi_\nu(L)$ is either empty or a
  submanifold of $\Pi_\nu(D)$ of dimension $\mu$.  Let also $L'
  \subseteq \Pi_\nu(Z) \cap \Pi_\nu(L)$ be a submanifold of dimension
  $\mu$.  Again by our assumption on $F^h_\lambda$, $\Pi^\nu_\lambda
  \rest{L'}$ is an immersion; we also assume here that
  $\Pi^\nu_\lambda(L')$ is a submanifold of $\RR^{\mu+1}$.  In this
  situation, we define
  \begin{equation*}
    Z(L') := Z \cap \left(\Pi^m_\lambda\right)^{-1}
    \left(\Pi^\nu_\lambda(L')\right).
  \end{equation*}
  Then by the Rank Theorem, the set $Z(L')$ is an integral manifold of
  $f_{l+1}$.
\end{rmkdf}

\begin{lemma}
  The nested distribution $f'$ is integrable.
\end{lemma}

\begin{proof}
  The integrability of $f_l$ and $h$ and our assumptions imply that
  for $y \in N$, there are an integral manifold $Z$ of $f_l$
  containing $y$ and an integral manifold $L$ of $h$ containing
  $\Pi_n(y)$.  By the Rank Theorem, after shrinking $Z$ and $L$ if
  necessary, we may assume that $\Pi_\nu(Z)$ and $\Pi_\nu(L)$ are
  submanifolds of $\RR^\nu$, so that $L':= \Pi_\nu(Z) \cap \Pi_\nu(L)$
  is a nonempty submanifold of $\RR^\nu$.  Shrinking $Z$ and $L$ again
  if necessary, we may also assume by the Rank Theorem that
  $\Pi_\lambda(L')$ is a submanifold of $\RR^{\mu+1}$.  So by Remark
  \ref{lift}, the corresponding $Z(L')$ is an integral manifold of
  $f_{l+1}$ containing $y$.
\end{proof}

For the next proposition, we let $\S$ be an o-minimal
expansion of $\R$.

\begin{prop}  
  \label{lifting}
  Let $Z$ be an integral manifold of $f_l$ and $L$ be an integral
  manifold of $h$, and assume that both $Z$ and $L$ are definable in
  $\S$ and $\Pi_n(Z) \cap L \ne \emptyset$.  Then there are integral
  manifolds $Z'_1, \dots, Z'_q$ of $f_{l+1}$ contained in $Z$ and
  definable in $\S$ such that $\Pi_n(Z) \cap L \subseteq \Pi_n(Z'_1)
  \cup \cdots \cup \Pi_n(Z'_q)$ and $\dim \Pi_\nu(Z'_p) < \dim
  \Pi_\nu(Z)$ for each $p$.
\end{prop}

\begin{proof}
  By Lemma \ref{o-min_lemma} with $\S$ in place of $\R$, we may assume
  that $\Pi_\nu(Z)$ and $\Pi_\nu(L)$ are submanifolds of $\RR^\nu$.
  By the above and Lemma \ref{o-min_lemma}, we have $\Pi_\nu(Z) \cap
  \Pi_\nu(L) = L'_1 \cup \cdots \cup L'_q$, where each $L'_p$ is an
  open subset of $\Pi_\nu(Z) \cap \Pi_\nu(L)$ such that
  $\Pi_\lambda(L'_p)$ is a submanifold of $\RR^{\mu+1}$.  Now we take
  $Z'_p := Z(L'_p)$, and we claim that these $Z'_p$ work.  To see
  this, we let $x \in \Pi_n(Z) \cap L$ and let $p \in \{1, \dots q\}$
  be such that $\Pi_\nu(x) \in L'_p$; we show that $x \in
  \Pi_n(Z'_p)$.  Since $Z \cap (\Pi^m_\lambda)^{-1}(\Pi_\lambda(x))
  \subseteq Z'_p$ by definition, and since $(\Pi^m_n)^{-1}(x)
  \subseteq (\Pi^m_\lambda)^{-1}(\Pi_\lambda(x))$, it follows that $Z
  \cap (\Pi^m_n)^{-1}(x) \subseteq Z'_p$, that is, $x \in
  \Pi_n(Z'_p)$.  The proposition now follows from Remark and
  Definition \ref{lift}.
\end{proof}

\section{Rewriting pfaffian limits}  \label{rewrite}

The goal of this section is to prove Proposition 1 of the introduction
using Corollary \ref{frontier_of_pfaffian_limit} and

\begin{prop}
  \label{recover_1} 
  Let $K \subseteq \RR^n$ be a pfaffian limit over $\R$.  Then there
  is a $q \in \NN$ and, for $p=1, \dots, q$, there are $n_p \ge n$ and
  an integral manifold $U_p \subseteq \RR^{n_p}$ over $\R$ such that
  $K \subseteq \Pi_n(U_1) \cup \cdots \cup \Pi_n(U_q)$ and for each
  $p$, the set $U_p$ is definable in $\P(\R)$ and $\dim \Pi_n(U_p)
  \le \dim K$.
\end{prop}

\begin{proof}
  Let $M \subseteq \RR^n$ be a definable $C^2$-manifold of dimension
  $m$, $d = (d_0, \dots, d_k)$ a nested distribution on $M$ and
  $(V_\iota)$ an admissible sequence of integral manifolds of $d$ such
  that $K = \lim_\iota V_\iota$.  We proceed by induction on the pair
  $(\deg d,\dim K)$, where we consider $\NN^2$ with its lexicographic
  ordering.  If $\deg d = 0$, then $K$ is definable in $\R$ by
  \cite[Theorem 3.1]{vdd:hausdorff} or \cite[Theorem
  1]{lio-spe:hausdorff}, so the proposition follows from cell
  decomposition and Example 2.  If $\dim K = 0$, then $K$ is finite
  and the proposition follows again from Example 2.  Therefore, we
  assume that $\deg d>0$ and $\dim K >0$ and that the proposition
  holds for lower values of $(\deg d, \dim K)$.  Moreover, by
  Proposition \ref{NPL_cutting}, we may assume that $K$ is a proper
  pfaffian limit over $\R$, that is, $\dim K = \dim d = m-k$.  By
  $C^{m+2}$-cell decomposition, Remark \ref{limit_rmks} and the
  inductive hypothesis, we may assume that $M$ and $d$ are of class
  $C^{m+2}$.  If $\sim_{d_k}$ is definable, then we are done as in the
  case $\deg d = 0$.  We therefore assume from now on that
  $\sim_{d_k}$ is not definable.

  We now blow up $m+1$ times along $d_k$, that is, we put $n_0 := n$,
  $M^0 := M$ and $d^0 := d$, and we put $V^0 := V$ and $F^0 V :=
  \bigcup_{\sigma \in \Sigma_n} \fr V_\sigma$ for every integral
  manifold $V$ of $d_l$ with $l \in \{0, \dots, k\}$.  By induction on
  $j = 1, \dots, m+1$ we define $n_j := (n_{j-1})_1 = n_{j-1} +
  n_{j-1}^2$, $M^j := \left(M^{j-1}\right)^1 = \gr d^{j-1}_k$ and $d^j
  := \left(d^{j-1}\right)^1$, and we define the corresponding liftings
  $V^j := \left(V^{j-1}\right)^1$ and $F^j V := F^1V^{j-1}$ for every
  integral manifold $V$ of $d_l$ with $l \in \{0, \dots, k\}$.  For $0
  \le i \le j \le m+1$, we also let $\pi^j_i:\RR^{n_j} \into
  \RR^{n_i}$ be the projection on the first $n_i$ coordinates.

  Passing to a subsequence if necessary, we may assume that $K^j:=
  \lim_\iota V^j_\iota$ and $\lim_\iota F^j V_\iota$ exist for $j=0,
  \dots, m$ (so $K^0 = K$).  Then $\pi^j_0\big(K^j\big) = K$ for each
  $j$, and since $K$ is proper, each $K^j$ is also proper.  It follows
  from Remark \ref{limit_rmks}, Corollary \ref{frontier_lemma} and
  the inductive hypothesis that
  \begin{itemize} 
  \item[(I)] the proposition holds with each $\lim_\iota F^j V_\iota$
    in place of $K$.
  \end{itemize}
  For $j=0, \dots, m+1$, we set $M^j_\sigma := (M^j)_{\sigma,2}$ as in
  Lemma \ref{eta-bdd-cover} with $M$, $d$ and $\eta$ there equal to
  $M^j$, $d^j_k$ and $2$ here.  We let $\C^j$ be a $C^2$-cell
  decomposition of $\,\RR^{n_j}$ compatible with $M^j$, $\fr M^j$,
  $\set{M^j_\sigma:\ \sigma \in \Sigma_{n_j}}$ and $\set{\fr
    M^j_\sigma:\ \sigma \in \Sigma_{n_j}}$, and we put
  \begin{equation*}
    \C^j_M:= \set{C \in \C^j:\ C \subseteq \cl M^j}.
  \end{equation*}
  Refining each $\C^j$ in order of decreasing $j \in \{0, \dots, m\}$
  if necessary, we may assume for each such $j$ that
  \begin{renumerate}
  \item $\C^j$ is a stratification compatible with
    $\set{\pi^{j+1}_j(C):\ C \in \C^{j+1}}$;
  \end{renumerate}
  and for every $D \in \C^{j+1}_M$ that is the graph of a map
  $g:C \into G^{m-k}_{n_j}$, where $C:= \pi^{j+1}_j(D)$, that
  \begin{renumerate}
  \item[(ii)] the map $g \cap g_C$ has dimension and hence is a
    distribution on $C$;
  \item[(iii)] if $g = g \cap g_C$, then either $g$ is
    integrable or $g$ is nowhere integrable.
  \end{renumerate}
  By Corollary \ref{frontier_of_pfaffian_limit} and the inductive
  hypothesis,
  \begin{itemize} 
  \item[(II)] for $j = 0, \dots, m$ and $E \in \C^j_M$, the
    proposition holds with $\fr\left(K^j \cap E\right)$ in place of
    $K$.
  \end{itemize}

  We now fix $j \in \{0, \dots, m\}$ and a cell $C \in \C^j_M$ such
  that $\dim C \ge j$.  \medskip

  \noindent\textbf{Claim 1:} 
  There is a $q \in \NN$ and, for $p=1, \dots, q$, there are $m_p \ge
  n_j$ and an integral manifold $U_p \subseteq \RR^{m_p}$ over $\R$ such
  that $K^j \cap C \subseteq \Pi_{n_j}(U_1) \cup \cdots \cup
  \Pi_{n_j}(U_q)$ and for each $p$, the set $U_p$ is definable in
  $\P(\R)$ and $\dim \Pi_{n_j} (U_p) \le \dim K$.  \medskip

  Assuming Claim 1 holds, the proposition follows by applying Claim 1
  to each $C \in C^0_M$.  To prove Claim 1, we proceed by reverse
  induction on $\dim C \le m$.  Let
  \begin{equation*}
    \D_C:= \set{D' \cap \left(\pi^{j+1}_j\right)^{-1}(C):\ D' \in
      \C^{j+1}_M,\ C \subseteq \pi^{j+1}_j(D')},
  \end{equation*}
  and fix an arbitrary $D \in \D_C$; it suffices to prove Claim 1 with
  $K^{j+1}$ and $D$ in place of $K^j$ and $C$.  Let $D' \in
  \C^{j+1}_M$ be such that $D \subseteq D'$; if $\dim D' >\dim C$,
  then Claim 1 with $K^{j+1}$ and $D$ in place of $K^j$ and $C$
  follows from the inductive hypothesis, so we also assume that $\dim
  D' = \dim C$.  Then $D$ is open in $D'$, and since $M^{j+1}
  \subseteq \RR^{n_j} \times G_{n_j}^{m-k}$ and $G_{n_j}^{m-k}$ is
  compact, there is a definable map $g:C \into G^{m-k}_{n_j}$ such
  that $D = \gr g$.  Let
  \begin{equation*}
    W := \bigcup \set{E \in \C^{j+1}_M:\ \dim E > \dim C};
  \end{equation*} 
  since $\C^{j+1}$ is a stratification, both $W$ and $W \cup D'$ are
  open in $\cl M^{j+1}$, and since $D$ is open in $D'$, the set $W
  \cup D$ is also open in $\cl M^{j+1}$.  Hence by Proposition
  \ref{local_2}, the set $\left(K^{j+1} \cap D\right) \setminus
  \left(\lim_\iota F^{j+1} V_\iota \cup \fr\left(W \cap
      K^{j+1}\right)\right)$ is an embedded integral manifold of
  $g^1$, where $g_1 : D \into G^{m-k}_{n_{j+1}}$ is the pullback of
  $g$ via the restriction of $\pi^{j+1}_j$ to $D$.  But $D \cap
  \fr\left(W \cap K^{j+1}\right) \subseteq F$, where
  \begin{equation*}
    F := \bigcup \set{\fr\left(E \cap K^{j+1}\right):\ E \in
      \C^{j+1}_M \text{ and } \dim E > \dim C}
  \end{equation*} 
  is compact, so the set
  \begin{equation*}
    L:= \left(K^{j+1} \cap D\right) \setminus \left(\lim_\iota
      F^{j+1}(V_\iota) \cup F\right)
  \end{equation*}
  is a finite union of connected integral manifolds of $g^1$ definable
  in $\P(\R)$.  Thus by (I) and (II), to prove Claim 1 with $K^{j+1}$
  and $D$ in place of $K^j$ and $C$, it now suffices to prove:
  \medskip

  \noindent\textbf{Claim 2:} There is a $q \in \NN$ and, for $p=1,
  \dots, q$, there are $m_p \ge n_{j+1}$ and an integral manifold $U_p
  \subseteq \RR^{m_p}$ over $\R$ such that $L \subseteq
  \Pi_{n_{j+1}}(U_1) \cup \cdots \cup \Pi_{n_{j+1}}(U_q)$ and for each
  $p$, the set $U_p$ is definable in $\P(\R)$ and $\dim \Pi_{n_{j+1}}
  (U_p) \le \dim L$. \medskip

  For the proof of Claim 2, we write $n$, $d$ and $V_\iota$ in place
  of $n_{j+1}$, $d^{j+1}$ and $V_\iota^{j+1}$.  Since $\sim_{d_k}$ is
  not definable, the nested distribution $d':= (d_0, \dots, d_{k-1})$
  satisfies $\deg d' < \deg d$.  Moreover, there is an admissible
  sequence $(U_\iota)$ of integral manifolds of $d'$ such that
  $V_\iota \subseteq U_\iota$ for all $\iota$.  Passing to a
  subsequence if necessary, we may assume that $Y:= \lim_\iota
  U_\iota$ exists; in particular, note that $L \subseteq Y$.  By the
  inductive hypothesis, there is a $q \in \NN$ and, for $p=1, \dots,
  q$, there exist $m_p \ge n$, a definable $C^2$-manifold $N_{p}
  \subseteq \RR^{m_p}$, a definable nested distribution $f_{p}=
  \big(f_{p,0}, \dots, f_{p,k(p)}\big)$ on $N_{p}$ and a nested
  integral manifold $Z_{p} = \big(Z_{p,0}, \dots, Z_{p,k(p)}\big)$ of
  $f_{p}$ definable in $\P(\R)$ such that $$Y \subseteq
  \Pi_{n}\big(Z_{1,k(1)}\big) \cup \cdots \cup
  \Pi_{n}\big(Z_{q,k(q)}\big)$$ and $\dim \Pi_n\big(Z_{p,k(p)}\big)
  \le \dim Y$ for each $p$.  After refining and pruning the collection
  $\big\{N_{1}, \dots, N_{q}\big\}$ if necessary, we may assume that
  that
  \begin{itemize}
  \item[($\dagger$)] $\Pi_{n}(N_{p}) \subseteq D$ and
    $\Pi_n(Z_{p,k(p)}) \cap L \ne \emptyset$ for each $p$, and that $L
    \subseteq \Pi_n(Z_{1,k(1)}) \cup \cdots \cup \Pi_n(Z_{q,k(q)})$.
  \end{itemize}
  We now prove Claim 2 by induction on $$\delta := \max\set{\dim
    \Pi_n \left(Z_{p,k(p)}\right):\ p = 1, \dots, q},$$
  simultaneously for all such collections $\set{(N_p, f_p, Z_p)}$
  satisfying ($\dagger$).  If $\delta = \dim L$, we take $U_p :=
  Z_{p,k(p)}$ for each $p$ and are done.  So assume $\delta > \dim L$
  and Claim 2 holds for lower values of $\delta$.  Refining and
  pruning the collection $\big\{N_{1}, \dots, N_{q}\big\}$ again if
  necessary, we may assume for each $p$, each $\mu \le n$ and each
  strictly increasing $\lambda:\{1, \dots, \mu\} \into \{1, \dots,
  n\}$ that the dimension of the spaces $$\Pi_\lambda(f_{p,k(p)}(y))
  \quad\text{and}\quad \Pi_\lambda\left(\Pi_n(f_{p,k(p)}(y)) \cap
    g^1(\Pi_n(y))\right)$$ is constant as $y$ ranges over $N_{p}$.
  Then by Proposition \ref{lifting} with $\nu = n$, for each $p$ such
  that $\dim \Pi_n \left(Z_{p,k(p)}\right) = \delta$, the integral
  manifold $Z_{p,k(p)}$ over $\R$ can be replaced by finitely many
  integral manifolds $Z'$ over $\R$ contained in $N_p$ such that $\dim
  \Pi_n(Z') < \delta$.  Claim 2 then follows from the inductive
  hypothesis, and the proposition is proved.
\end{proof}

\section{Normal sets and regular closure}  
\label{analytic}

We assume from now on that $\R$ admits analytic cell decomposition.
In this section, we denote by $\|\cdot\|$ the Euclidean norm on
$\RR^k$ (for $k = 1,2,\dots$).

\begin{df}  
  \label{normal}
  A definable open set $U \subseteq \RR^n$ is \textbf{normal} if there
  exists an analytic, definable carpeting function on $U$.
\end{df}

\begin{rmk}
  If $U,V \subseteq \RR^n$ are normal, then so are $U \cap V$ and
  $U \times V$.
\end{rmk}

\begin{expl}  
  \label{open_cells_are_normal}
  Every analytic, open definable cell is normal.
\end{expl}

\begin{prop}  
  \label{normal_covering}
  Let $A \subseteq \RR^n$ be open and definable.  Then $A$ is a
  union of finitely many normal sets.
\end{prop}

\begin{proof}
  By induction on $n$; we may assume that $A \ne \RR^n$.  The case
  $n=0$ is trivial, so we assume that $n>0$ and that the proposition
  holds for lower values of $n$.  By analytic cell decomposition, it
  suffices to show that every definable analytic cell contained in $A$
  is contained in a finite union of normal sets contained in $A$.

  So we let $C \subseteq A$ be a definable analytic cell; we proceed by
  induction on the dimension $c$ of $C$.  If $c=0$, then $C$ is a
  singleton and any ball centered at $C$ and contained in $A$ will do.
  So we assume that $c>0$ and that every analytic cell of dimension
  less than $c$ contained in $A$ is contained in a finite union of
  normal sets contained in $A$.  If $c=n$, then $C$ is open and hence
  normal by Example \ref{open_cells_are_normal}; so we also assume
  that $c<n$.

  After permuting coordinates if necessary, there is an open, analytic
  cell $D \subseteq \RR^c$ and a definable, analytic map $g:D \into
  \RR^{n-c}$ such that $C = \gr g$.  For $x \in \RR^n$, we write $x =
  (y,z)$ with $y \in \RR^c$ and $z \in \RR^{n-c}$.  Define $\beta:D
  \into (0,\infty)$ by $\beta(y):= \dist\big((y,g(y)), \RR^n \setminus
  A\big)$, where this distance is computed using the Euclidean norm.
  By the inductive hypothesis and analytic cell decomposition, we may
  assume that $\beta$ is analytic.  Now we put
  \begin{equation*}
    U:= \set{(y,z) \in A:\ y \in D \text{ and } \|z - g(y)\|^2 <
      \beta^2(y)}. 
  \end{equation*}
  This $U$ is normal: given an analytic, definable carpeting function
  $\gamma:D \into (0,\infty)$, we define $\phi:U \into (0,\infty)$ by
  \begin{equation*}
    \phi(y,z):= \gamma(y) \big(\beta^2(y) - \|z-g(y)\|^2\big),
  \end{equation*}
  which is easily seen to be an analytic, definable carpeting function on
  $U$.
\end{proof}

\begin{df}  
  \label{normal_in}
  Let $U \subseteq \RR^n$ be normal and $A \subseteq U$.  
  \begin{enumerate}
  \item $A$ is \textbf{normal in $U$} if $A$ is a finite union of sets
    of the form $$\set{x \in U:\ g(x) = 0,\ h(x) > 0},$$ where $g:U
    \into \RR^q$ and $h:U \into \RR^r$ are definable and analytic and
    ``$h(x) > 0$'' means ``$h_s(x) > 0$ for $s = 1, \dots, r$''.
  \item $A$ is a \textbf{normal leaflet in $U$} (\textbf{of
      codimension $p$}) if
    \begin{equation*}
      A = \set{x \in U:\ f(x) = g(x) = 0,\ h(x)>0},
    \end{equation*}
    where $f:U \into \RR^p$, $g:U \into \RR^q$ and $h:U \into \RR^r$
    are definable and analytic and for all $x \in A$, the rank of $f$
    at $x$ is $p$ and $\ker df(x) \subseteq \ker dg(x)$.
  \end{enumerate}
\end{df}

\begin{nrmk}
  \label{leaflet_rmk}
  In the situation of Definition \ref{normal_in}(2), the set $A$ is an
  analytic submanifold of $\,U$ of dimension $n-p\,$; in fact, we have
  $T_xA = \ker df(x)$ for all $x \in A$.  Moreover, if $\delta:U \into
  (0,\infty)$ is an analytic, definable carpeting function on $U$, then the
  restriction of $\phi:= \delta \cdot \prod_{s=1}^r H_s:U \into \RR$
  to $A$ is an analytic carpeting function on $A$, where $H_s
  := h_s/\sqrt{1+h_s^2}$ for each $s$.
\end{nrmk}

\begin{expl}
  Let $U \subseteq \RR^n$ be normal; then $U$ is a normal leaflet in
  $U$.  Let also $f:U \into \RR^p$ and $h:U \into \RR^r$ be analytic
  and definable.  Then the set
  \begin{multline*}
    \set{x \in U:\ f(x) = 0,\ h(x) > 0 \text{ and } f \text{ has rank
      } p \text{ at } x} \\ = \set{x \in U:\ f(x) = 0,\ h(x) > 0,\
      |df|^2(x) > 0} \\ = \set{x \in U:\ f(x) = 0,\ h(x) > 0,\ df_1(x)
      \wedge \cdots \wedge df_p(x) \ne 0}
  \end{multline*}
  is a normal leaflet in $U$, where $|df|^2$ denotes the sum of the
  squares of all $(p \times p)$-subdeterminants of $df$.
\end{expl}

The following lemma is elementary, and its proof is left to the
reader.

\begin{lemma}  
  \label{normal_closures}
  Let $U \subseteq \RR^n$ and $V \subseteq \RR^m$ be normal.
  \begin{enumerate}
  \item If $A$ and $B$ are normal in $U$, then so are $A \cup B$, $A
    \cap B$ and $A \setminus B$.
  \item Let $A$ be normal in $U$ and $B$ be normal in $V$.  Then $A
    \times B$ is normal in $U \times V$.  Moreover, if $A$ and $B$ are
    normal leaflets in $U$ and $V$, respectively, then $A \times B$ is
    a normal leaflet in $U \times V$.
  \item Let $\phi:U \into V$ be definable and analytic, and let $B$ be
    normal in $V$.  Then $\phi^{-1}(B)$ is normal in $U$.
  \item Let $A$ be a normal leaflet in $U$ of codimension $p$, and
    assume that there is a definable, analytic embedding $\phi:A \into
    \RR^{n-p}$.  Then $\phi(A)$ is normal.  Moreover, if $B \subseteq
    A$ is normal (resp., a normal leaflet) in $U$, then $\phi(B)$ is
    normal (resp., a normal leaflet) in $\phi(A)$. \qed
  \end{enumerate}
\end{lemma}

The following lemma will be used several times in this section:

\begin{lemma}  
  \label{pfaffian_puiseux}
  Let $\eta:X \into (0,\infty)$ be definable in $\P(\R)$ and put
  \begin{equation*}
    Y:= \set{(x,t):\ x \in X \text{ and } 0<t<\eta(x)}.
  \end{equation*}
  Let also $\alpha:Y \into [0,\infty)$ be definable in $\P(\R)$, and
  assume that for every $x \in X$, the function $\alpha_x:(0,\eta(x))
  \into [0,\infty)$ defined by $\alpha_x(t):= \alpha(x,t)$ is
  semianalytic.  Then there exists an $N \in \NN$ such that for all $x
  \in X$, either ultimately $\alpha_x(t) = 0$ or ultimately
  $\alpha_x(t) > t^N$ (where ``ultimately'' abbreviates ``for all
  sufficiently small $t>0$'').
\end{lemma}

\begin{proof}
  By monotonicity, for every $x \in X$ the function $\alpha_x$ is
  ultimately of constant sign.  By Puiseux's Theorem, for every $x \in
  X$ such that $\alpha_x$ is ultimately positive, there are $c_x > 0$
  and $r_x \in \QQ$ such that ultimately $\alpha_x(t) = c_x t^{r_x} +
  o(t^{r_x})$.  However, the set $R_X:= \set{r_x \in \RR:\ x \in X}$
  is definable, since for all $x \in X$ such that $\alpha_x$ is
  ultimately positive, we have $r_x = \lim_{t \to 0^+}
  t\alpha_x'(t)/\alpha_x(t)$.  Since each $r_x$ is rational, it
  follows that $R_X$ is finite, so any $N>\max R_X$ will do.
\end{proof}

\begin{prop}[Gabrielov \cite{gab:frontier}]  
  \label{leaflet_dec}
  Let $U \subseteq \RR^n$ be normal and $A$ be normal in $U$.  Then
  $A$ is a finite union of normal leaflets in $U$.
\end{prop}

\begin{proof}
  Let $g:U \into \RR^q$ and $h:U \into \RR^r$ be definable and
  analytic such that $A = \set{x \in U:\ g(x) = 0,\ h(x) > 0}$; we
  proceed by induction on $d:= \dim A$.  If $d=0$, the proposition is
  trivial, so we assume that $d>0$ and the proposition holds for lower
  values of $d$.

  Let $\C$ be a finite stratification of $A$ into analytic cells, and
  put $\C_d := \set{C \in \C:\ \dim C = d}$.  We show that for $C \in
  \C_d$, there is a normal leaflet $A_C \subseteq A$ such that $\dim(C
  \setminus A_C) < d$.  The proposition then follows from the
  inductive hypothesis, since $A \setminus \bigcup_{C \in \C_d} A_C$
  is a normal subset of $U$ of dimension less than $d$.

  Fix $C \in \C_d$, and let $\G$ be the set of all partial derivatives
  (of all orders) of $g_1, \dots, g_q$.  Let $M$ be the set of all
  natural numbers $m$ for which there exist $f_1, \dots, f_m \in \G$
  such that $C \subseteq \set{x \in U:\ f_1(x) = \cdots = f_m(x) = 0}$
  and $df_1(a) \wedge \cdots \wedge df_m(a) \ne 0$ for some $a \in C$.

  Put $p:= \sup M \le n-d$; we claim that $p = n-d$.  To see this, let
  $f_1, \dots, f_p \in \G$, $a \in C$ and an open ball $B$ centered at
  $a$ be such that $C \cap B$ is a connected submanifold of
  \begin{equation*}
    \Gamma:= \set{x \in U:\ f_1(x) = \cdots = f_p(x) = 0} \cap B,
  \end{equation*}
  and such that $\Gamma$ is a connected, analytic submanifold of
  codimension $p$ contained in $\set{x \in U:\ h(x) > 0,\ df_1(x)
    \wedge \cdots \wedge df_p(x) \ne 0}$.  The maximality of $p$ now
  implies that $g(x) = 0$ for all $x \in \Gamma$.  Since $\C$ is a
  stratification of $A$ and $\dim C = d$, the cell $C$ is an open
  subset of $A$; so after shrinking $B$ if necessary, we may assume
  that $C \cap B = A \cap B$.  It follows that $A \cap B \subseteq
  \Gamma$, that is, $C \cap B = A \cap B = \Gamma$, which proves the
  claim.

  Put $f:= 0$ if $p=0$ and $f:= (f_1, \dots, f_p)$ otherwise, where
  $f_1, \dots, f_p \in \G$ are as in the previous paragraph.  Let
  \begin{equation*}
    V:= \set{x \in U:\ f(x) = 0,\ h(x) > 0,\ |df|^2(x) > 0},
  \end{equation*}
  a normal leaflet in $U$.  Note that $C \setminus V$ is contained in
  the set $\{x \in U:\ f(x) = 0,\ |df|^2(x) = 0\}$, and the latter has
  dimension strictly less than $d$.  The only remaining problem,
  therefore, is that $V$ is not necessarily a subset of $A$.  To
  address this issue, we let $\eta:V \into (0,\infty)$ be a definable
  function such that $B(x,2\eta(x)) \subseteq \set{x \in U:\ h(x) >
    0,\ |df|^2(x) > 0}$ for all $x \in V$.  For $x \in V$ and $t \in
  (0,\eta(x))$, we put
  \begin{equation*}
    \alpha(x,t):= \max\set{|g|(y):\ y \in V,\ \|y-x\| \le t}.
  \end{equation*}
  Note that $\alpha$ is definable, and for $x \in V$ the function
  $\alpha_x:(0,\eta(x)) \into [0,\infty)$ is semianalytic.  Hence by
  Lemma \ref{pfaffian_puiseux}, there is an $N \in \NN$ such that for
  all $x \in V$, either ultimately $\alpha_x(t) = 0$ or ultimately
  $\alpha_x(t) > t^N$.

  We now let $A_C$ be the subset of $V$ where all partials of $g
  \rest{V}$ up to order $N$ vanish.  More precisely, we let $\S$ be
  the set of all $\phi:U \into \RR$ for which there exist $\nu \le N$
  and functions $\phi_0, \dots, \phi_{\nu}:U \into \RR$ such that
  $\phi = \phi_{\nu}$, $\phi_0 \in \{g_1, \dots, g_q\}$ and for
  $i \in \{1, \dots, \nu\}$, the function $\phi_{i}$ is one of the
  coefficient functions of $d\phi_{i-1} \wedge df_1 \wedge \cdots
  \wedge df_p$.  Then we put $$A_C:= \set{x \in V:\ \phi(x) = 0 \text{
      for all } \phi \in \S}.$$ 

  It follows for all $x \in V$ that $x \in A_C$ if and only if
  ultimately $\alpha_x(t) \le t^N$, that is, if and only if ultimately
  $\alpha_x(t) = 0$, that is, if and only if $g = 0$ in a neighbourhood
  of $x$ in $V$.  Hence $A_C \subseteq A$, and the only points $x \in V
  \cap A$ that are not contained in $A_C$ are those where $g$ is not
  identically $0$ on any neighbourhood of $x$ in $V$.  Thus, $\dim(C
  \setminus A_C) < d$, and since $A_C$ is a leaflet in $U$, the
  proposition is proved.
\end{proof}

\begin{cor}  
  \label{constant_rank}
  Let $U \subseteq \RR^n$ be normal and $\D$ be a finite collection of
  definable, analytic distributions on $U$, and let $A$ be a normal
  subset of $U$.  Then $A$ is a finite union of normal leaflets in $U$
  that are compatible with $\D$.
\end{cor}

\begin{proof}
  By induction on $\dim A$; if $\dim A = 0$, there is nothing to do,
  so we assume $\dim A > 0$ and the corollary is true for lower values
  of $\dim A$.  By Proposition \ref{leaflet_dec} and the inductive
  hypothesis, we may assume that $A$ is a normal leaflet in $U$.
  Proceeding as in the proof of Proposition \ref{partition}, we note
  that for $\E \subseteq \D$ and $i \in \{0, \dots, \dim A\}$, the set
  $A_{\E,i}:= \set{x \in A:\ \dim T_x\E = i}$ is a normal subset of
  $U$.  Again by Proposition \ref{leaflet_dec}, there is a finite
  collection $\C$ of normal leaflets in $U$ such that each $A_{\E,i}$
  is a union of leaflets in $\C$.  Arguing as in the proof of
  Proposition \ref{partition}, we see that every leaflet in $\C$ of
  dimension $\dim A$ is compatible with $\D$.  The corollary now
  follows from the inductive hypothesis.
\end{proof}

We next obtain a fiber cutting lemma for nested Rolle leaves
(Corollary \ref{fiber_cutting} below), using normal leaflets and
building on the techniques found in Moussu and Roche
\cite{mou-roc:khovtheory}, Lion and Rolin \cite{lio-rol:feuilletages}
and \cite{spe:pfaffian}.  We let $U \subseteq \RR^n$ be a normal set
and $A$ a normal subset of $U$.  We also let $\Delta = \{d^1, \dots,
d^q\}$ be a set of definable, analytic nested distributions on $U$; we
write $d^p = \big(d^p_0, \dots, d^p_{k(p)}\big)$ for $p=1, \dots, q$
and associate $\D_\Delta$ to $\Delta$ as in Section \ref{khovsection}.

\begin{lemma}  
  \label{leaflet_morse}
  Assume that $A$ is a normal leaflet in $U$ compatible with
  $\D_\Delta$, and suppose that $\dim d^{\Delta,A}_{k(\Delta,A)} > 0$.
  Then there is an analytic, definable $\phi:U \into \RR$ such that
  $\phi\rest{A}$ is a carpeting function on $A$ and the definable set
  \begin{equation*}
    B:= \set{x \in A:\ d^{\Delta,A}_{k(\Delta,A)}(x) \subseteq \ker
      d(\phi\rest{A})(x)} 
  \end{equation*}
  has dimension less than $\dim A$.
\end{lemma}

\begin{proof}
  Let $\psi:U \into \RR$ be analytic and definable such that
  $\psi\rest{A}$ is a carpeting function on $A$ (as obtained in Remark
  \ref{leaflet_rmk}, say).  For $u \in (0,\infty)^n$, we define
  $\psi_u:U \into \RR$ by $\psi_u(x):= \psi(x) \phi_u(x)$, where
  $\phi_u(x)$ is as in Section \ref{prelim}; note that
  $\psi_u\rest{A}$ is an analytic, definable carpeting function on
  $A$.  Now consider the definable set
  \begin{equation*}
    D := \set{(u,x) \in \RR^n \times A:\ d^{\Delta,A}_{k(\Delta,A)}(x)
      \subseteq \ker d(\psi_u\rest{A})(x)}.
  \end{equation*}
  Arguing as in the proof of Lemma \ref{lowering dimension}, we
  conclude that $\dim D_u < \dim A$ for some $u \in (0,\infty)^n$,
  so we take $\phi := \psi_u$.
\end{proof}

For $I \subseteq \{1, \dots, q\}$ we put $\Delta(I) := \{d^p:\ p \in
I\}$.

\begin{lemma}
  \label{cell_cutting}
  Let $I \subseteq \{1, \dots, q\}$.  Then there is a finite
  collection $\P$ of normal leaflets in $U$ contained in $A$ such that
  $\P$ is compatible with $\D_{\Delta(J)}$ for every $J \subseteq \{1,
  \dots, q\}$ and
  \begin{renumerate}
  \item $\dim d^{\Delta(I),N}_{k(\Delta(I),N)} = 0$ for every $N
    \in \P$;
  \item whenever $V_p$ is a Rolle leaf of $d^p$ for $p \in I$, every
    component of $A \cap \bigcap_{p \in I} V_p$ intersects some
    leaflet in $\P$.
  \end{renumerate}
\end{lemma}

\begin{proof}
  By induction on $\dim A$; if $\dim A = 0$, there is nothing to do,
  so we assume $\dim A > 0$ and the corollary is true for lower values
  of $\dim A$.  By Corollary \ref{constant_rank} and the inductive
  hypothesis, we may assume that $A$ is a normal leaflet in $U$
  compatible with $\D_{\Delta(J)}$ for $J \subseteq \{1, \dots, q\}$.
  Thus, if $\dim d^{\Delta(I),A}_{k(\Delta(I),A)} = 0$, we are done;
  otherwise, we let $\phi$ and $B$ be as in Lemma \ref{leaflet_morse}
  with $\Delta(I)$ in place of $\Delta$.

  Let $V_p$ be a Rolle leaf of $d^p$ for each $p$.  By Proposition
  \ref{leaflet_dec} and the inductive hypothesis, it now suffices to
  show that every component of $X:= A \cap \bigcap_{p \in I} V_p$
  intersects $B$.  However, since $d^{\Delta(I),A}_{k(\Delta(I),A)}$
  has dimension, $X$ is a closed, embedded submanifold of $A$.  Thus,
  $\phi$ attains a maximum on every component of $X$, and any point in
  $X$ where $\phi$ attains a local maximum belongs to $B$.
\end{proof}

\begin{cor}  
  \label{normal_fiber_cutting}
  Let $d$ be a definable, analytic nested distribution on $U$ and $m
  \le n$.  Then there is a finite collection $\P$ of normal leaflets
  in $U$ contained in $A$ such that for every Rolle leaf $\,V$ of
  $d$, we have
  \begin{equation*}
    \Pi_m(A \cap V) = \bigcup_{N \in \P} \Pi_m(N \cap V)
  \end{equation*}
  and for every $N \in \P$, the set $N \cap V$ is an analytic
  submanifold of $\,U$, $\Pi_m\rest{(N \cap V)}$ is an immersion and
  for every $n' \le n$ and every strictly increasing $\lambda:\{1,
  \dots, n'\} \into \{1, \dots, n\}$, the projection
  $\Pi_{\lambda}\rest{(N \cap V)}$ has constant rank.
\end{cor}

\begin{proof}
  Apply Lemma \ref{cell_cutting} with $q:= n+1$, $d^p:= \ker
  dx_p$ for $p=1, \dots, n$, $d^q:= d$ and $I:= \{1, \dots, m,n+1\}$.
\end{proof}

\begin{cor}
  \label{fiber_cutting}
  Let $M \subseteq \RR^n$ be a definable $C^2$-manifold, $d$ be a
  definable nested distribution on $M$, $B \subseteq \RR^n$ be a
  definable set and $m \le n$.  Then there is a finite collection $\P$
  of definable analytic manifolds contained in $B$ such that for every
  Rolle leaf $\,V$ of $d$, we have
  \begin{equation*}
    \Pi_m(B \cap V) = \bigcup_{N \in \P} \Pi_m(N \cap V)
  \end{equation*}
  and for $N \in \P$, the set $N \cap V$ is an analytic manifold and a
  submanifold of $M$, $\Pi_m\rest{(N \cap V)}$ is an immersion and for
  every $n' \le n$ and every strictly increasing $\lambda:\{1, \dots,
  n'\} \into \{1, \dots, n\}$, the projection $\Pi_{\lambda}\rest{(N
    \cap V)}$ has constant rank.
\end{cor}

\begin{proof}
  By analytic cell decomposition and Corollary \ref{khovprop}(1), we
  may assume that $M = B$ is an analytic cell and $d$ is analytic.
  Let $\mu:\{1, \dots, \dim M\} \into \{1, \dots, n\}$ be such that
  $\Pi_\mu(M)$ is open and $\Pi_\mu\rest{M}:M \into \Pi_\mu(M)$ is an
  analytic diffeomorphism.  Then $\Pi_\mu(M)$ is normal and $M$ is the
  graph of an analytic function $g:\Pi_\mu(M) \into \RR^{n-\dim M}$.
  The corollary now follows from Corollary \ref{normal_fiber_cutting}
  applied to the push-forward of $d$ via $\Pi_\mu\rest{M}$.
\end{proof}

\subsection*{Regular closure}  \label{regular_closure}

Let $U \subseteq \RR^n$ be a normal set and $A$ a normal subset of
$U$.  Let also $d = (d_0, \dots, d_k)$ be a definable, integrable,
analytic nested distribution on $U$ and $V$ be a Rolle leaf of $d$.
Following \cite{gab:frontier}, we study below the \textbf{closure in
  $U$ of $A \cap V$}, defined by $\cl_U(A \cap V) := U \cap \cl(A\cap
V)$, and the \textbf{frontier in $U$ of $A \cap V$}, defined by
$\fr_U(A \cap V) := U \cap \fr(A\cap V)$.

\begin{prop}  
  \label{reg_closure}
  There are normal subsets $B$ and $C$ of $\,U$ such that
  \begin{equation*}
    \cl_U(A \cap V) = B \cap V \quad\text{and}\quad \fr_U(A
    \cap V) = C \cap V.
  \end{equation*}
\end{prop}

For the proof of Proposition \ref{reg_closure}, we need the following
preliminary observations; here we consider $\Sigma_n$ as a definable
subset of $\RR^{2n}$.  For every $\sigma \in \Sigma_n$, we let
$U_\sigma$ be the set of all $x \in U$ such that $\sigma(d_k(x))$ is
the graph of a linear map $L:\RR^{n-k} \into \RR^k$ satisfying $|L| <
2$.  Then $U = \bigcup_{\sigma \in \Sigma_n} U_\sigma$ by Lemma
\ref{eta-bdd-cover}, and by definable choice there is a definable map
$x \mapsto \sigma_x:U \into \Sigma_n$ such that $x \in U_{\sigma_x}$
for all $x \in U$.  Since each $U_\sigma$ is open, there is a
definable map $x \mapsto \epsilon(x):U \into (0,\infty)$ such that $y
\in U_{\sigma_x}$ for all $y \in B(x,4(n-k)\epsilon(x))$.  We
put $$G:= \set{(x,y) \in U \times U:\ y \in B(x,2(n-k)\epsilon(x))}$$
and define $L:G \into GL(\RR^{n-k}, \RR^k)$ such that
$\sigma_x(d_k(y))$ is the graph of $L(x,y)$ for all $(x,y) \in G$.
Then $L$ is definable and for every $x \in U$, the map
$L_x:B(x,2(n-k)\epsilon(x)) \into GL(\RR^{n-k},\RR^k)$ defined by
$L_x(y):= L(x,y)$ is analytic.

Below, we set $x_- := x_{\le n-k}$ and $x_+:= x_{>n-k}$ and 
\begin{equation*}
  B^\sigma(x,\epsilon):= \sigma_x^{-1}\left(B(x_-,\epsilon) \times
    B(x_+,2(n-k)\epsilon)\right) 
\end{equation*}
for $x \in U$ and $\epsilon > 0$.  Also for $x \in U$, we denote by
$V_x$ the leaf of $d_k^{B^\sigma(x,\epsilon(x))}$ containing $x$.  By
Lemma \ref{local_single} and because $V_x$ is an analytic manifold,
the set $\sigma_x(V_x)$ is the graph of an analytic,
$2(n-k)$-Lipschitz map $F_x:B(x_-,\epsilon(x)) \into
B(x_+,2(n-k)\epsilon(x))$ such that
\begin{equation}
  \label{fiber_diffeq}
  dF_x(z) = L_x(z,F_x(z)) \quad\text{ for all } z \in
  B(x_-,\epsilon(x)). 
\end{equation}
For the next lemma, we put
\begin{equation*}
  G_0:= \set{(x,y) \in U \times \RR^n:\ |y| < 2(n-k)\epsilon(x)}.
\end{equation*}

\begin{lemma}  
  \label{4.1}
  For $\nu \in \NN$ there is a definable map $P_\nu:G_0 \into \RR^k$
  such that for every $x \in U$,
  \begin{renumerate}
  \item the map $P_{\nu,x}:B(0,2(n-k)\epsilon(x)) \into \RR^k$ defined
    by $P_{\nu,x}(y):= P_\nu(x,y)$ is a homogeneous polynomial of
    degree $\nu$ in $y$;
  \item $\sum_{\nu=0}^\infty P_{\nu,x}(y-x)$ converges to an analytic
    map $\phi_x: B^\sigma(x,\epsilon(x)) \into \RR^k$ definable in $\P(\R)$
    such that $\sigma_x(V_x) = \set{y \in B^\sigma(x,\eta_x):\ \phi_x(y) =
      0}$.
  \end{renumerate}
\end{lemma}

\begin{proof}
  Differentiating \eqref{fiber_diffeq} with respect to $z$, one finds
  by induction on $|\alpha| = \alpha_1+ \cdots + \alpha_{n-k}$ for
  $\alpha \in \NN^{n-k}$ that there is a definable function
  $L^\alpha:G \into \RR$ such that for all $x \in U$ and $z \in
  B(x_-,\epsilon(x))$,
  \begin{equation*}
    L^\alpha(x,z,F_x(z)) = \frac1{\alpha!} \frac
    {\partial^\alpha}{\partial z^\alpha} F_x(z),  
  \end{equation*}
  and such that for all $x \in U$, the function
  $L^\alpha_x:B(x,2\eta_x) \into \RR$ defined by $L^\alpha_x(y):=
  L^\alpha(x,y)$ is analytic.  For $x \in U$, we now define
  $\phi_x:B^\sigma(x,\eta_x) \into \RR^k$ by $\phi_x(y):= \sigma_x(y)_+ -
  F_x(\sigma_x(y)_-)$; then $\phi_x$ is analytic and definable in
  $\P(\R)$, and $\sigma_x(V_x) = \set{y \in B^\sigma(x,\eta_x):\ \phi_x(y)
    = 0}$.  Moreover, from Taylor expansion we get $\phi_x(y) =
  \sum_{q=0}^\infty P_{q,x}(y-x)$, where $P_{0,x}(y):= \sigma_x(x)_+ -
  F_x(\sigma_x(x)_-)$, $P_{1,x}(y):= \sigma_x(y)_+ -
  L_x(\sigma_x(x)_-) \cdot \sigma_x(y)_-$ and
  \begin{equation*}
    P_{\nu,x}(y) := \sum_{|\alpha| = \nu}
    L^\alpha_x(\sigma_x(x)_-,F_x(\sigma_x(x)_-)) 
    \cdot (y_-)^\alpha \quad\text{for } \nu > 1;
  \end{equation*}
  hence $P_\nu(x,y):= P_{\nu,x}(y)$ will do.
\end{proof}

Given an analytic map $h = (h_1, \dots, h_l):U \into \RR^l$, $\nu \in
\NN$ and $x \in U$, we denote by $h_x^\nu:U - x \into \RR^l$ the
Taylor expansion of order $\nu$ of $h$ at $x$ and by $h_{\min}:U \into
\RR$ the function $h_{\min}(x):= \min\{h_1(x), \dots, h_l(x)\}$.

\begin{proof}[Proof of Proposition \ref{reg_closure}]
  Assume that $$A = \set{x \in U:\ g(x) = 0,\ h(x) > 0},$$ with $g:U
  \into \RR^q$ and $h:U \into \RR^r$ definable and analytic, and put
  $Z(h):= \set{x \in U:\ h(x) = 0}$.  It suffices to find a normal set
  $C \subseteq Z(h)$ in $U$ such that $\fr_U(A \cap V) = C \cap
  V$, since then $B:= A \cup C$ will do.  Below we work with the
  notations from Lemma \ref{4.1} and the paragraphs preceding it.

  Let $Y:= \set{(x,t) \in U \times (0,\infty):\ x \in V \cap Z(h),\ 0
    < t < \epsilon(x)}$.  First, we define $\alpha:Y \into [0,\infty)$
  by
  \begin{equation*}
    \alpha(x,t):= \max\big(\set{h_{\min} (y):\ y \in V,\ g(y) = 0,\ \|y-x\|
      \le t} \cup \{0\}\big).
  \end{equation*}
  By Lemma \ref{pfaffian_puiseux}, there exists an $N \in \NN$ such
  that for all $x \in V \cap Z(h)$, either ultimately $\alpha_x(t) =
  0$ or ultimately $\alpha_x(t) > t^N$.

  Fix an arbitrary $x \in V \cap Z(h)$.  Then $x \in \fr(A \cap V)$ if
  and only if ultimately $\alpha_x(t) > t^N$.  However, we have
  ultimately $\alpha_x(t) > t^N$ if and only if $x$ belongs to the
  closure of $\set{y \in V:\ g(y) = 0,\ h(y) > \|y-x\|^N > 0}$, and
  the latter  holds if and only if $x$ belongs to the closure
  of $\set{y \in D_x:\ g(y) = 0}$, where
  \begin{equation*}
    D_x := \set{y \in V:\ 2 h_x^N(y-x) \ge \|y-x\|^N > 0}.
  \end{equation*}

  Second, we define $\beta:Y \into [0,\infty)$ by
  \begin{equation*}
    \beta(x,t) := \min\set{|(g(y),\phi_x(y))|:\ 2 h_x^N(y-x) \ge
        \|y-x\|^N,\ \|y-x\| = t}.
  \end{equation*}
  Again by Lemma \ref{pfaffian_puiseux}, there exists an $M \in \NN$
  such that for all $x \in V \cap Z(h)$, either ultimately $\beta_x(t)
  = 0$ or ultimately $\beta_x(t) > t^M$.

  Fix again an arbitrary $x \in V \cap Z(h)$.  Then $x$ is
  \textit{not} in the closure of $\set{y \in D_x:\ g(y) = 0}$ if and
  only if ultimately $\beta_x(t) > t^M$.  However, if ultimately
  $\beta_x(t) = 0$, then $x$ is in the closure of $$\set{y \in D_x:\ 4
    |(g(y),\phi_x(y))| < \|y-x\|^M},$$ which implies that $x$ is in the
  closure of
  \begin{equation*}
    E_x:= \set{y \in D_x:\ 2
      \left|\left(g_x^M(y-x),(\phi_x)_x^M(y-x)\right)\right| <
      \|y-x\|^M}. 
  \end{equation*}
  Conversely, if $x \in \cl E_x$, then $x$ is in the closure
  of $$\set{y \in D_x:\ |(g(y),\phi_x(y))| < \|y-x\|^M},$$ which
  implies that ultimately $\beta_x = 0$.  It follows from the above
  that
  \begin{itemize}
  \item[($\ast$)] for all $x \in V \cap Z(h)$, $x \in \fr(A \cap V)$
    if and only if $x \in \cl E_x$.
  \end{itemize}

  Let $G \in \RR[a,y]^q$ be the general $q$-tuple of polynomials in
  $y$ of degree $M$ and coefficients $a \in \RR^{m_1}$, $H \in
  \RR[b,y]^r$ be the general $r$-tuple of polynomials in $y$ of degree
  $N$ and coefficients $b \in \RR^{m_2}$ and $\Phi \in \RR[c,y]^k$ be
  the general $k$-tuple of polynomials in $y$ of degree $M$ and
  coefficients $c \in \RR^{m_3}$.  Let $S \subseteq \RR^n \times
  \RR^{m_1} \times \RR^{m_2} \times \RR^{m_3} \times \RR^n$ be the
  semialgebraic set
  \begin{multline*}
    S:= \big\{(x,a,b,c,y):\ 2 H(b,y-x) \ge \|y-x\|^N > 0,\\ 2
    |(G(a,y-x), \Phi(c,y-x))| < \|y-x\|^M\big\}.
  \end{multline*}
  Then there are definable, analytic functions $a:U \into \RR^{m_1}$,
  $b:U \into \RR^{m_2}$ and $c:U \into \RR^{m-3}$ such that for all $x
  \in U$,
  \begin{equation*}
    E_x = \set{y \in V:\ (x,a(x),b(x),c(x),y) \in S}.
  \end{equation*}
  Thus by ($\ast$), we have for $x \in V \cap Z(h)$ that $x \in \fr(A
  \cap V)$ if and only if $x \in \cl S_{(x,a(x),b(x),c(x))}$.  By
  Tarski's Theorem, there is a semialgebraic set $T \subseteq
  \RR^{n+m_1+m_2+m_3+n}$ such that for all $(x,a,b,c) \in
  \RR^{n+m_1+m_2+m_3}$, we have $T_{(x,a,b,c)} = \cl S_{(x,a,b,c)}$.
  Therefore, the set
  \begin{equation*}
    C:= \set{x \in U:\ (x,a(x), b(x), c(x), x) \in T}
  \end{equation*}
  is normal in $U$ and satisfies $\fr(A \cap V) = C \cap V$.
\end{proof}

Combining Corollary \ref{normal_fiber_cutting} with Proposition
\ref{reg_closure}, we obtain

\begin{cor}  
  \label{closure_cutting}
  Let $m \le n$.  Then there is a finite collection $\P$ of normal
  leaflets in $U$ contained in $A$ such that
  \begin{renumerate}
  \item $\Pi_m(\cl_U(A \cap V)) = \bigcup_{N \in \P} \Pi_m(N \cap V)$;
  \item $\Pi_m(A \cap V)$ and $\Pi_m(\fr_U(A \cap V))$ are unions
    of some of the $\Pi_m(N \cap V)$ with $N \in \P$;
  \item for every $N \in \P$, the set $N \cap V$ is an analytic
    submanifold of $U$, the restriction of $\,\Pi_m$ to $N \cap V$ is
    an immersion, and for every $m' \le m$ the restriction of
    $\,\Pi_{m'}$ to $N \cap V$ has constant rank. \qed
  \end{renumerate}
\end{cor}

\section{Proper nested subpfaffian sets}  \label{proper}

In this section, we put $I:= [-1,1]$ and $I':= I \setminus \{0\}$ and
let $m,n \ge 1$.  To simplify terminology in this section, we write
``pfaffian'' in place of ``pfaffian over $\R$''.

\begin{df} 
  \label{proper_subpfaffian} 
  Let $Y \subseteq \RR^n$ be closed, and assume there is a $U
  \subseteq \RR^n$ such that $U$ is normal, $I^n \setminus Y \subseteq
  U$ and $U \cap Y = \emptyset$.  Let $d = (d_0, \dots, d_k)$ be a
  definable, analytic nested distribution on $U$, $V$ a Rolle leaf of
  $d$ and $A$ a normal subset of $U$.  In this situation, we say that
  the basic nested pfaffian set $V \cap A \cap I^n$ is
  \textbf{restricted off} $Y$.  A nested pfaffian set is
  \textbf{restricted off $Y$} if it is a finite union of basic nested
  pfaffian sets that are restricted off $Y$.

  Let $Z \subseteq \RR^m$ be closed.  A nested subpfaffian set $W
  \subseteq I^m$ is \textbf{proper off} $Z$ if $W$ is a finite union
  of sets of the form $\Pi^n_m(X)$, where $X \subseteq I^n$ is
  restricted nested pfaffian off $Z \times \RR^{n-m}$.
\end{df}

\begin{expl}
  \label{proper_example}
  Let $X \subseteq I^n$ be restricted nested pfaffian off $\{0\}$.
  Then $X \setminus (\{0\} \times \RR^{n-1})$ is restricted nested
  pfaffian off $\{0\} \times \RR^{n-1}$.
\end{expl}


\begin{prop}  
  \label{proper_cells}
  Let $W_1, \dots, W_q \subseteq I^m$ be proper nested subpfaffian
  off $\{0\} \times \RR^{m-1}$.  Then there is a finite partition $\C$
  of $I' \times I^{m-1}$ into analytic cells compatible with $W_1,
  \dots, W_q$ such that each $C \in \C$ is proper nested subpfaffian
  off $\{0\} \times \RR^{m-1}$.
\end{prop}

To prove Proposition \ref{proper_cells}, we first need to establish a
few closure properties for the collection of all proper nested
subpfaffian sets off $\{0\} \times \RR^{m-1}$, for $m=1,2,\dots$.

\begin{lemma}  
  \label{5.3}
  The collection of all proper nested subpfaffian sets off $\{0\}
  \times \RR^{m-1}$, for $m=1,2,\dots$, is closed under taking finite
  unions, projections on the first $k$ coordinates for $1 \le k \le
  m$, permutations of the last $m-1$ coordinates and topological
  closure inside $I^m \setminus (\{0\} \times \RR^{m-1})$.
\end{lemma}

\begin{proof}
  Closure under taking finite unions, projections on the first $k$
  coordinates with $k \in \{1, \dots, m\}$ and permutations of the
  last $m-1$ coordinates is obvious from the definition.  Closure with
  respect to taking topological closure inside $I^m \setminus (\{0\}
  \times \RR^{m-1})$ follows from Corollary \ref{closure_cutting}.
\end{proof}

The collection of all proper nested subpfaffian sets off $\{0\}
\times \RR^{m-1}$, for $m =1,2,\dots$, is obviously not closed under
taking cartesian products.  However, we have the following weaker
statement:

\begin{lemma}
  \label{fiber_product}
  Let $W \subseteq I^m$ be proper nested subpfaffian off $\{0\}
  \times \RR^{m-1}$, and let $m' \ge 1$ and $W' \subseteq I^{m'}$ be
  proper nested subpfaffian off $\{0\} \times \RR^{m'-1}$.
  \begin{enumerate}
  \item $W \times I$ is proper nested subpfaffian off $\{0\} \times
    \RR^m$.
  \item Let $1 \le k \le \min\{m,m'\}$.  Write $(x,y)$ and $(x,y')$
    for the elements of $\RR^m$ and $\RR^{m'}$, respectively, where $x
    \in \RR^k$, $y \in \RR^{m-k}$ and $y' \in \RR^{m'-k}$.  Then the
    fiber product
    \begin{equation*}
      W \times_k W' := \set{(x,y,y') \in \RR^{m+m'-k}:\ (x,y) \in W,
        (x,y') \in W'}
    \end{equation*}
    is proper nested subpfaffian off $\{0\} \times \RR^{m+m'-k-1}$.
  \end{enumerate}
\end{lemma}

\begin{proof}[Sketch of proof]
  It suffices to consider the case where $W = \Pi_m(X)$ for some
  restricted nested pfaffian set $X \subseteq I^n$ off $\{0\} \times
  \RR^{n-1}$.  

  (1) Arguing as in Corollary \ref{pfaff_intersection}, we see that $W
  \times I = \Pi_{m+1}(Y)$, where
  \begin{equation*}
    Y:= \set{(x,t,y) \in I^{n+1}:\ x \in I^m,\ y \in I^{n-m},\ t \in I
      \text{ and } (x,y) \in X}
  \end{equation*}
  is restricted nested pfaffian off $\{0\} \times \RR^n$.

  (2) We may also assume that $W' = \Pi_{m'}(X')$ for some restricted
  nested pfaffian set $X' \subseteq I^{n'}$ off $\{0\} \times
  \RR^{n'-1}$.  Below, we let $z$ range over $I^{n-m}$ and $z'$ range
  over $I^{n'-m'}$.  Since $k \ge 1$, we have $W \times_k W' =
  \Pi_{m+m'-k}(Y)$, where
  \begin{equation*}
    Y :=  \set{(x,y,y',z,z') \in I^{n+n'-k}:\ (x,y,z) \in X \text{
        and } (x,y',z') \in X'} 
  \end{equation*}
  is restricted nested pfaffian off $\{0\} \times \RR^{n+n'-k-1}$.
\end{proof}

\begin{cor}
  \label{proper_intersection}
  Let $W, W' \subseteq I^m$ be proper nested subpfaffian off $\{0\}
  \times \RR^{m-1}$.  Then $W \cap W'$ is proper nested subpfaffian
  off $\{0\} \times \RR^{m-1}$.
\end{cor}

\begin{proof}
  $W \cap W' = W \times_m W'$.
\end{proof}

Also using Lemma \ref{fiber_product}, we obtain the following:

\begin{lemma}
  \label{proper_fiber}
  Let $W_1, \dots, W_q \subseteq I^m$ be proper nested subpfaffian
  off $\{0\} \times \RR^{m-1}$.  Then 
  \begin{enumerate}
  \item the set
    \begin{equation*}
      W':= \set{x' \in I^{m-1}:\ \exists y_1 < \cdots < y_q,\
        (x',y_p) \in W_p,\ p=1, \dots, q}
    \end{equation*}
    is nested subpfaffian off $\{0\} \times \RR^{m-2}$;
  \item for $p < q$ the set
    \begin{multline*}
      W := \{(x',y) \in I^m:\ \exists y_1 < \cdots < y_p < y <
      y_{p+1} < \cdots < y_q, \\ (x',y_l) \in W_l,\ l=1, \dots, q\}
    \end{multline*}
    is nested subpfaffian off $\{0\} \times \RR^{m-1}$;
  \item for $p \le q$ the set
    \begin{multline*}
      W:= \{(x',y) \in I^m:\ \exists y_1 < \cdots < y_p = y < \cdots
      < y_q, \\ (x',y_l) \in W_l,\ l=1, \dots, q\}
    \end{multline*}
    is nested subpfaffian off $\{0\} \times \RR^{m-1}$.
  \end{enumerate}
\end{lemma}

\begin{proof}
  We leave the details to the reader.
\end{proof}

\begin{proof}[Proof of Proposition \ref{proper_cells}]
  By induction on $m$; the case $m=1$ follows from the
  \hbox{o-minimality} of $\P(\R)$, so we assume that $m>1$ and the
  theorem holds for lower values of $m$.  Increasing $q$ if necessary,
  we may assume that the singleton set $\{0\}$ and the sets $I^{m-1}
  \times \{-1\}$ and $I^{m-1} \times \{1\}$ are among the $W_i$.
  Decomposing each $W_i$ if necessary, we may also assume that for $p
  \in \{1, \dots, q\}$ there are $n_p \ge m$, a normal set $U_p
  \subseteq \RR^{n_p} \setminus \left(\{0\} \times \RR^{n_p-1}\right)$
  containing $I^{n_p} \setminus \left(\{0\} \times I^{n_p-1}\right)$,
  a definable, analytic nested distribution $d^p = (d^p_0, \dots,
  d^p_{k(p)})$ on $U_p$, a Rolle leaf $V_p$ of $d^p$ and a normal
  subset $A_p$ of $U_p$ such that $W_p = \Pi_m(A_p \cap V_p)$.

  For $p \in \{1, \dots, q\}$, we now apply Corollary
  \ref{closure_cutting} with $n_p$, $U_p$, $d^p$, $V_p$ and $A_p$ in
  place of $n$, $U$, $d$, $V$ and $A$.  (Here we use the fact that the
  collection of all proper nested subpfaffian sets off $\{0\} \times
  \RR^{m-1}$ is closed with respect to taking topological closure
  inside $I^m \setminus (\{0\} \times \RR^{m-1})$.)  We let $\P_p$ be
  the corresponding collection of normal leaflets in $U_p$ obtained
  for $m$ and $\P'_p$ be the corresponding collection of normal
  leaflets in $U_p$ obtained with $m-1$ in place of $m$, and we put
  \begin{equation*}
    \Q := \set{\Pi_m(N \cap V_p):\ p \in \{1, \dots, q\},\ N \in
      \P_p,\ \dim(N \cap V_p) < m},
  \end{equation*}
  \begin{equation*}
    \Q' := \set{\Pi_{m-1}(N \cap V_p):\ p \in \{1, \dots, q\},\ N \in
      \P'_p,\ \dim(N \cap V_p) < m-1}.
  \end{equation*}
  By definition, the elements of $\Q$ and $\Q'$ are proper nested
  subpfaffian sets off $\{0\} \times \RR^{m-1}$ and $\{0\} \times
  \RR^{m-2}$, respectively.  Each $Z \in \Q$ is an immersed, analytic
  manifold in $\RR^m$ with empty interior such that the restriction of
  $\Pi^m_{m-1}$ to $Z$ has constant rank.  We put $\Q_0:= \set{Z \in
    \Q:\ \Pi_{m-1}\rest{Z} \text{ is an immersion}}$ and let $F$ be
  the union of all sets in $\Q_0$.  Similarly, each $Z' \in \Q'$ is an
  immersed, analytic manifold in $\RR^{m-1}$ with empty interior.

  Since every $Z \in \Q$ is definable in $\P(\R)$, there exists an $N
  \in \NN$ such that $Z_{x'}$ has at most $N$ components for every $x'
  \in \RR^{m-1}$ and every $Z \in \Q$.  For $k \le N|\Q|$ and $Z_1,
  \dots, Z_k \in \Q$, we put
  \begin{equation*}
    Z'(Z_1, \dots, Z_k):= \set{x' \in I^{m-1}: \exists y_1 < \cdots
    < y_k,\ (x',y_j) \in Z_j \text{ for  each } j},
  \end{equation*}
  and we denote by $\Q''$ the collection of these sets.  (Note in
  particular that the projections on the first $m-1$ coordinates of
  all $Z \in \Q$ belong to $\Q''$.)  By Lemma \ref{proper_fiber}, each
  set in $\Q''$ is proper nested subpfaffian off $\{0\} \times
  \RR^{m-2}$.  Hence by the inductive hypothesis applied to the
  collection $\Q' \cup \Q''$, there is a finite partition $\C'$ of $I'
  \times I^{m-2}$ into analytic cells such that $\C'$ is compatible
  with $\Q' \cup \Q''$ and each $C' \in \C'$ is proper nested
  subpfaffian off $\{0\} \times \RR^{m-2}$.

  We now fix $C' \in \C'$; it suffices to show that $C' \times I$
  admits a finite partition $\C$ into analytic cells such that $\C$ is
  compatible with $\{W_1, \dots, W_q\}$ and each $C \in \C$ is proper
  nested subpfaffian off $\{0\} \times \RR^{m-2}$.  However, for $p
  \in \{1, \dots, q\}$, the set $W_p \cap (C' \times I)$ is the union
  of some of the sets $Z \cap (C' \times I)$ with $Z \in \Q_0$ and
  some of the components of $(C' \times I) \setminus F$.  Therefore,
  it suffices to show that $C' \times I$ admits a finite partition
  $\C$ into analytic cells such that $\C$ is compatible with $\Q_0$
  and each $C \in \C$ is proper nested subpfaffian off $\{0\} \times
  \RR^{m-2}$.

  By construction, Lemma \ref{fiber_product} and Corollary
  \ref{proper_intersection}, if $Z \in \Q_0$ then the set $Z \cap (C'
  \times I)$ is proper nested subpfaffian off $\{0\} \times
  \RR^{m-1}$ and an analytic submanifold of $I' \times I^{m-1}$, and
  each of its components is the graph of an analytic function from
  $C'$ to $\RR$.  In particular, $F \cap (C' \times I)$ is a closed
  subset of $C' \times I$.  Moreover, if $Y \in \Q_0$ also, then $(Z
  \cap Y) \cap (C' \times I)$ is the union of some of the components
  of $Z \cap (C' \times I)$.  On the other hand, each component of $Z
  \cap (C' \times I)$ is of the form
  \begin{equation*}
    \big\{(x',y) \in C' \times I:\ \exists y_1 < \cdots < y_k,\ y=y_l
    \text{ and } (x',y_j) \in Z_j \text{ for each } j\big\},
  \end{equation*}
  where $k \le N|\Q|$, $l \le k$ and $Z_1, \dots, Z_k \in \Q$.  Hence
  by Lemma \ref{proper_fiber}, each such component is proper nested
  subpfaffian off $\{0\} \times \RR^{m-1}$ and an analytic cell.  It
  follows that each component of $(C' \times I) \setminus F$ is an
  open analytic cell, and each such component is proper nested
  subpfaffian off $\{0\} \times \RR^{m-1}$ by Lemma
  \ref{proper_fiber} again, because it is of the form
  \begin{equation*}
    \big\{(x',y) \in C' \times I:\ \exists y_1 < \cdots < y_k,\ y_l <
    y < y_{l+1} \text{ and } (x',y_j) \in Z_j \text{ for each }
    j\big\}
  \end{equation*}
  with $k \le N|\Q|$, $l<k$ and $Z_1, \dots, Z_k \in \Q$.
\end{proof}

\section{Rewriting nested integral manifolds} \label{main}

Assume that $\R$ admits analytic cell decomposition.  Let $M \subseteq
\RR^n$ be an analytic, definable manifold, and let $d = (d_0, \dots,
d_k)$ be an analytic, definable nested distribution on $M$.  Let also
$A \subseteq M$ be definable.  

\begin{prop}  
  \label{inside_out}
  There are $n_1, \dots, n_s \in \NN$ and, for $j = 1, \dots, s$,
  there exist an analytic, definable nested distribution $e_j$ on
  $C_j:= \set{y \in \RR^{n_j}:\ 0 < |y| < 2}$ and a definable,
  analytic embedding $\psi_j:C_j \into M$ such that, with $B_j:=
  \set{y \in \RR^{n_j}:\ 0 < |y| < 1}$,
  \begin{renumerate}
  \item $\psi_j(C_j) \subseteq A$ for each $j$ and the collection
    $\set{\psi_j(B_j): j=1, \dots, s}$ covers $A$;
  \item for every Rolle leaf $\,V$ of $d$, we have $A \cap V =
    \bigcup_{j=1}^s \psi_j(B_j \cap V_j)$, where each $V_j$ is either
    empty or a Rolle leaf of $e_j$.
  \end{renumerate}
\end{prop}

\begin{rmk}
  Each $B_j \cap V_j$ is a restricted nested pfaffian set off $\{0\}$.
\end{rmk}

\begin{proof}
  By Corollary \ref{pfaff_rolle}, we may assume that $A = M = \RR^n$.
  If $d$ has no Rolle leaves, the proposition is now trivial.  So we
  also assume that $d$ has a Rolle leaf; in particular, $d_1$ has a
  Rolle leaf $V_1$, say.  Then $V_1$ is embedded, closed and of
  codimension $1$ in $\RR^n$, so $V_1$ separates $\RR^n$.  Let $D_1$
  and $D_2$ be closed boxes in $\RR^n \setminus V_1$ with nonempty
  interior and contained in different components of $\RR^n \setminus
  V_1$, and denote by $c_1$ and $c_2$ their centers and by $U_1$ and
  $U_2$ their complements in $\RR^n$.  For $j=1,2$, let $D_j'$ be a
  closed box with center $c_j$ such that $D_j'$ is contained in the
  interior of $D_j$, and put $U_j':= \RR^n \setminus D_j'$.  Let
  $\phi_j:\RR^n\setminus\{c_j\} \into \RR^n \setminus\{0\}$ be a
  definable, analytic diffeomorphism such that $\phi_j(U'_j) = C :=
  \set{x \in \RR^n:\ 0 < |x| < 2}$ and $\phi_j(U_j) = B := \set{x \in
    \RR^n:\ 0 < |x| < 1}$.  We let $e_j$ be the push-forward of the
  restriction of $d$ to $\RR^n \setminus \{c_j\}$ via $\phi_j$ and put
  $\psi_j := \phi_j^{-1}$.  Note that $V_1 \subseteq \psi_j(B)$ for
  each $j$, and each component of $\RR^n \setminus V_1$ is contained
  in $\psi_j(B)$ for one $j$.  Therefore, any Rolle leaf of $d$ is
  contained in $\psi_j(B)$ for at least one $j$.  The corollary now
  follows.
\end{proof}

\begin{cor}  
  \label{pfaffian_complement}
  Let $V$ be a Rolle leaf of $d$.  Then $A \cap V$ is a finite union
  of simply connected nested subpfaffian sets over $\R$ that are
  analytic manifolds.
\end{cor}

\begin{proof}
  By induction on $\delta := \dim A$; if $\delta = 0$, there is
  nothing to do, so we assume $\delta > 0$ and the corollary holds for
  lower values of $\delta$.  By the previous proposition and the
  inductive hypothesis, we may assume that $M = \{x \in \RR^\delta:\ 0
  < |x| < 2\}$ and $A = \{x \in \RR^\delta:\ 0 < |x| < 1\}$, and we
  let $V$ be a Rolle leaf of $d$.  Since $A \cap V$ is restricted
  nested pfaffian off $\{0\}$, the corollary now follows with $(A \cap
  V) \setminus \big(\{0\} \times \RR^{\delta-1}\big)$ in place of $A
  \cap V$ from Example \ref{proper_example} and Proposition
  \ref{proper_cells}.  Since $A \cap \big(\{0\} \times
  \RR^{\delta-1}\big)$ has dimension $\delta-1$, the corollary with $A
  \cap V \cap \big(\{0\} \times \RR^{\delta-1}\big)$ in place of $A
  \cap V$ follows from the inductive hypothesis.
\end{proof}

For Lemma \ref{Rolle_to_subpfaffian} and Proposition \ref{Rolle}
below, we let $\S$ be an o-minimal expansion of $\R$.

\begin{lemma}
  \label{Rolle_to_subpfaffian}
  Let $(W_0, \dots, W_k)$ be a nested integral manifold of $d$
  definable in $\S$ and $\nu \le n$, and assume that $d$ is analytic
  and integrable and $(W_0, \dots, W_{k-1})$ is a nested Rolle leaf
  over $\R$.  Then there is a nested subpfaffian set $Y \subseteq
  \RR^n$ such that $W_k \subseteq Y$ and $\dim \Pi_\nu(Y) \le \dim
  \Pi_\nu(W_k)$.
\end{lemma}

\begin{proof}
  By induction on $\dim W_k$; if $\dim W_k = 0$, there is nothing to
  do, so we assume $\dim W_k > 0$ and the lemma holds for lower values
  of $\dim W_k$.  By Corollary \ref{pfaffian_complement}, there are
  simply connected nested subpfaffian sets $Y_1, \dots, Y_q \subseteq
  \RR^n$ such that each $Y_p$ is an analytic submanifold of $\RR^n$
  and $W_{k-1} = Y_1 \cup \cdots \cup Y_q$.  By Corollary
  \ref{fiber_cutting}, for each $p$, there are $n_p \ge n$ and nested
  pfaffian sets $Z_{p,1}, \dots, Z_{p,k(p)} \subseteq \RR^{n_p}$ over
  $\R$ such that $Y_p = \Pi_n(Z_{p,1}) \cup \cdots \cup
  \Pi_n(Z_{p,k(p)})$ and each $Z_{p,j}$ is an analytic manifold,
  $\Pi_n\rest{Z_{p,j}}$ is an immersion and $\Pi_\nu\rest{Z_{p,j}}$
  has constant rank.  Refining the $Z_{p,j}$ further if necessary and
  using Corollary \ref{khovprop}(1), we may assume in addition that
  each $Z_{p,j}$ is a Rolle leaf over $\R$ obtained from a definable,
  analytic, integrable nested distribution $d_{p,j} = (d_{p,j,0},
  \dots, d_{p,j,k(p,j)})$ on a definable analytic manifold $M_{p,j}
  \subseteq \RR^{n_p}$ such that $\Pi_n(M_{p,j}) \subseteq M$, and
  that the dimension of the spaces $F_{p,j}(y) :=
  \Pi_n(d_{p,j,k(p,j)}(y))$ and $E_{p,j}(y) := F_{p,j}(y) \cap
  d_k(\Pi_n(y))$ and of their projections $\Pi_\nu(F_{p,j}(y))$ and
  $\Pi_\nu(E_{p,j}(y))$ are constant as $y$ ranges over $M_{p,j}$.
  Note in particular that $F_{p,j}(y) \subseteq d_{k-1}(\Pi_n(y))$ and
  $\dim F_{p,j}(y) = \dim d_{p,j,k(p,j)}(y)$ for $y \in M_{p,j}$; we
  denote by $\Lambda$ the set of all pairs $(p,j)$ such that
  $F_{p,j}(y) \ne E_{p,j}(y)$ for $y \in M_{p,j}$.

  For $(p,j) \in \Lambda$ and $y \in M_{p,j}$, we now define
  $d_{p,j,k(p,j)+1}(y):= d_{p,j,k(p,j)}(y) \cap
  \Pi_n^{-1}(E_{p,j}(y))$.  Then $d_{p,j}':= (d_{p,j,0}, \dots,
  d_{p,j,k(p,j)+1})$ is a definable nested distribution on $M_{p,j}$,
  and the set $Z'_{p,j} := \Pi_n^{-1}(W_k) \cap Z_{p,j}$ is an
  integral manifold of $d_{p,j}'$ definable in $\S$.  Similarly, for
  $(p,j) \notin \Lambda$, the set $Z'_{p,j}:= \Pi_n^{-1}(W_k) \cap
  Z_{p,j}$ is an integral manifold of $d'_{p,j} := d_{p,j}$ definable
  in $\S$.

  We now fix $(p,j)$; it suffices to prove the lemma with $Z'_{p,j}$
  in place of $W_k$.  If $(p,j) \notin \Lambda$, we conclude with $Y:=
  Z_{p,j}$; so we assume from now on that $(p,j) \in \Lambda$.  If
  $F_{p,j}(y) \ne d_{k-1}(\Pi_n(y))$ for $y \in M_{p,j}$, then $\dim
  Z'_{p,j} < \dim W_k$ and we conclude using the inductive hypothesis;
  so we assume also from now on that $F_{p,j}(y) = d_{k-1}(\Pi_n(y))$
  for $y \in M_{p,j}$.  Let $L$ be a leaf of $d'_{p,j}$ such that $L
  \cap Z'_{p,j} \ne \emptyset$; then $L \subseteq Z_{p,j}$.  We claim
  that $L$ is a Rolle leaf of $d'_{p,j}$; this claim finishes the
  proof of the lemma, because $\dim \Pi_\nu(L) \le \dim \Pi_\nu(W_k)$
  by construction.  To see the claim, note that $Y_p$ is an open
  subset of $W_{k-1}$ and $\Pi_n(L)$ is a connected integral manifold
  of $d_k^{Y_p}$.  By Haefliger's Theorem, the leaf $L'$ of
  $d_k^{Y_p}$ containing $\Pi_n(L)$ is a Rolle leaf.  On the other
  hand, since $Z_{p,j}$ is a Rolle leaf and
  $\Pi_n\rest{d_{p,j,k(p,j)}(y)}$ is an immersion for $y \in M_{p,j}$,
  the set $Z_{p,j}$ is the graph of a function $g_{p,j}:\Pi_n(Z_{p,j})
  \into \RR^{n_p-n}$.  Thus, if $\gamma:[0,1] \into Z_{p,j}$ is a
  curve intersecting $L$ in two distinct points, then $\Pi_n \circ
  \gamma:[0,1] \into Y_p$ is a curve intersecting $L'$ in two distinct
  points, so $\Pi_n \circ \gamma$ is tangent to $d_k$ at some point.
  It follows that $\gamma$ is tangent to $d_{p,j,k(p,j)+1}$ at some
  point, as required.
\end{proof}

\begin{prop}  
  \label{Rolle}
  Let $W \subseteq \RR^n$ be an integral manifold over $\R$ definable
  in $\S$ and $\nu \le n$.  Then there is a nested subpfaffian set $Y
  \subseteq \RR^n$ such that $W \subseteq Y$ and $\dim \Pi_\nu(Y) \le
  \dim \Pi_\nu(W)$.
\end{prop}

\begin{proof}
  We proceed by induction on $k$.  If $k=0$ there is nothing to do, so
  we assume $k>0$ and the proposition holds for lower values of $k$.
  By the inductive hypothesis and Corollary \ref{khovprop}(1), there
  exists a $q \in \NN$, and for $p=1, \dots, q$ there exist $n_p \ge
  n$, a definable manifold $N_p \subseteq \RR^{n_p}$, a definable
  nested distribution $f_p = (f_{p,0}, \dots, f_{p,k(p)})$ on $N_p$
  and a nested Rolle leaf $V_p = (V_{p,0}, \dots, V_{p,k(p)})$ of
  $f_p$ such that $ W_{k-1} \subseteq \bigcup_{p=1}^q
  \Pi_n(V_{p,k(p)})$ and $\dim \Pi_\nu(V_{p,k(p)}) \le \dim
  \Pi_\nu(W_{k-1})$ for each $p$.  After refining and pruning the
  collection $\big\{N_{1}, \dots, N_{q}\big\}$ if necessary, we may
  assume that that
  \begin{itemize}
  \item[($\dagger$)] $\Pi_{n}(N_{p}) \subseteq M$ and
    $\Pi_n(V_{p,k(p)}) \cap W_k \ne \emptyset$ for each $p$, and that
    $W_k \subseteq \Pi_n(V_{1,k(1)}) \cup \cdots \cup
    \Pi_n(V_{q,k(q)})$.
  \end{itemize}
  We now finish the proof of the proposition by induction on $$\delta
  := \max\set{\dim \Pi_\nu \left(V_{p,k(p)}\right):\ p = 1, \dots,
    q},$$ simultaneously for all such collections $\set{(N_p, f_p,
    Z_p)}$ satisfying ($\dagger$).  If $\delta = \dim \Pi_\nu(W_k)$,
  we are now done.  So assume $\delta > \dim \Pi_\nu(W_k)$ and the
  proposition holds for lower values of $\delta$.  Refining and
  pruning the collection $\big\{N_{1}, \dots, N_{q}\big\}$ again if
  necessary, we may assume for each $p$, each $\mu \le \nu$ and each
  strictly increasing $\lambda:\{1, \dots, \mu\} \into \{1, \dots,
  \nu\}$ that the dimension of the spaces $$\Pi_\lambda(f_{p,k(p)}(y))
  \quad\text{and}\quad \Pi_\lambda\left(\Pi_\nu(f_{p,k(p)}(y)) \cap
    \Pi_\nu( d_k(\Pi_n(y)))\right)$$ is constant as $y$ ranges over
  $N_{p}$.  Then by Proposition \ref{lifting} and Lemma
  \ref{Rolle_to_subpfaffian}, for each $p$ such that $\dim \Pi_\nu
  \left(V_{p,k(p)}\right) = \delta$, the Rolle leaf $V_{p,k(p)}$ over
  $\R$ can be replaced by finitely many Rolle leaves $V'$ over $\R$
  such that $\dim \Pi_\nu(V') < \delta$.  The proposition now follows
  from the inductive hypothesis.
\end{proof}

\begin{proof}[Proofs of Proposition 1 and Theorem 1]
  Let $W \subseteq \RR^n$ be a Rolle leaf over $\R$; we proceed by
  induction on $\dim W$.  If $\dim W = 0$, there is nothing to do, so
  we assume $\dim W > 0$ and Theorem 1 holds for lower values of $\dim
  W$.  Using the definable, analytic diffeomorphism $t \mapsto
  t/\sqrt{1+t^2}: \RR \into (-1,1)$ in each coordinate, we may assume
  that $M$ is bounded.  By analytic cell decomposition and Remark
  \ref{limit_rmks}, we may assume that $M$ is an analytic cell.  By
  Proposition \ref{frontier_reduction}, there are pfaffian limits
  $K_1, \dots, K_q \subseteq \RR^n$ over $\R$ such that $\fr W
  \subseteq K_1 \cup \cdots \cup K_q$ and $\dim K_p < \dim W$ for each
  $p$.  By Proposition \ref{recover_1}, we may assume for each $p$
  that $K_p = \Pi_n(U_p)$ for some integral manifold $U_p \subseteq
  \RR^{n_p}$ over $\R$ definable in $\P(\R)$ such that $n_p \ge n$
  (this proves Proposition 1).  By Proposition \ref{Rolle} with $U_p$
  in place of $W$, we may assume that each $K_p$ is nested
  subpfaffian over $\R$, say $K_p = \Pi_n(Y_p)$ for some nested
  pfaffian set $Y_p \subseteq \RR^{m_p}$ over $\R$ such that $m_p \ge
  n_p$.  By Corollaries \ref{khovprop}(1) and \ref{fiber_cutting}, we
  may assume that each $Y_p$ is a Rolle leaf over $\R$ such that
  $\Pi_n\rest{Y_p}$ is an immersion; in particular, $\dim Y_p < \dim
  W$.  Hence by the inductive hypothesis, there is a closed nested
  subpfaffian set $Z_p \subseteq \RR^{m_p}$, for each $p$, such that
  $\fr Y_p \subseteq Z_p$ and $\dim Z_p < \dim W$.  Then the union of
  the sets $\Pi_n(Y_p \cup Z_p)$, for $p=1, \dots, q$, is a closed
  nested subpfaffian set containing $\fr W$ and of dimension strictly
  less than $\dim W$, as required.
\end{proof}

\section{Conclusion}  \label{conclusion}

We conclude by proving the corollary in the introduction; we continue
to assume that $\R$ admits analytic cell decomposition.  Let $L
\subseteq \RR^n$ be one of the Rolle leaves added to $\R$ in the
construction of $\P(\R)$ in \cite{spe:pfaffian}; it suffices to
establish the following:

\begin{prop}
  \label{interdefinable}
  $L$ is a nested subpfaffian set over $\R$.
\end{prop}

\begin{proof}
  By construction of $\P(\R)$ and Example 1, there are an $l \in \NN$
  and an $(n-1)$-distribution $e$ on $\RR^n$ definable in $\R_l$ such
  that $L$ is a Rolle leaf of $e$.  We proceed by induction on $l$; if
  $l=0$, we are done, so we assume that $l>0$ and the proposition
  holds for lower values of $l$; in particular, every set definable in
  $\R_{l}$ is definable in $\N(\R)$.  Thus, by the Main Theorem,
  analytic cell decomposition and Corollary \ref{khovprop}(1), we may
  assume that there are $n' \ge n+n^2$, a definable, analytic manifold
  $M \subseteq \RR^{n'}$, a definable, analytic nested distribution
  $d= (d_0, \dots, d_k)$ on $M$ and a Rolle leaf $V$ of $d$ such that
  $\gr e = \Pi_{n+n^2}(V)$.  By Corollary \ref{fiber_cutting}, we may
  further assume that $\Pi_{n+n^2}\rest{V}$ is an immersion and $\dim
  \Pi_{n+n^2}(d_k(y)) = \dim d_k(y) = n$ for all $y \in M$.  Since
  $\Pi_{n+n^2}(V) \subseteq \RR^n \times G^{n-1}_n$, we may also
  assume that $\Pi_{n+n^2}(M) \subseteq \RR^n \times G^{n-1}_n$.  Let
  now $d_{k+1}$ be the $(n-1)$-distribution on $M$ defined by
  \begin{equation*}
    d_{k+1}(y) := d_k(y) \cap
    \left(\Pi^{n'}_n\right)^{-1}(\pi(\Pi_{n+n^2}(y))), 
  \end{equation*}
  where $\pi:\RR^{n+n^2} \into \RR^{n^2}$ is the projection on the
  last $n^2$ coordinates and $\pi(\Pi_{n+n^2}(y))$ is identified with
  the $(n-1)$-dimensional subspace of $\RR^n$ that it represents.
  Then $d' := (d_0, \dots, d_{k+1})$ is a definable nested
  distribution on $M$ and $V':= V \cap (\Pi^{n'}_{n})^{-1}(L)$ is
  a Rolle leaf of $d_{k+1}$, as required.
\end{proof}


\end{document}